\documentclass[12pt,twoside,reqno,openany]{amsart}

\pdfoutput=1

\usepackage{amsbsy,amscd,amsfonts,amsmath,amssymb,amsthm,color,
fancybox,fancyhdr,footmisc,graphics,graphicx,hyperref,ifthen,mathrsfs,
multicol,multirow,nccrules,pdfpages,pifont,rotating,shuffle,stmaryrd,
textcomp,times,wasysym,wrapfig}

\usepackage[dvipsnames,svgnames,x11names]{xcolor}

\usepackage[all]{xy}
\usepackage[utf8]{inputenc}
\usepackage[T1]{fontenc}
\usepackage{mathtools}
\newtagform{EngelLie}[\scriptstyle]{$}{$}
\makeatletter\newcommand{\leqnomode}{\tagsleft@true}
\newcommand{\reqnomode}{\tagsleft@false}\makeatother

\newtheorem{Theorem}[equation]{Theorem}

\newtheorem{Proposition}[equation]{Proposition}

\newtheorem{Lemma}[equation]{Lemma}

\newtheorem{Assertion}[equation]{Assertion}

\theoremstyle{definition}

\newtheorem{Hypothesis}[equation]{Hypothesis}

\newtheorem{Definition-Notation}[equation]{D\'efinition-Notation}

\newtheorem{Principle}[equation]{Principle}

\newcommand{\N}{\mathbb{N}}

\newcommand{\R}{\mathbb{R}}

\definecolor{blue}{cmyk}{1.,1.,0.,0.63}
\definecolor{red}{cmyk}{0.,1.,1.,0.63}
\definecolor{green}{cmyk}{1.,0.,1.,0.63}
\definecolor{black}{cmyk}{1.,1.,1.,1.}

\newcommand{\green}{\textcolor{green}}
\newcommand{\red}{\textcolor{red}}

\makeatletter
\renewcommand{\@fnsymbol}[1]
{\ensuremath{\ifcase#1\or $*$\or $**$\or $***$\or $****$\or $*****$
\else\@ctrerr\fi}}
\makeatother

\newcommand{\HEAD}[2]{%
\pagestyle{fancy}
\fancyhead[RO]{\tiny\sf\thepage}
\fancyhead[CO]{{\tiny\sf #1}}
\fancyhead[LE]{\tiny\sf\thepage}
\fancyhead[CE]{{\tiny\sf #2}}
\fancyfoot{}}

\numberwithin{equation}{section}

\newcommand{\Section}[1]{
\renewcommand{\thesection}{\bf\arabic{section}}
\section{#1}
\renewcommand{\thesection}{\arabic{section}}}

\newcommand{\Subsection}[1]{
\refstepcounter{equation}
\medskip\noindent{\bf\arabic{section}.\arabic{equation}.~#1.}}

\newcommand{\SectionHead}[2]{
\Section{\bf #1}
\label{#2}
\HEAD{\ref{#2}.~{\sf 
#1}}{
Julien {\sc Heyd} and 
Jo\"el {\sc Merker}, 
D\'epartement de Math\'ematiques d'Orsay, 
Universit\'e Paris-Saclay, France}}

\newcommand{\style}[1]{{\sf #1}}

\newcommand{\Aff}{\style{Aff}}

\newcommand{\eqFG}{\style{eqFG}}

\newcommand{\eqLF}{\style{eqLF}}

\newcommand{\GL}{\style{GL}}

\newcommand{\Hessian}{\style{Hessian}}

\newcommand{\Jac}{\style{Jac}}

\renewcommand{\lim}{\style{lim}}

\newcommand{\order}{\style{order}}

\newcommand{\rank}{\style{rank}}

\newcommand{\Span}{\style{Span}}

\newcommand{\Hall}{\Hall}

\newcommand{\smallbullet}{{\scriptscriptstyle{\bullet}}}

\newcommand{\vf}{\vfill\end{document}}

\makeatletter
\def\hlinethick#1{%
\noalign{\ifnum0=`}\fi\hrule \@height #1 \futurelet
\reserved@a\@xhline}
\makeatother

\newlength{\myskip}
\begin{document}

\setcounter{section}{0}

\bigskip\bigskip


\begin{center}

{\large\bf Classification of Affinely Homogeneous} 

\medskip

{\large\bf 
Hessian Rank 2 Hypersurfaces 
$S^3 \subset \R^4$}
\label{S3-R4-Hessian-rank-2}

\bigskip\bigskip

Julien {\sc Heyd}\footnotemark[1] 
and 
Jo\"el~{\sc Merker}\footnotemark[1]

\end{center}\bigskip

\footnotetext[1]{\,\,
D\'epartement de Math\'ematiques d'Orsay,
CNRS, Universit\'e Paris-Saclay, 91405 Orsay Cedex,
France, 
{\bf julien.heyd@universite-paris-saclay.fr},
{\bf joel.merker@universite-paris-saclay.fr}}

\footnotetext[2]{\,
This research was supported
in part by the Polish National Science Centre (NCN) 
via the grant number 2018/29/B/ST1/02583,
and by the Norwegian Financial Mechanism
2014--2021 via the project registration number 2019/34/H/ST1/00636.}

\footnotetext[3]{\,
2020 Mathematics Subject Classification: 
22E05, 53A55, 20G05, 14R05, 22E45, 20-08, 53B25, 14R20, 53A15, 13J05, 
13A50, 32M17, 32A05. 
Keywords: 
Equivalence Method. 
Power Series. 
Linear Representations.
Affine Homogeneity. 
Group Reduction and Classification.
Absolute and Relative Invariants.}

\begin{center}
\begin{minipage}[t]{12.5cm}
\parindent 0.53cm
\footnotesize
\noindent
{\sc Abstract}.
We determine all affinely homogeneous 
hypersurfaces $S^3 \subset \R^4$ whose Hessian is (invariantly)
of constant rank 2,
{\em including the simply transitive ones}.

We find 34 inequivalent terminal branches yielding each to a nonempty
moduli space of homogeneous models of 
hypersurfaces $S^3 \subset \R^4$, 
sometimes parametrized by a certain complicated algebraic variety,
especially for the 15 (over 34) families of models which are 
{\em simply transitive}.

We employ the {\sl power series method of equivalence},
which captures invariants at the origin,
creates branches, and infinitesimalizes calculations.

In Lie's original classification spirit, we describe 
the found homogeneous models by listing explicit
Lie algebras of infinitesimal transformations,
sometimes parametrized by absolute invariants
satisfying certain algebraic equations.

\end{minipage}
\end{center}

\SectionHead{Introduction}
{introduction-S2-R4}

The goal of this article is to determine all affinely homogeneous
Hessian rank 2
hypersurfaces $S^3 \subset \R^4$,
{\em including the simply transitive ones}.
This problem has been studied, among
other things, 
by Mozhei~{\cite{Mozhei-2000}} who found
the {\em multiply transitive models},
and also by Wermann~{\cite{Wermann-2001}} for the special affine group.

About affinely homogeneous submanifolds, 
the literature is extensive, 
a few references being~{\cite{
Doubrov-Komrakov-Rabinovich-1996,
Doubrov-Komrakov-1998,
Eastwood-Ezhov-1999,
Wermann-2001,
Eastwood-Ezhov-2001-2,
Chen-Merker-2019,
Chen-Merker-2020,
Merker-2022,
Heyd-Merker-2024}}.
Information about history, context, perspective, problems, 
can be found there\,\,---\,\,{\em see} also these
articles' bibliographies.

To be specific, in $\R^4 \ni (x,y,z, u)$, a local analytic
hypersurface $S^3$ can be graphed, 
after an affine transformation, as:
\[
\aligned
u
\,=\,
F\big(x,y,z\big)
\,=\,
\sum_{i+j+k\geqslant0}\,
F_{i,j,k}\,
x^iy^jz^k
\,=\,
\sum_{\mu\geqslant0}\,
\sum_{i+j+k=\mu}\,
F_{i,j,k}\,
x^iy^jz^k
\,=\,
F_0
+
F_1
+
F_2
+\cdots,
\endaligned
\]
with $F$ real-analytic in some
neighborhood of the origin $(0,0,0) \in \R^3$.
The Hessian matrix:
\[
\Hessian_F
\,:=\,
\left[\!
\begin{array}{ccc}
F_{xx} & F_{xy} & F_{xz}
\\
F_{yx} & F_{yy} & F_{yz}
\\
F_{zx} & F_{zy} & F_{zz}
\end{array}
\!\right],
\]
has invariant rank~{\cite[Sec.~2]{Merker-2022}}.

\begin{Hypothesis}
\label{Hyp-Hessian-rank-2}
At every point $(x,y,z)$, assume 
$\rank\, \Hessian_F = 2$.
\end{Hypothesis}

After an elementary affine transformation
(Section~{\ref{hypersurfaces-S3-R4}}), 
order 1 terms in $F$ disappear, and
two inequivalent normalizations for the quadratic terms $F_2$
exist at order 2:
\[
\begin{array}{cc}
{} & F_2 
\\
\green{\bf 2a} & x\,y 
\\
\green{\bf 2b} & x^2+y^2
\end{array}
\]
with $x^2 - y^2$ equivalent to the more convenient $x\,y$.
At order 2, 
these are the 2 existing branches,
labelled {\green{\bf 2a}}, {\green{\bf 2b}},
similarly as in~{\cite{Heyd-Merker-2024}}
where 7 branches exist at order 2.

Beyond order 2, further branches exist, exhibited
just before stating the main result.

\begin{center}
\includegraphics[scale=0.21]{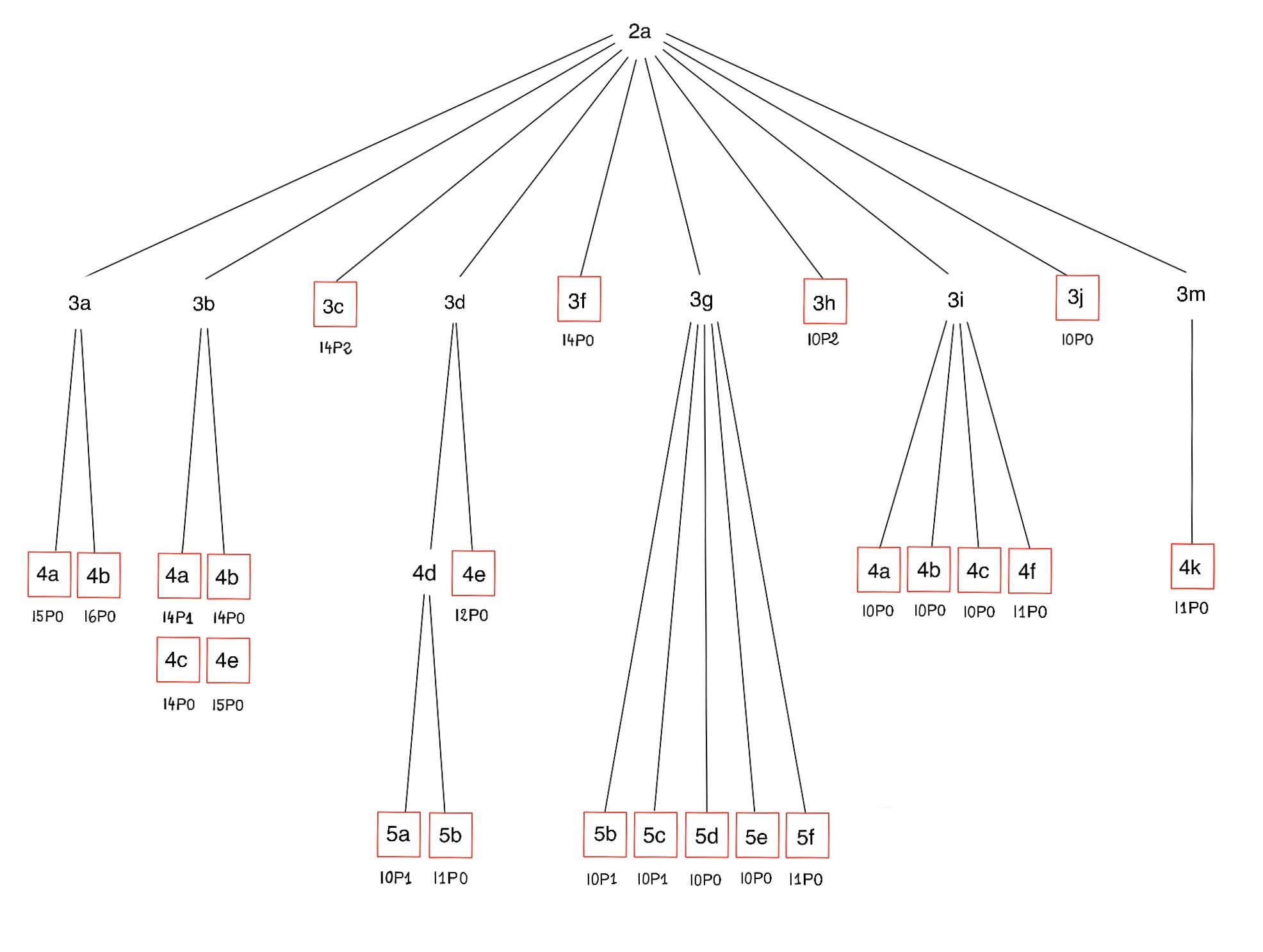}
\end{center}

\begin{center}
\includegraphics[scale=0.25]{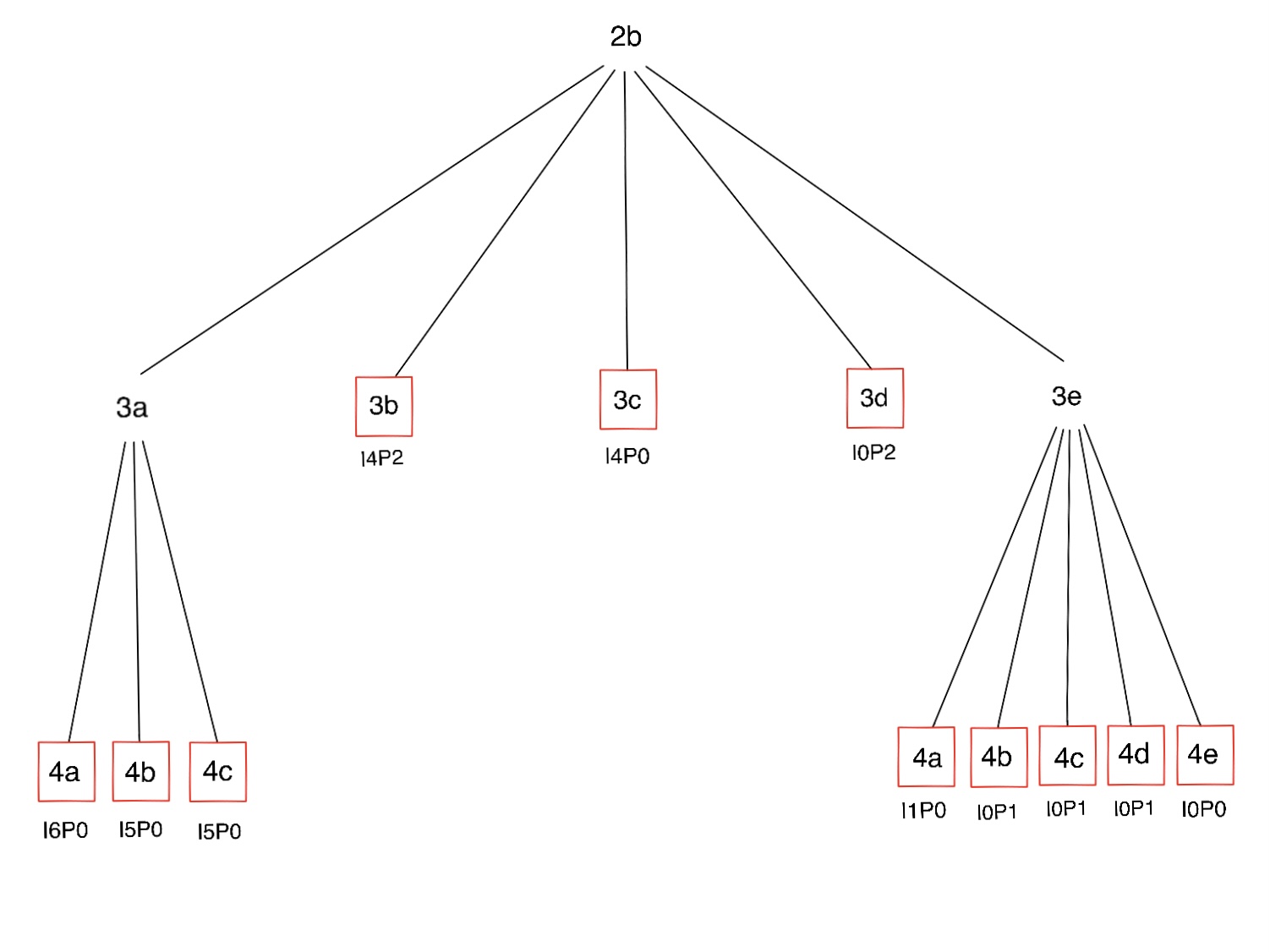}
\end{center}

\begin{Theorem}
Up to the action of a finite subgroup of $\Aff(\R^4)$, all the
possible affinely homogeneous hypersurfaces $S^3 \subset \R^4$ having
constant Hessian rank 2 are represented by means of 3 lists.

\smallskip\noindent$\bullet$\,
The list of branching trees appearing above.

\smallskip\noindent$\bullet$\,
The list of occurring linear representations appearing
in Section~{\ref{linear-representations-branches}}.

\smallskip\noindent$\bullet$\,
The list of (truncated) normal forms together with
their associated
transitive Lie algebras of vector fields,
appearing in 
Sections~{\ref{2a-models}} and~{\ref{2b-models}}.

\end{Theorem}

Our principal technique is the
{\sl Power Series Method of Equivalence},
inspired from \'Elie Cartan's method. 
Part~I of~{\cite{Heyd-Merker-2024}} describes in details
the features of this (new, computationally economical) method.

Problem~15.2 in~{\cite{Heyd-Merker-2024}} emphasizes
that there is currently no satisfactory theory
which provides necessary and sufficient conditions
on a given transitive Lie algebra of vector fields 
that are symmetries of a given type of geometric structures,
in order that, after a change of coordinates belonging
to the initial group, the geometric
structure becomes neatly expressed by means of either 
polynomials or\big/and usual transcendental functions, like
exponentials, logarithms, trigonometric functions.

This is why, similarly as in~{\cite{Heyd-Merker-2024}},
we abandon in this paper the search for closed forms
of the found homogeneous models.

\SectionHead{Hypersurfaces $S^3 \subset \R^4$ 
Under $\Aff(4,\R)$}
{hypersurfaces-S3-R4}

So, let us start.
After translation, an affine transformation of $\R^4$
fixes the origin. Consider therefore an invertible linear map
$(x,y,z,u) \longmapsto (p,q,r,v)$:
\[
\aligned
p
&
\,:=\,
a_{1,1}\,x+a_{1,2}\,y+a_{1,3}\,z+b_1\,u,
\\
q
&
\,:=\,
a_{2,1}\,x+a_{2,2}\,y+a_{2,3}\,z+b_2\,u,
\\
r
&
\,:=\,
a_{3,1}\,x+a_{3,2}\,y+a_{3,3}\,z+b_3\,u,
\\
v
&
\,:=\,
c_1\,x+c_2\,y+c_3\,z+d\,u,
\endaligned
\ \ \ \ \ \ \ \ \ \ \ \ \ \ \ \ \ \ \ \
\text{with}
\ \ \ \ \ \ \ \ \ \ \ \ \ \ \ \ \ \ \ \
0
\,\neq\,
\left\vert\!
\begin{array}{cccc}
a_{1,1} & a_{1,2} & a_{1,3} & b_1
\\
a_{2,1} & a_{2,2} & a_{2,3} & b_2
\\
a_{3,1} & a_{3,2} & a_{3,3} & b_3
\\
c_1 & c_2 & c_3 & d
\end{array}
\!\right\vert.
\]

Also, consider two local analytic hypersurfaces passing through 
the origin, graphed as:
\[
u
\,=\,
F(x,y,z)
\ \ \ \ \ \ \ \ 
{\scriptstyle{(F(0,0,0)\,=\,0)}}
\ \ \ \ \ \ \ \ \ \ \ \
\text{and}
\ \ \ \ \ \ \ \ \ \ \ \
v
\,=\,
G(p,q,r)
\ \ \ \ \ \ \ \ 
{\scriptstyle{(0\,=\,G(0,0,0))}},
\]
with convergent series:
\[
F
\,=\,
\sum_{i+j+k\geqslant 1}\,
F_{i,j,k}\,
x^i\,y^j\,z^k,
\ \ \ \ \ \ \ \ \ \ \ \ \ \ \ \ \ \ \ \
\text{and}
\ \ \ \ \ \ \ \ \ \ \ \ \ \ \ \ \ \ \ \
G
\,=\,
\sum_{i+j+k\geqslant 1}\,
G_{i,j,k}\,
p^i\,q^j\,r^k.
\]

The concerned general linear transformation of $\GL(4,\R)$ 
maps the left hypersurface $\{u = F\}$ 
to the right hypersurface $\{v = G\}$ 
if and only if a {\sl fundamental equation:}
\leqnomode\usetagform{default}
\begin{align}
\label{eqG-S3-R4}
0
\,\equiv\,
\eqFG(x,y,z),
\end{align}
holds identically in $\R\{x,y,z\}$, where:
\[
\aligned
\eqFG
&
\,:=\,
-\,c_1\,x-c_2\,y-c_2\,z-d\,F(x,y,z)
+
G
\Big(
a_{1,1}x+a_{1,2}y+a_{1,3}z+b_1F(x,y,z),
\\
&
\ \ \ \ \ \ \ \ \ \ \ \ \ \ \ \ \ \ \ \ \ \ \ \ \ \ \ \ \ \ \ \ \ \ \
\ \ \ \ \ \ \ \ \ \ \ \ \ \ \ \ \ \ \ \ \ \ \ \ \ \ \ \ \ \ \ \ \ \ \ 
\ \ \ \ 
a_{2,1}x+a_{2,2}y+a_{2,3}z+b_2F(x,y,z),
\\
&
\ \ \ \ \ \ \ \ \ \ \ \ \ \ \ \ \ \ \ \ \ \ \ \ \ \ \ \ \ \ \ \ \ \ \
\ \ \ \ \ \ \ \ \ \ \ \ \ \ \ \ \ \ \ \ \ \ \ \ \ \ \ \ \ \ \ \ \ \ \ 
\ \ \ \ 
a_{3,1}x+a_{3,2}y+a_{3,3}z+b_3F(x,y,z)
\Big),
\endaligned
\]
so that:
\[
0
\,=\,
\eqFG
\,=\,
\sum_{i,j\in\N}\,
\mathcal{N}_{i,j,k}
\Big(
a_{\smallbullet,\smallbullet},\,
b_{\smallbullet},\,
c_{\smallbullet},\,
d,\,\,
F_{\smallbullet,\smallbullet,\smallbullet},\,
G_{\smallbullet,\smallbullet,\smallbullet}
\Big)\,
x^i\,y^j\,z^k.
\]

For all $i,j,k \in \N$, the coefficient of $x^i y^j z^k$ in
$\eqFG$ can, in a standard way, be denoted as:
\[
\aligned
\big[x^i\,y^j\,z^k\big]
\eqFG
\,:=\,
\mathcal{N}_{i,j,k}
\,=\,
0,
\endaligned
\]
and we will indicate the corresponding indices 
$\green{\bf i}$, $\green{\bf j}$, $\green{\bf k}$ over the equal 
sign:
\[
0
{\overset{\green{\bf i,j,k}}{\,\,=\,\,}}
\mathcal{N}_{i,j,k}.
\]
We will proceed by increasing:
\[
\order
\,:=\,
i+j+k.
\]

Two simple affine transformations of $\R_{x,y,z,u}^4$
and of $\R_{p,q,r,v}^4$ erase $\ast\, x^1 + \ast\, y^1 + \ast\, z^1$
and $\ast\, p^1 + \ast\, q^1 + \ast\, r^1$:
\[
u
\,=\,
\red{\bf 0}
+
{\rm O}_{x,y,z}(2),
\ \ \ \ \ \ \ \ \ \ \ \ \ \ \ \ \ \ \ \
\text{and}
\ \ \ \ \ \ \ \ \ \ \ \ \ \ \ \ \ \ \ \
\aligned
v
\,=\,
\red{\bf 0}
+
{\rm O}_{p,q,r}(2),
\endaligned
\]
so that order 1 terms are {\sl normalized} to be $\red{\bf 0}$.

\begin{Lemma}
\label{Lm-stab-order-1}
Stabilization of order 1 terms holds if and only if 
$0 = c_1 = c_2 = c_3$:
\[
\left[
\begin{array}{cccc}
a_{1,1} & a_{1,2} & a_{1,3} & b_1
\\
a_{2,1} & a_{2,2} & a_{2,3} & b_2
\\
a_{3,1} & a_{3,2} & a_{3,3} & b_3
\\
c_1 & c_2 & c_3 & d
\end{array}
\right]^{\green{\bf 0}}
\,\,\,\leadsto\,\,\,
\left[
\begin{array}{cccc}
a_{1,1} & a_{1,2} & a_{1,3} & b_1
\\
a_{2,1} & a_{2,2} & a_{2,3} & b_2
\\
a_{3,1} & a_{3,2} & a_{3,3} & b_3
\\
\red{\bf 0} & \red{\bf 0} & \red{\bf 0} & d
\end{array}
\right]^{\green{\bf 1}}.
\]
\end{Lemma}

\proof
Just read~{\eqref{eqG-S3-R4}} at order 1:
\[
0
{\overset{\green{\bf 1,0,0}}{\,\,=\,\,}}
-\,c_1,
\ \ \ \ \ \ \ \ \ \ \ \ \ \ \ \ \ \ \ \
0
{\overset{\green{\bf 0,1,0}}{\,\,=\,\,}}
-\,c_2,
\ \ \ \ \ \ \ \ \ \ \ \ \ \ \ \ \ \ \ \
0
{\overset{\green{\bf 0,0,1}}{\,\,=\,\,}}
-\,c_3.
\qedhere
\]
\endproof

\SectionHead{Quadratic Normalizations and Isotropy Reductions}
{quadratic-normalizations-isotropy-reductions}

At order 2, a quadratic form $F_2$ appears:
\[
u
\,=\,
F_2
+
{\rm O}_{x,y,z}(3)
\,=\,
F_{2,0,0}\,x^2
+
F_{1,1,0}\,x\,y
+
F_{1,0,1}\,x\,z
+
F_{0,2,0}\,y^2
+
F_{0,1,1}\,y\,z
+
F_{0,0,2}\,z^2
+
\cdots,
\]
At the origin,
the Hessian matrix is:
\[
\left[\!
\begin{array}{ccc}
2F_{200} & F_{110} & F_{101}
\\
F_{110} & 2F_{020} & F_{011}
\\
F_{101} & F_{011} & 2F_{002}
\end{array}
\!\right]
\,\,\,\sim\,\,\,
\left[\!
\begin{array}{ccc}
2F_{200} & F_{110} & 0
\\
F_{110} & 2F_{020} & 0
\\
0 & 0 & 0
\end{array}
\!\right],
\]
and due to Hypothesis~{\ref{Hyp-Hessian-rank-2}},
after an elementary affine transformation
in the $(x,y,z)$-space,
only its upper-left $2 \times 2$ block is nonzero,
and of rank 2.

It is well known that, under the group 
$\GL(2,\R) \ni \big( \begin{smallmatrix} \alpha & \beta \\
\gamma & \delta \end{smallmatrix} \big)$, with $x \longmapsto
\alpha\,x + \beta\,y$ and $y \longmapsto \gamma\, x + \delta\,y$,
any quadratic form $F_{2,0}\, x^2 + F_{1,1}\, xy + F_{0,2}\, y^2$
of rank 2 can be normalized in two inequivalent ways as:
\[
\def\arraystretch{1.25}
\begin{array}{cc}
\,\,{\sf Hyperbolic}\,\, & \,\,{\sf Elliptic}\,\,
\\
x\,y & x^2+y^2
\end{array}
\]
the quadratic form $x\,y$, of signature $(1, 1)$,
being equivalent to $x^2 - y^2$.

\begin{Proposition}
At order 2, two nonequivalent branches exist:
\reqnomode\usetagform{EngelLie}
\begin{align}
\text{\bf Branch 2a}
\ \ \ \ \ \ \ \ \ \ \ \ \ \ \ \ \ \ \ \
u
&
\,=\,
x\,y
\ \ \ \ \ \ \ \,
+
{\rm O}_{x,y,z}(3),
\notag
\\
\text{\bf Branch 2b}
\ \ \ \ \ \ \ \ \ \ \ \ \ \ \ \ \ \ \ \
u
&
\,=\,
x^2+y^2
+
{\rm O}_{x,y,z}(3).
\tag{\qed}
\end{align}
\end{Proposition}

Once a normalization has been made in the left
space $\R^4 \ni (x,y,z,u)$,
always, it can also be made {\em exactly the same} in the right
space $(p,q,r,v) \in \R^4$:
\[
\aligned
\text{\bf Branch 2a}
\ \ \ \ \ \ \ \ \ \ \ \ \ \ \ \ \ \ \ \
v
&
\,=\,
p\,q
\ \ \ \ \ \ \ \,
+
{\rm O}_{p,q,r}(3),
\notag
\\
\text{\bf Branch 2b}
\ \ \ \ \ \ \ \ \ \ \ \ \ \ \ \ \ \ \ \
v
&
\,=\,
p^2+q^2
+
{\rm O}_{p,q,r}(3).
\endaligned
\]

\begin{Principle}
\label{Principle-FG}
{\sl At each order, every performed normalization will always 
be instantly achieved on both hypersurfaces
$\{u = F\}$ and $\{v = G\}$.}\qed
\end{Principle}

To terminate order 2, we must determine the group reductions
in the two branches.

\begin{Lemma}
\label{Lm-stab-order-2}
Stabilization of order 2 terms
holds with the group reductions:
\[
\left[
\begin{array}{cccc}
a_{1,1} & a_{1,2} & a_{1,3} & b_1
\\
a_{2,1} & a_{2,2} & a_{2,3} & b_2
\\
a_{3,1} & a_{3,2} & a_{3,3} & b_3
\\
\red{\bf 0} & \red{\bf 0} & \red{\bf 0} & d
\end{array}
\right]^{\green{\bf 1}}
\,\,\,\leadsto\,\,\,
\aligned
{}
&
\green{\bf 2a:}
\ \ \ \ \ \ \ \ \ 
\left[
\begin{array}{cccc}
a_{1,1} & \red{\bf 0} & \red{\bf 0} & b_1
\\
\red{\bf 0} & a_{2,2} & \red{\bf 0} & b_2
\\
a_{3,1} & a_{3,2} & a_{3,3} & b_3
\\
\red{\bf 0} & \red{\bf 0} & \red{\bf 0} & \red{a_{1,1}a_{2,2}}
\end{array}
\right]^{\green{\bf 2}}
\\
{}
&
\green{\bf 2b:}
\ \ \ \ \ \ \ \ \ 
\left[
\begin{array}{cccc}
a_{1,1} & \red{-a_{2,1}} & \red{\bf 0} & b_1
\\
a_{2,1} & \red{a_{1,1}} & \red{\bf 0} & b_2
\\
a_{3,1} & a_{3,2} & a_{3,3} & b_3
\\
\red{\bf 0} & \red{\bf 0} & \red{\bf 0} & \red{a_{1,1}^2+a_{2,1}^2}
\end{array}
\right]^{\green{\bf 2}}
\endaligned
\]
with nonzero determinants:
\[
a_{1,1}^2\,a_{2,2}^2\,a_{3,3}
\,\neq\,
0,
\ \ \ \ \ \ \ \ \ \ \ \ \ \ \ \ \ \ \ \ \ \ \ \ \ \ \ \ \ \ \ \ \ \ \
0
\,\neq\,
\big(
a_{1,1}^2
+
a_{2,1}^2
\big)^2\,
a_{3,3}.
\]
\end{Lemma}

\proof
Read~{\eqref{eqG-S3-R4}} at order 2, for branch $\green{\bf 2a}$
and for branch $\green{\bf 2b}$:
\[
\aligned
&
0
{\overset{\green{\bf 2,0,0}}{\,\,=\,\,}}
a_{1,1}\,a_{2,1},
\\
&
0
{\overset{\green{\bf 1,1,0}}{\,\,=\,\,}}
a_{1,1}\,a_{2,2}
+
a_{1,2}\,a_{2,1}
-
d,
\\
&
0
{\overset{\green{\bf 0,2,0}}{\,\,=\,\,}}
a_{1,2}\,a_{2,2},
\\
&
0
{\overset{\green{\bf 1,0,1}}{\,\,=\,\,}}
a_{1,1}\,a_{2,3}
+
a_{1,3}\,a_{2,1},
\\
&
0
{\overset{\green{\bf 0,1,1}}{\,\,=\,\,}}
a_{1,2}\,a_{2,3}
+
a_{1,3}\,a_{2,2},
\\
&
0
{\overset{\green{\bf 0,2,0}}{\,\,=\,\,}}
a_{1,3}\,a_{2,3},
\endaligned
\ \ \ \ \ \ \ \ \ \ \ \ \ \ \ \ \ \ \ \ \ \ \ \ \ \ \ \ \ \ \ \ \ \ \
\aligned
&
0
{\overset{\green{\bf 2,0,0}}{\,\,=\,\,}}
a_{1,1}^2
+
a_{2,1}^2
-
d,
\\
&
0
{\overset{\green{\bf 1,1,0}}{\,\,=\,\,}}
2\,a_{1,1}\,a_{1,2}
+
2\,a_{2,1}\,a_{2,2},
\\
&
0
{\overset{\green{\bf 0,2,0}}{\,\,=\,\,}}
a_{1,2}^2
+
a_{2,2}^2
-
d,
\\
&
0
{\overset{\green{\bf 1,0,1}}{\,\,=\,\,}}
2\,a_{1,1}\,a_{1,3}
+
2\,a_{2,1}\,a_{2,3},
\\
&
0
{\overset{\green{\bf 0,1,1}}{\,\,=\,\,}}
2\,a_{1,2}\,a_{1,3}
+
2\,a_{2,2}\,a_{2,3},
\\
&
0
{\overset{\green{\bf 0,2,0}}{\,\,=\,\,}}
a_{1,3}^2
+
a_{2,3}^2,
\endaligned
\]
and solve for the $a_{\smallbullet, \smallbullet}$ and $d$,
in the connected component of the identity.
\endproof

\SectionHead{Infinitesimal Affine Transformations}
{infinitesimal-affine-transformations}

A general affine vector field writes:
\[
\aligned
L
&
\,=\,
\ \ 
\big(
T_1+A_{1,1}\,x+A_{1,2}\,y+A_{1,3}\,z+B_1\,u
\big)\,\frac{\partial}{\partial x}
\\
&
\ \ \ \ \
+
\big(
T_2+A_{2,1}\,x+A_{2,2}\,y+A_{2,3}\,z+B_2\,u
\big)\,\frac{\partial}{\partial y}
\\
&
\ \ \ \ \
+
\big(
T_3+A_{3,1}\,x+A_{3,2}\,y+A_{3,3}\,z+B_3\,u
\big)\,\frac{\partial}{\partial z}
\\
&
\ \ \ \ \
+
\big(
U_0+C_1\,x+C_2\,y+C_3\,z+D\,u
\big)\,\frac{\partial}{\partial u}.
\endaligned
\]
It is tangent to $\big\{u = F(x,y,z) \big\}$ 
if and only if:
\[
0
\,\equiv\,
\eqLF(x,y,z)
\,=:\,
L
\big(
-\,u+F(x,y,z)
\big)
\Big\vert_{u=F(x,y,z)},
\]
identically as power series in $\R\{x,y,z\}$.
With increasing orders $\mu = 0, 1, 2, 3, \dots$, 
this $\eqLF$ may be expanded:
\[
0
\,=\,
\eqLF
\,=\,
\sum_{\mu=0}^\infty\,\,
\sum_{i+j+k=\mu}\,
\mathcal{V}_{i,j,k}
\Big(
T_{\smallbullet},\,
U_{\smallbullet},\,
A_{\smallbullet,\smallbullet},\,
B_{\smallbullet},\,
C_{\smallbullet,\smallbullet},\,
D_{\smallbullet},\,\,
F_{\smallbullet,\smallbullet,\smallbullet}
\Big)\,
x^i\,y^j\,z^k,
\]
so that, for all $i, j, k \in \N$:
\[
0
{\overset{\green{\bf i,j,k}}{\,\,=\,\,}}
\mathcal{V}_{i,j,k}.
\]

Such a vector field $L$ is tangent to: 
\[
u
\,=\,
\red{\bf 0}
+
{\rm O}_{x,y,z}(2),
\]
if and only if:
\[
0
{\overset{\green{\bf 0,0,0}}{\,\,=\,\,}}
-\,U_0.
\]

The key constraint of transitivity:
\[
\Span\,
\big(
\tfrac{\partial}{\partial x},\,
\tfrac{\partial}{\partial y},\,
\tfrac{\partial}{\partial z}
\big)
\,=\,
T_{\sf origin}S
\,=\,
\Span\,
L
\big\vert_{\sf origin}
\,=\,
\Span\,
\Big(
T_1\,
\tfrac{\partial}{\partial x}
+
T_2\,
\tfrac{\partial}{\partial y}
+
T_3\,
\tfrac{\partial}{\partial z}
\Big),
\]
forces to always keep $T_1$, $T_2$, $T_3$ 
absolutely free\,\,---\,\,never solved.

\SectionHead{Order 3 Linear Representations for Branches 2a and 2b}
{order-3-linear-representations}

If an affine infinitesimal transformation $L$ with $U_0 := 0$ 
is tangent to:
\[
u
\overset{\green{\bf 2a}}{\,=\,}
x\,y
+
{\rm O}_{x,y,z}(3),
\ \ \ \ \ \ \ \ \ \ \ \ \ \ \ \ \ \ \ \
\ \ \ \ \ \ \ \ \ \ \ \ \ \ \ \ \ \ \ \
u
\overset{\green{\bf 2b}}{\,=\,}
x^2+y^2
+
{\rm O}_{x,y,z}(3),
\]
then at order 1:
\[
\aligned
0
&
{\overset{\green{\bf 1,0,0}}{\,\,=\,\,}}
-\,C_1+T_2,
\\
0
&
{\overset{\green{\bf 0,1,0}}{\,\,=\,\,}}
-\,C_2+T_1,
\\
0
&
{\overset{\green{\bf 0,0,1}}{\,\,=\,\,}}
-\,C_3,
\endaligned
\ \ \ \ \ \ \ \ \ \ \ \ \ \ \ \ \ \ \ \
\ \ \ \ \ \ \ \ \ \ \ \ \ \ \ \ \ \ \ \
\aligned
0
&
{\overset{\green{\bf 1,0,0}}{\,\,=\,\,}}
-\,C_1+2\,T_1,
\\
0
&
{\overset{\green{\bf 0,1,0}}{\,\,=\,\,}}
-\,C_2+2\,T_2,
\\
0
&
{\overset{\green{\bf 0,0,1}}{\,\,=\,\,}}
-\,C_3.
\endaligned
\]

Next, after assigning:
\[
C_1
\,:=\,
T_2,
\ \ \ \ \ \ \
C_2
\,:=\,
T_1,
\ \ \ \ \ \ \
C_3
\,:=\,
0
\ \ \ \ \ \ \ \ \ \ \ \ \ \ \ \ \ \ \ \
\ \ \ \ \ \ \ \ \ \ \ \ \ \ \ \ \ \ \ \
C_1
\,:=\,
2\,T_1,
\ \ \ \ \ \ \
C_2
\,:=\,
2\,T_2,
\ \ \ \ \ \ \
C_3
\,:=\,
0,
\]
the tangency of the affine infinitesimal transformation $L$ to:
\[
u
\,=\,
x\,y
+
\sum_{i+j+k=3}\,
F_{i,j,k}\,
x^i\,y^j\,z^k
+
{\rm O}_{x,y,z}(4),
\ \ \ \ \ 
\ \ \ \ \ 
u
\,=\,
x^2+y^2
+
\sum_{i+j+k=3}\,
F_{i,j,k}\,
x^i\,y^j\,z^k
+
{\rm O}_{x,y,z}(4),
\]
gives at order 2:
\[
\footnotesize
\aligned
0
&
{\overset{\green{\bf 2,0,0}}{\,\,=\,\,}}
A_{2,1}
+
3\,F_{3,0,0}\,T_1
+
F_{2,1,0}\,T_2
+
F_{2,0,1}\,T_3,
\\
0
&
{\overset{\green{\bf 1,1,0}}{\,\,=\,\,}}
A_{1,1}+A_{2,2}-D
+
2\,F_{2,1,0}\,T_1
+
2\,F_{1,2,0}\,T_2
+
F_{1,1,1}\,T_3,
\\
0
&
{\overset{\green{\bf 1,0,1}}{\,\,=\,\,}}
A_{2,3}
+
2\,F_{2,0,1}\,T_1
+
F_{1,1,1}\,T_2
+
2\,F_{1,0,2}\,T_3,
\\
0
&
{\overset{\green{\bf 0,2,0}}{\,\,=\,\,}}
A_{1,2}
+
F_{1,2,0}\,T_1
+
3\,F_{0,3,0}\,T_2
+
F_{0,2,1}\,T_3,
\\
0
&
{\overset{\green{\bf 0,1,1}}{\,\,=\,\,}}
A_{1,3}
+
F_{1,1,1}\,T_1
+
2\,F_{0,2,1}\,T_2
+
2\,F_{0,1,2}\,T_3,
\\
0
&
{\overset{\green{\bf 0,0,2}}{\,\,=\,\,}}
F_{1,0,2}\,T_1
+
F_{0,1,2}\,T_2
+
3\,F_{0,0,3}\,T_3,
\endaligned
\ \ \ \ \
\footnotesize
\aligned
0
&
{\overset{\green{\bf 2,0,0}}{\,\,=\,\,}}
2\,A_{1,1}-D
+
3\,F_{3,0,0}\,T_1
+
F_{2,1,0}\,T_2
+
F_{2,0,1}\,T_3,
\\
0
&
{\overset{\green{\bf 1,1,0}}{\,\,=\,\,}}
2\,A_{1,2}+2\,A_{2,1}
+
2\,F_{2,1,0}\,T_1
+
2\,F_{1,2,0}\,T_2
+
F_{1,1,1}\,T_3,
\\
0
&
{\overset{\green{\bf 1,0,1}}{\,\,=\,\,}}
2\,A_{1,3}
+
2\,F_{2,0,1}\,T_1
+
F_{1,1,1}\,T_2
+
2\,F_{1,0,2}\,T_3,
\\
0
&
{\overset{\green{\bf 0,2,0}}{\,\,=\,\,}}
2\,A_{2,2}-D
+
F_{1,2,0}\,T_1
+
3\,F_{0,3,0}\,T_2
+
F_{0,2,1}\,T_3,
\\
0
&
{\overset{\green{\bf 0,1,1}}{\,\,=\,\,}}
2\,A_{2,3}
+
F_{1,1,1}\,T_1
+
2\,F_{0,2,1}\,T_2
+
2\,F_{0,1,2}\,T_3,
\\
0
&
{\overset{\green{\bf 0,0,2}}{\,\,=\,\,}}
F_{1,0,2}\,T_1
+
F_{0,1,2}\,T_2
+
3\,F_{0,0,3}\,T_3.
\endaligned
\]
In each branch (column) \green{\bf 2a} and \green{\bf 2b},
the first five equations can be solved, progressively, as:
\[
\footnotesize
\aligned
D
&
\,:=\,
A_{1,1}+A_{2,2}
+
2\,F_{2,1,0}\,T_1
+
2\,F_{1,2,0}\,T_2
+
F_{1,1,1}\,T_3,
\\
A_{2,1}
&
\,:=\,
-\,3\,F_{3,0,0}\,T_1
-
F_{2,1,0}\,T_2
-
F_{2,0,1}\,T_3,
\\
A_{1,2}
&
\,:=\,
-\,F_{1,2,0}\,T_1
-
3\,F_{0,3,0}\,T_2
-
F_{0,2,1}\,T_3,
\\
A_{2,3}
&
\,:=\,
-\,2\,F_{2,0,1}\,T_1
-
F_{1,1,1}\,T_2
-
2\,F_{1,0,2}\,T_3,
\\
A_{1,3}
&
\,:=\,
-\,F_{1,1,1}\,T_1
-
2\,F_{0,2,1}\,T_2
-
2\,F_{0,1,2}\,T_3,
\endaligned
\ \ \ \ \ 
\footnotesize
\aligned
A_{2,3}
&
\,:=\,
-\,\tfrac{1}{2}\,F_{1,1,1}\,T_1
-
F_{0,2,1}\,T_2
-
F_{1,0,2}\,T_3,
\\
A_{1,3}
&
\,:=\,
-\,F_{2,0,1}\,T_1
-
\tfrac{1}{2}\,F_{1,1,1}\,T_2
-
F_{1,0,2}\,T_3,
\\
D
&
\,:=\,
2\,A_{2,2}
+
F_{1,2,0}\,T_1
+
3\,F_{0,3,0}\,T_2
+
F_{0,2,1}\,T_3,
\\
A_{1,2}
&
\,:=\,
-\,A_{2,1}
-
F_{2,1,0}\,T_1
-
F_{1,2,0}\,T_2
-
\tfrac{1}{2}\,F_{1,1,1}\,T_3,
\\
A_{2,2}
&
\,:=\,
A_{1,1}
+
\Big(
\tfrac{3}{2}\,F_{3,0,0}
-
\tfrac{1}{2}\,F_{1,2,0}
\Big)\,T_1
\\
&
\ \ \ \ \
+
\Big(
\tfrac{1}{2}\,F_{2,1,0}
-
\tfrac{3}{2}\,F_{0,3,0}
\Big)\,T_2
+
\Big(
\tfrac{1}{2}\,F_{2,0,1}
-
\tfrac{1}{2}\,F_{0,2,1}
\Big)\,T_3.
\endaligned
\]

Concerning the two sixths equations 
${\overset{\green{\bf 0,0,2}}{\,=\,}}$, 
since the transitivity parameters $T_1$, $T_2$, $T_3$ must remain
absolutely free of any linear relations, 
it necessarily holds in the source space that:
\[
0
\,=\,
F_{1,0,2}
\,=\,
F_{0,1,2}
\,=\,
F_{0,0,3},
\ \ \ \ \ \ \ \ \ \ \ \ \ \ \ \ \ \ \ \
\ \ \ \ \ \ \ \ \ \ \ \ \ \ \ \ \ \ \ \
0
\,=\,
F_{1,0,2}
\,=\,
F_{0,1,2}
\,=\,
F_{0,0,3}.
\]
Of course simililary in the target space, according
to Principle~{\ref{Principle-FG}}, it holds for the two branches
{\green{\bf 2a}} and {\green{\bf 2b}}
that:
\[
0
\,=\,
G_{1,0,2}
\,=\,
G_{0,1,2}
\,=\,
G_{0,0,3},
\ \ \ \ \ \ \ \ \ \ \ \ \ \ \ \ \ \ \ \
\ \ \ \ \ \ \ \ \ \ \ \ \ \ \ \ \ \ \ \
0
\,=\,
G_{1,0,2}
\,=\,
G_{0,1,2}
\,=\,
G_{0,0,3}.
\]

Next, in Branch {\green{\bf 2a}} (only here), at order~3,
the equations of $\eqFG$ are:
\[
\aligned
0
&
{\overset{\green{\bf 3,0,0}}{\,\,=\,\,}}
-\,a_{1,1}a_{2,2}\,F_{3,0,0}
+
a_{1,1}^3\,G_{3,0,0}
+
a_{1,1}^2\,a_{3,1}\,G_{2,0,1},
\\
0
&
{\overset{\green{\bf 2,1,0}}{\,\,=\,\,}}
-\,a_{1,1}a_{2,2}\,F_{2,1,0}
+
a_{1,1}^2a_{2,2}\,G_{2,1,0}
+
a_{1,1}^2a_{3,2}\,G_{2,0,1}
+
a_{1,1}a_{2,2}a_{3,1}\,G_{1,1,1}
+
a_{1,1}\,\boxed{b_2},
\\
0
&
{\overset{\green{\bf 1,2,0}}{\,\,=\,\,}}
-\,a_{1,1}a_{2,2}\,F_{1,2,0}
+
a_{1,1}a_{2,2}^2\,G_{1,2,0}
+
a_{1,1}a_{2,2}a_{3,2}\,G_{1,1,1}
+
a_{2,2}^2a_{3,1}\,G_{0,2,1}
+
a_{2,2}\,\boxed{b_1},
\\
0
&
{\overset{\green{\bf 0,3,0}}{\,\,=\,\,}}
-\,a_{1,1}a_{2,2}\,F_{0,3,0}
+
a_{2,2}^3\,G_{0,3,0}
+
a_{2,2}^2a_{3,2}\,G_{0,2,1},
\\
0
&
{\overset{\green{\bf 2,0,1}}{\,\,=\,\,}}
-\,a_{1,1}a_{2,2}\,F_{2,0,1}
+
a_{1,1}^2a_{3,3}\,G_{2,0,1},
\\
0
&
{\overset{\green{\bf 1,1,1}}{\,\,=\,\,}}
-\,a_{1,1}a_{2,2}\,F_{1,1,1}
+
a_{1,1}a_{2,2}a_{3,3}\,G_{1,1,1},
\\
0
&
{\overset{\green{\bf 0,2,1}}{\,\,=\,\,}}
-\,a_{1,1}a_{2,2}\,F_{0,2,1}
+
a_{2,2}^2a_{3,3}\,G_{0,2,1}.
\endaligned
\]
(In Branch {\green{\bf 2b}}, the equations are more complicated.)
Since Lemma~{\ref{Lm-stab-order-2}} showed that the coefficients
$a_{1,1}$ and $a_{2,2}$ of the boxed $b_2$ and $b_1$ are nonzero, 
we can use these two group parameters $b_2$ and $b_1$ to normalize:
\[
G_{2,1,0}
\,:=\,
0
\,=:\,
F_{2,1,0},
\ \ \ \ \ \ \ \ \ \ \ \ \ \ \ \ \ \ \ \ \ \ \ \ \ \ \ \ \ \ \ \ \ \ \
G_{1,2,0}
\,:=\,
0
\,=:\,
F_{1,2,0}.
\]
Once this is done, we may solve from the five remaining equations:
\[
\left(\!
\begin{array}{c}
G_{2,0,1}
\\
G_{1,1,1}
\\
G_{0,2,1}
\\
G_{3,0,0}
\\
G_{0,3,0}
\end{array}
\!\right)
\,=\,
\left(\!
\def\arraystretch{1.25}
\begin{array}{ccccc}
\frac{a_{2,2}}{a_{1,1}a_{3,3}} & 0 & 0 & 0 & 0
\\
0 & \frac{1}{a_{3,3}} & 0 & 0 & 0
\\
0 & 0 & \frac{a_{1,1}}{a_{2,2}a_{3,3}} & 0 & 0
\\
-\frac{a_{2,2}a_{3,1}}{a_{1,1}^2a_{3,3}} & 0 & 0 & 
\frac{a_{2,2}}{a_{1,1}^2} & 0
\\
0 & 0 & -\frac{a_{1,1}a_{3,2}}{a_{2,2}^2a_{3,3}} & 0 &
\frac{a_{1,1}}{a_{2,2}^2}
\end{array}
\!\right)\,
\left(\!
\begin{array}{c}
F_{2,0,1}
\\
F_{1,1,1}
\\
F_{0,2,1}
\\
F_{3,0,0}
\\
F_{0,3,0}
\end{array}
\!\right),
\]
to receive a linear representation of a certain matrix-group
on a 5-dimensional vector space. At order 3, skipping details,
this leads us to the creation of 19 branches,
mutually inequivalent and having empty intersection: 
\[
\def\arraystretch{1.25}
\begin{array}{rccccc}
\green{\bf 2a}\,\,\,\,
\green{\downarrow}\,\,
& 
F_{2,0,1} & F_{1,1,1} & F_{0,2,1} & F_{3,0,0} & F_{0,3,0}
\\
\green{\bf 3a} & 
0 & 0 & 0 & 0 & 0
\\
\green{\bf 3b} & 
0 & 0 & 0 & 0 & 1
\\
\green{\bf 3c} & 
0 & 0 & 0 & 1 & 0
\\
\green{\bf 3d} &
0 & 0 & 0 & 1 & 1
\\
\green{\bf 3e} & 
0 & 0 & 1 & 0 & 0
\\
\green{\bf 3f} & 
0 & 0 & 1 & 1 & 0
\\
\green{\bf 3g} & 
0 & 1 & 0 & 0 & 0
\\
\green{\bf 3h} & 
0 & 1 & 0 & 0 & 1
\\
\green{\bf 3i} & 
0 & 1 & 0 & 1 & 0
\\
\green{\bf 3j} & 
0 & 1 & 0 & 1 & 1
\\
\green{\bf 3k} & 
0 & 1 & 1 & 0 & 0
\\
\green{\bf 3l} & 
0 & 1 & 1 & 1 & 0
\\
\green{\bf 3m} & 
1 & 0 & 0 & 0 & 0 
\\
\green{\bf 3n} &
1 & 0 & 0 & 0 & 1
\\
\green{\bf 3o} & 
1 & 0 & 1 & 0 & 0
\\
\green{\bf 3p} & 
1 & 0 & -1 & 0 & 0
\\
\green{\bf 3q} & 
1 & 1 & F_{0,2,1}^{\neq0} & 0 & 0
\\
\green{\bf 3r} & 
1 & 1 & 0 & 0 & 0
\\
\green{\bf 3s} & 
1 & 1 & 0 & 0 & 1
\end{array}
\]
In branch \green{\bf 2a3q}, the absolute invariant
$F_{0,2,1} = G_{0,2,1}$ is assumed to be $\neq 0$.

\medskip

Next, in branch {\green{\bf 2b}}, using $b_1$ and $b_2$,
we can similarly normalize:
\leqnomode\usetagform{default}
\begin{align}
\label{normalize-F210-F120}
G_{2,1,0}
\,:=\,
0
\,=:\,
F_{2,1,0},
\ \ \ \ \ \ \ \ \ \ \ \ \ \ \ \ \ \ \ \ \ \ \ \ \ \ \ \ \ \ \ \ \ \ \
G_{1,2,0}
\,:=\,
0
\,=:\,
F_{1,2,0}.
\end{align}

Once this is done, we may solve 
$G_{2,0,1}$, $G_{1,1,1}$, $G_{0,2,1}$, $G_{3,0,0}$,
$G_{0,3,0}$ from the five remaining equations:
\[
\left(\!
\begin{array}{c}
G_{2,0,1}
\\
G_{1,1,1}
\\
G_{0,2,1}
\\
G_{3,0,0}
\\
G_{0,3,0}
\end{array}
\!\right)
\,=\,
\left(\!
\def\arraystretch{1.25}
\begin{array}{ccccc}
\frac{a_{1,1}^2}{(a_{1,1}^2+a_{2,1}^2)\,a_{3,3}} & 
\frac{-\,a_{1,1}\,a_{2,1}}{(a_{1,1}^2+a_{2,1}^2)\,a_{3,3}} & 
\frac{a_{2,1}^2}{(a_{1,1}^2+a_{2,1}^2)\,a_{3,3}} & 
0 & 0
\\
\frac{2\,a_{1,1}\,a_{2,1}}{(a_{1,1}^2+a_{2,1}^2)\,a_{3,3}} & 
\frac{a_{1,1}^2-a_{2,1}^2}{(a_{1,1}^2+a_{2,1}^2)\,a_{3,3}} & 
\frac{-\,2\,a_{1,1}\,a_{2,1}}{(a_{1,1}^2+a_{2,1}^2)\,a_{3,3}} & 
0 & 0
\\
\frac{a_{2,1}^2}{(a_{1,1}^2+a_{2,1}^2)\,a_{3,3}} & 
\frac{a_{1,1}\,a_{2,1}}{(a_{1,1}^2+a_{2,1}^2)\,a_{3,3}} & 
\frac{a_{1,1}^2}{(a_{1,1}^2+a_{2,1}^2)\,a_{3,3}} & 
0 & 0
\\
Q_{4,1} & Q_{4,2} & Q_{4,3} & 
\frac{a_{1,1}(a_{1,1}^2-3\,a_{2,1}^2)}{
(a_{1,1}^2+a_{2,1}^2)^2} & 
\frac{a_{2,1}(3\,a_{1,1}^2-a_{2,1}^2)}{
(a_{1,1}^2+a_{2,1}^2)^2}
\\
Q_{5,1} & Q_{5,2} & Q_{5,3} &
\frac{-\,a_{2,1}(3\,a_{1,1}^2-a_{2,1}^2)}{
(a_{1,1}^2+a_{2,1}^2)^2} &
\frac{a_{1,1}(a_{1,1}^2-3\,a_{2,1}^2)}{
(a_{1,1}^2+a_{2,1}^2)^2}
\end{array}
\!\right)\,
\left(\!
\begin{array}{c}
F_{2,0,1}
\\
F_{1,1,1}
\\
F_{0,2,1}
\\
F_{3,0,0}
\\
F_{0,3,0}
\end{array}
\!\right),
\]
where:
\[
\footnotesize
\aligned
Q_{4,1}
&
\,=\,
\frac{
-\,a_{1,1}^3\,a_{3,1}-3\,a_{1,1}^2\,a_{2,1}\,a_{3,2}
-3\,a_{1,1}\,a_{2,1}^2\,a_{3,1}+a_{2,1}^3\,a_{3,2}
}{\big(a_{1,1}^2+a_{2,1}^2\big)^2\,a_{3,3}},
\\
Q_{4,2}
&
\,=\,
\frac{
a_{1,1}^3\,a_{3,2}+3\,a_{1,1}^2\,a_{2,1}\,a_{3,1}
-3\,a_{1,1}\,a_{2,1}^2\,a_{3,2}-a_{2,1}^3\,a_{3,1}
}{\big(a_{1,1}^2+a_{2,1}^2\big)^2\,a_{3,3}},
\\
Q_{4,3}
&
\,=\,
\frac{
a_{1,1}^3\,a_{3,1}-3\,a_{1,1}^2\,a_{2,1}\,a_{3,2}
-3\,a_{1,1}\,a_{2,1}^2\,a_{3,1}+a_{2,1}^3\,a_{3,2}
}{\big(a_{1,1}^2+a_{2,1}^2\big)^2\,a_{3,3}},
\\
Q_{5,1}
&
\,=\,
\frac{
a_{1,1}^3\,a_{3,2}+3\,a_{1,1}^2\,a_{2,1}\,a_{3,1}
-3\,a_{1,1}\,a_{2,1}^2\,a_{3,2}-a_{2,1}^3\,a_{3,1}
}{\big(a_{1,1}^2+a_{2,1}^2\big)^2\,a_{3,3}},
\\
Q_{5,2}
&
\,=\,
\frac{
a_{1,1}^3\,a_{3,1}-3\,a_{1,1}^2\,a_{2,1}\,a_{3,2}
-3\,a_{1,1}\,a_{2,1}^2\,a_{3,1}+a_{2,1}^3\,a_{3,2}
}{\big(a_{1,1}^2+a_{2,1}^2\big)^2\,a_{3,3}},
\\
Q_{5,3}
&
\,=\,
\frac{
-\,a_{1,1}^3\,a_{3,2}+3\,a_{1,1}^2\,a_{2,1}\,a_{3,1}
-3\,a_{1,1}\,a_{2,1}^2\,a_{3,2}-a_{2,1}^3\,a_{3,1}
}{\big(a_{1,1}^2+a_{2,1}^2\big)^2\,a_{3,3}},
\endaligned
\]
and the determinant of this $5 \times 5$ matrix is equal to:
\[
\frac{1}{\big(a_{1,1}^2+a_{2,1}^2\big)\,a_{3,3}^3}
\,\neq\,
0.
\]

In this branch~\green{\bf 2b}, the creation of order 3 branches
is more subtle than in branch~\green{\bf 2a}, 
hence we provide details. 

At first, with the 3 group parameters 
$\big(a_{1,1}, a_{2,1}, a_{3,3} \big)$,
there is a 3D sub-representation:
\[
\left(\!
\begin{array}{c}
G_{2,0,1}
\\
G_{1,1,1}
\\
G_{0,2,1}
\end{array}
\!\right)
\,=\,
\left(\!
\def\arraystretch{1.25}
\begin{array}{ccc}
\frac{a_{1,1}^2}{(a_{1,1}^2+a_{2,1}^2)\,a_{3,3}} & 
\frac{-\,a_{1,1}\,a_{2,1}}{(a_{1,1}^2+a_{2,1}^2)\,a_{3,3}} & 
\frac{a_{2,1}^2}{(a_{1,1}^2+a_{2,1}^2)\,a_{3,3}} 
\\
\frac{2\,a_{1,1}\,a_{2,1}}{(a_{1,1}^2+a_{2,1}^2)\,a_{3,3}} & 
\frac{a_{1,1}^2-a_{2,1}^2}{(a_{1,1}^2+a_{2,1}^2)\,a_{3,3}} & 
\frac{-\,2\,a_{1,1}\,a_{2,1}}{(a_{1,1}^2+a_{2,1}^2)\,a_{3,3}} 
\\
\frac{a_{2,1}^2}{(a_{1,1}^2+a_{2,1}^2)\,a_{3,3}} & 
\frac{a_{1,1}\,a_{2,1}}{(a_{1,1}^2+a_{2,1}^2)\,a_{3,3}} & 
\frac{a_{1,1}^2}{(a_{1,1}^2+a_{2,1}^2)\,a_{3,3}} 
\end{array}
\!\right)\,
\left(\!
\begin{array}{c}
F_{2,0,1}
\\
F_{1,1,1}
\\
F_{0,2,1}
\end{array}
\!\right),
\]
having determinant $\frac{1}{a_{3,3}^3}$ that must be nonzero.

It can be verified that the initial 5D linear representation
is {\em not} decomposable as a direct sum of
a 3D representation plus a 2D representation, but this fact
will not be useful. The identity $3 \times 3$ matrix
corresponds to the value $(1, 0, 1)$ of 
$(a_{1,1}, a_{2,1}, a_{3,3})$. Setting:
\[
\aligned
a_{1,1}
&
\,=\,
1
+
1\cdot\varepsilon,
\\
a_{2,1}
&
\,=\,
0
+
0\cdot\varepsilon,
\\
a_{3,3}
&
\,=\,
1
+
0\cdot\varepsilon,
\endaligned
\ \ \ \ \ \ \ \ \ \ \ \ \ \ \ \ \ \ \ \ \ \ \ \ \ \
\aligned
a_{1,1}
&
\,=\,
1
+
0\cdot\varepsilon,
\\
a_{2,1}
&
\,=\,
0
+
1\cdot\varepsilon,
\\
a_{3,3}
&
\,=\,
1
+
0\cdot\varepsilon,
\endaligned
\ \ \ \ \ \ \ \ \ \ \ \ \ \ \ \ \ \ \ \ \ \ \ \ \ \
\aligned
a_{1,1}
&
\,=\,
1
+
0\cdot\varepsilon,
\\
a_{2,1}
&
\,=\,
0
+
0\cdot\varepsilon,
\\
a_{3,3}
&
\,=\,
1
+
1\cdot\varepsilon,
\endaligned
\]
and differentiating $\frac{d}{d\varepsilon} \big\vert_{\varepsilon=0}$
this matrix-action,
we receive the three infinitesimal generators 
(group-invariant vector fields): 
\[
\aligned
L_1
&
\,:=\,
0,
\\
L_2
&
\,:=\,
-\,F_{1,1,1}\,
\frac{\partial}{\partial F_{2,0,1}}
+
\big(
2\,F_{2,0,1}-2\,F_{0,2,1}
\big)\,
\frac{\partial}{\partial F_{1,1,1}}
+
F_{1,1,1}\,
\frac{\partial}{\partial F_{0,2,1}},
\\
L_3
&
\,:=\,
-\,F_{2,0,1}\,
\frac{\partial}{\partial F_{2,0,1}}
-
F_{1,1,1}\,
\frac{\partial}{\partial F_{1,1,1}}
-
F_{0,2,1}\,
\frac{\partial}{\partial F_{0,2,1}}.
\endaligned
\]
The first vector field can be disregarded. 
Then the locus of points 
$\big(F_{2,0,1}, F_{1,1,1}, F_{0,2,1}\big)$ at which
$L_2$, $L_3$ span a vector subspace of rank $\leqslant 1$ 
is described by equating to $0$ 
all the $2 \times 2$ minors of the coefficients
matrix:
\[
\left[\!
\def\arraystretch{1.25}
\begin{array}{ccc}
-F_{1,1,1} & 2F_{2,0,1}-2F_{0,2,1} & F_{1,1,1}
\\
-F_{2,0,1} & -F_{1,1,1} & -F_{0,2,1}
\end{array}
\!\right].
\]
After some algebraic manipulations (using Gröbner bases), 
appropriate generators of the zero-set of this ideal are:
\[
\aligned
0
&
\overset{1}{\,=\,}
\big(
F_{2,0,1}-F_{0,2,1}
\big)\,
\big(
F_{2,0,1}+F_{0,2,1}
\big),
\\
0
&
\overset{2}{\,=\,}
F_{1,1,1}\,
\big(
F_{2,0,1}+F_{0,2,1}
\big),
\\
0
&
\overset{3}{\,=\,}
F_{1,1,1}^2
+
2\,F_{0,2,1}\,
\big(-F_{2,0,1}+F_{0,2,1}\big).
\endaligned
\]

\begin{Lemma}
The locus where $\rank\, \Span \big( L_1, L_2, L_3 \big) \leqslant 1$
is:
\[
\Big\{
F_{2,0,1}=F_{0,2,1}\,\,
\text{\rm and}\,\,
F_{1,1,1}=0
\Big\}.
\]
\end{Lemma}

\proof
Equation $\overset{1}{=}$ has 2 factors. If one solves 
from the second factor $F_{0,2,1} := - F_{2,0,1}$,
then $\overset{3}{=}$ becomes:
\[
0
\overset{3}{\,=\,}
F_{1,1,1}^2
+
4\,F_{2,0,1}^2,
\]
whence $0 = F_{1,1,1} = F_{2,0,1} = F_{0,2,1}$. 
At the origin $(0, 0, 0)$ indeed, the rank is $0 \leqslant 1$.

If one solves  $F_{0,2,1} := F_{2,0,1}$ from 
the first factor, then the remaining 
two equations become:
\[
\aligned
0
&
\overset{2}{\,=\,}
2\,F_{1,1,1}\,F_{2,0,1},
\\
0
&
\overset{3}{\,=\,}
F_{1,1,1}^2,
\endaligned
\]
whence $F_{1,1,1} = 0$, which concludes the proof.
\endproof

Back to the $3 \times 3$ sub-representation matrix above,
it can be verified that:
\[
\Big\{
F_{2,0,1}=F_{0,2,1}\,\,
\text{\rm and}\,\,
F_{1,1,1}=0
\Big\}
\ \ \ \ \ \ \ \ \ \ \ \ \ \ 
\Longleftrightarrow
\ \ \ \ \ \ \ \ \ \ \ \ \ \
\Big\{
G_{2,0,1}=G_{0,2,1}\,\,
\text{\rm and}\,\,
G_{1,1,1}=0
\Big\}
\]
(a consequence which in fact follows from the preceding lemma
by general arguments of the theory of Lie groups),
by direct computations as follows:
\[
\aligned
G_{2,0,1}
-
G_{0,2,1}
&
\,=\,
\frac{1}{(a_{1,1}^2+a_{2,1}^2)\,a_{3,3}}\,
\Big\{
\big(
a_{1,1}^2-a_{2,1}^2
\big)\,
\big(
F_{2,0,1}
-
F_{0,2,1}
\big)
-
2\,a_{1,1}\,a_{2,1}\,F_{1,1,1}
\Big\},
\\
G_{1,1}
&
\,=\,
\frac{1}{(a_{1,1}^2+a_{2,1}^2)\,a_{3,3}}\,
\Big\{
2\,a_{1,1}\,a_{2,1}\,
\big(
F_{2,0,1}-F_{0,2,1}
\big)
+
\big(
a_{1,1}^2-a_{2,1}^2
\big)\,
F_{1,1,1}
\Big\},
\endaligned
\]
with nonzero determinant:
\[
\left\vert\!
\begin{array}{cc}
a_{1,1}^2-a_{2,1}^2 & -2a_{1,1}a_{2,1}
\\
2a_{1,1}a_{2,1} & a_{1,1}^2-a_{2,1}^2
\end{array}
\!\right\vert
\,=\,
\big(
a_{1,1}^2+a_{2,1}^2
\big)^2
\,\neq\,
0.
\]

Consequently, we must distinguish two (non-overlapping) cases:

\smallskip\noindent{\footnotesize\sf C1}\,\,
$\big(F_{2,0,1},\, 0,\, F_{2,0,1}\big)$;

\smallskip\noindent{\footnotesize\sf C2}\,\,
$\big(F_{1,1,1} \neq 0\big)$ or $\big( F_{1,1,1} = 0\,\,
\text{but}\,\, F_{2,0,1} \neq F_{0,2,1} \big)$.

\smallskip
In case {\footnotesize\sf C1}, dropping $F_{1,1,1}$
which is equal to $0$ and $F_{0,2,1}$
which is equal to $F_{2,0,1}$, we come to the following 3D group 
linear representation:

\[
\left(\!
\begin{array}{c}
G_{2,0,1}
\\
G_{3,0,0}
\\
G_{0,3,0}
\end{array}
\!\right)
\,=\,
\left(\!
\def\arraystretch{1.25}
\begin{array}{ccc}
\frac{1}{a_{3,3}} & 
0 & 0
\\
0 &
\frac{a_{1,1}(a_{1,1}^2-3a_{2,1}^2)}{(a_{1,1}^2+a_{2,1}^2)^2} & 
\frac{a_{2,1}(3a_{1,1}^2-a_{2,1}^2)}{(a_{1,1}^2+a_{2,1}^2)^2} 
\\
0 &
-\frac{a_{2,1}(3a_{1,1}^2-a_{2,1}^2)}{(a_{1,1}^2+a_{2,1}^2)^2} & 
\frac{a_{1,1}(a_{1,1}^2-3a_{2,1}^2)}{(a_{1,1}^2+a_{2,1}^2)^2} 
\end{array}
\!\right)\,
\left(\!
\begin{array}{c}
F_{2,0,1}
\\
F_{3,0,0}
\\
F_{0,3,0}
\end{array}
\!\right),
\]
with determinant $= \frac{1}{(a_{1,1}^2+a_{2,1}^2)\,a_{3,3}}$.
A moment of reflection convinces that 
the lower-right $2 \times 2$ bloc corresponds to a similitude
in $\R^2 \ni (G_{3,0,0}, G_{0,3,0})$, which is confirmed by
the computation of the two infinitesimal generators
corresponding to the two group parameters $a_{1,1}$, $a_{2,1}$:
\[
\aligned
L_{1}
&
\,:=\,
-\,F_{3,0,0}\,\frac{\partial}{\partial F_{3,0,0}}
-
F_{0,3,0}\,\frac{\partial}{\partial F_{0,3,0}},
\\
L_{2}
&
\,:=\,
3\,F_{0,3,0}\,\frac{\partial}{\partial F_{3,0,0}}
-
3\,F_{3,0,0}\,\frac{\partial}{\partial F_{0,3,0}}.
\endaligned
\]
In $\R^2$, the similitude group has only 2 orbits: 
the origin $\{(0,0)\}$ and
its complement, $\R^2 \big\backslash \{(0,0)\}$,
with the point $(1, 0)$ being distinguished.
Thus there are 4 non-intersecting branches:
\[
\def\arraystretch{1.25}
\begin{array}{rccccc}
\green{\bf 2b}\,\,\,\,
\green{\downarrow}\,\,
& 
F_{2,0,1} & F_{3,0,0} & F_{0,3,0}
\\
\green{\bf 3a} & 
0 & 0 & 0
\\
\green{\bf 3b} & 
0 & 1 & 0
\\
\green{\bf 3c} & 
1 & 0 & 0
\\
\green{\bf 3d} & 
1 & 1 & 0
\end{array}
\]
With the values $F_{0,2,1} = F_{2,0,1}$ and $F_{1,1,1} = 0$,
this is:
\[
\def\arraystretch{1.25}
\begin{array}{rccccc}
\green{\bf 2b}\,\,\,\,
\green{\downarrow}\,\,
& 
F_{2,0,1} & F_{1,1,1} & F_{0,2,1} & F_{3,0,0} & F_{0,3,0}
\\
\green{\bf 3a} & 
0 & 0 & 0 & 0 & 0
\\
\green{\bf 3b} & 
0 & 0 & 0 & 1 & 0
\\
\green{\bf 3c} & 
1 & 0 & 1 & 0 & 0
\\
\green{\bf 3d} & 
1 & 0 & 1 & 1 & 0
\end{array}
\]

In case {\footnotesize\sf C2}, let us abbreviate:
\[
\alpha
\,:=\,
F_{2,0,1}-F_{0,2,1},
\ \ \ \ \ \ \ \ \ \ \ \ \ \ \ \ \ \ \ \ \ \ \ \ \ \
\beta
\,:=\,
F_{1,1,1},
\]
so that:
\[
\alpha^2
+
\beta^2
\,>\,
0.
\]

\begin{Assertion}
In case {\footnotesize\sf C2},
one can normalize $G_{1,1,1} := 0$.
\end{Assertion}

\proof
Recall:
\[
G_{1,1,1}
\,=\,
\frac{1}{(a_{1,1}^2+a_{2,1}^2)\,a_{3,3}}\,
\Big\{
2\,a_{1,1}\,a_{2,1}\,
\big(
F_{2,0,1}-F_{0,2,1}
\big)
+
\big(
a_{1,1}^2-a_{2,1}^2
\big)\,
F_{1,1,1}
\Big\}.
\]
To make $G_{1,1,1} := 0$, 
the full freedom of the group parameters $\big(a_{1,1}, a_{2,1},
a_{3,3} \big)$ is not needed. Hence, we set $a_{3,3} := 1$, and:
\[
a_{1,1} 
\,:=\,
c,
\ \ \ \ \ \ \ \ \ \ \ \ \ \ \ \ \ \ \ \
a_{2,1}
\,:=\,
s,
\ \ \ \ \ \ \ \ \ \ \ \ \ \ \ \ \ \ \ \
\text{with}
\ \ \ \ \ \ \ \ \ \ \ \ \ \ \ \ \ \ \ \
c^2+s^2
\,=\,
1,
\]
so that we should find a real $c$ with $-1 \leqslant c \leqslant 1$ 
which satisfies:
\[
\aligned
0
&
\overset{\text{\bf ?}}{\,=\,}
2\,c\,s\,\alpha
+
(c^2-s^2)\,\beta
\\
&
\,=\,
2\,c\,\sqrt{1-c^2}\,\alpha
+
(1-2\,c^2)\,\beta.
\endaligned
\]
The 4 solutions of this biquadratic equation in $c$ are:
\[
c
\,=\,
\pm\,\frac{1}{\sqrt{2}}\,
\sqrt{
\frac{
\alpha^2+\beta^2
\pm\sqrt{\alpha^4+\alpha^2\beta^2}}{\alpha^2+\beta^2}}.
\]
Indeed, these values of $c$ belong to $[-1, 1]$:
\[
\alpha^2+\beta^2
\pm
\sqrt{\alpha^4+\alpha^2\beta^2}
\overset{\sf Yes}{\,\,\leqslant\,\,}
2\,\alpha^2
+
2\,\beta^2.
\qedhere
\]
\endproof

By Principle~{\ref{Principle-FG}}, 
we may assume $G_{1,1,1} = 0 = F_{1,1,1}$.

\begin{Assertion}
In case {\footnotesize\sf C2},
the associated group reduction reads:
\[
a_{2,1}
\,:=\,
0.
\]
\end{Assertion}

\proof
Indeed:
\[
0\,\,
{\overset{\green{\bf 1,1,1}}{\,\,=\,\,}}\,\,
\frac{2\,a_{1,1}\,a_{2,1}\,(F_{2,0,1}-F_{0,2,1})}{
(a_{1,1}^2+a_{2,1}^2)\,a_{3,3}},
\]
and since the set of case~{\footnotesize\sf C2}:
\[
\big\{
F_{1,1,1}\neq 0\,\,
\text{or}\,\,
F_{2,0,1}\neq F_{0,2,1}
\big\},
\]
is invariant, we get $0 = a_{1,1}\,a_{2,1}$, 
whence $0 = a_{2,1}$ in a neighborhood of the identity.
\endproof

The 5D group representation then reduces to 4D:
\[
\left(\!
\begin{array}{c}
G_{2,0,1}
\\
G_{0,2,1}
\\
G_{3,0,0}
\\
G_{0,3,0}
\end{array}
\!\right)
\,=\,
\left(\!
\def\arraystretch{1.25}
\begin{array}{cccc}
\frac{1}{a_{3,3}} & 
0 & 
0 & 
0
\\
0 &
\frac{1}{a_{3,3}} & 
0 & 
0
\\
\frac{-a_{3,1}}{a_{1,1}a_{3,3}} & 
\frac{a_{3,1}}{a_{1,1}a_{3,3}} & 
\frac{1}{a_{1,1}} & 
0 
\\
\frac{a_{3,2}}{a_{1,1}a_{3,3}} &
\frac{-a_{3,2}}{a_{1,1}a_{3,3}} & 
0 &
\frac{1}{a_{1,1}}
\end{array}
\!\right)\,
\left(\!
\begin{array}{c}
F_{2,0,1}
\\
F_{0,2,1}
\\
F_{3,0,0}
\\
F_{0,3,0}
\end{array}
\!\right),
\]
and there are 2 non-intersecting branches:

\[
\def\arraystretch{1.25}
\begin{array}{rcccc}
\green{\bf 2b}\,\,\,\,
\green{\downarrow}\,\,
& 
F_{2,0,1} & F_{0,2,1} & F_{3,0,0} & F_{0,3,0}
\\
\green{\bf 3e} & 
0 & 1 & 0 & 0
\\
\green{\bf 3f} & 
1 & F_{0,2,1}^{\neq1} & 0 & 0
\end{array}
\]
Here, in \green{\bf 2b3f}, the absolute invariant
$F_{0,2,1} = G_{0,2,1}$ is assumed to be $\neq 1$,
consistent with case {\footnotesize\sf C2}.

\medskip

In conclusion, in branch~\green{\bf 2b}, 
6 non-intersecting branches exist:
\[
\def\arraystretch{1.25}
\begin{array}{rccccc}
\green{\bf 2b}\,\,\,\,
\green{\downarrow}\,\,
& 
F_{2,0,1} & F_{1,1,1} & F_{0,2,1} & F_{3,0,0} & F_{0,3,0}
\\
\green{\bf 3a} & 
0 & 0 & 0 & 0 & 0
\\
\green{\bf 3b} & 
0 & 0 & 0 & 1 & 0
\\
\green{\bf 3c} & 
1 & 0 & 1 & 0 & 0
\\
\green{\bf 3d} & 
1 & 0 & 1 & 1 & 0
\\
\green{\bf 3e} & 
0 & 0 & 1 & 0 & 0
\\
\green{\bf 3f} & 
1 & 0 & F_{0,2,1}^{\neq1} & 0 & 0
\end{array}
\]

Because it would be (much) too lengthy to expose
all calculations in all (sub)branches,
we will content ourselves with presenting only
a few selected aspects. 

\SectionHead{Branch~\green{\bf 2b3e}}
{branch-2b3e}

After having normalized the two 
(pairs of) coefficients~{\eqref{normalize-F210-F120}},
in Branch~{\green{\bf 2b}}, we have:
\[
u
\,=\,
F
\,=\,
x^2+y^2+y^2z
+
{\rm O}_{x,y,z}(4)
\ \ \ \ \ \ \ \ \ \ \ \ \ \ \ \ \ \ \ \
\text{and}
\ \ \ \ \ \ \ \ \ \ \ \ \ \ \ \ \ \ \ \
v
\,=\,
G
\,=\,
p^2+q^2+q^2r
+
{\rm O}_{p,q,r}(4),
\]
with the equations of $\eqFG$ at order $4$ being:
\[
\aligned
0
&
{\overset{\green{\bf 3,0,0}}{\,\,\,=\,\,\,}}
-\,a_{1,1}^2\,a_{3,1}
-
2\,a_{1,1}\,a_{2,1}\,a_{3,2}
+
a_{2,1}^2\,a_{3,1},
\\
0
&
{\overset{\green{\bf 2,0,1}}{\,\,\,=\,\,\,}}
a_{2,1}^2\,a_{3,3},
\\
0
&
{\overset{\green{\bf 1,1,1}}{\,\,\,=\,\,\,}}
2\,a_{1,1}\,a_{2,1}\,a_{3,3},
\\
0
&
{\overset{\green{\bf 0,3,0}}{\,\,\,=\,\,\,}}
a_{1,1}^2\,a_{3,2}
-
2\,a_{1,1}\,a_{2,1}\,a_{3,1}
-
a_{2,1}^2\,a_{3,2},
\\
0
&
{\overset{\green{\bf 0,2,1}}{\,\,\,=\,\,\,}}
a_{1,1}^2\,a_{3,3}
-
a_{1,1}^2
-
a_{2,1}^2.
\endaligned
\]

The group reduction to stabilize such an
order 3 normal form, {\em i.e.} to fulfill these equations, is:
\[
a_{2,1}
\,:=\,
0,
\ \ \ \ \ \ \ \ \ \ \ \ \ \ \ \ \ \ \ \
a_{3,1}
\,:=\,
0,
\ \ \ \ \ \ \ \ \ \ \ \ \ \ \ \ \ \ \ \
a_{3,2}
\,:=\,
0,
\ \ \ \ \ \ \ \ \ \ \ \ \ \ \ \ \ \ \ \
a_{3,3}
\,=\,
1,
\]
with reduced matrix:
\[
\left[\,
\begin{array}{cccc}
a_{1,1} & \red{\bf 0} & \red{\bf 0} & \red{\bf 0}
\\
\red{\bf 0} & \red{a_{1,1}} & \red{\bf 0} & \red{\bf 0}
\\
\red{\bf 0} & \red{\bf 0} & \red{\bf 1} & b_3
\\
\red{\bf 0} & \red{\bf 0} & \red{\bf 0} & \red{a_{1,1}^2}
\end{array}
\!\right].
\]

Then at order 3, the
tangency equations (not written here)
$\eqLF$ 
require to solve similarly the corresponding
infinitesimal group parameters as:
\[
\aligned
A_{2,1}
&
\,:=\,
-\,F_{2,1,1}\,T_1-F_{1,2,1}\,T_2-F_{1,1,2}\,T_3,
\\
A_{3,1}
&
\,:=\,
\big(-\,2\,F_{2,2,0}+4\,F_{4,0,0}\big)\,T_1
+
\big(-\,3\,F_{1,3,0}+F_{3,1,0}\big)\,T_2
+
\big(-\,F_{1,2,1}+F_{3,0,1}\big)\,T_3,
\\
A_{3,2}
&
\,:=\,
\big(-\,F_{1,3,0}+3\,F_{3,1,0}\big)\,T_1
+
\big(-\,4\,F_{0,4,0}+2\,F_{2,2,0}\big)\,T_2
+
\big(-\,F_{0,3,1}+F_{2,1,1}\big)\,T_3,
\\
A_{3,3}
&
\,:=\,
-\,F_{1,2,1}\,T_1
-
3\,F_{0,3,1}\,T_2
+
\big(-\,2\,F_{0,2,1}+1\big)\,T_3.
\endaligned
\]

It then remains:
\[
\aligned
0
&
{\overset{\green{\bf 2,0,1}}{\,\,\,=\,\,\,}}
3\,F_{3,0,1}\,T_1
+
F_{2,1,1}\,T_2
+
2\,F_{2,0,2}\,T_3,
\\
0
&
{\overset{\green{\bf 1,0,2}}{\,\,\,=\,\,\,}}
2\,F_{2,0,2}\,T_1
+
F_{1,1,2}\,T_2
+
3\,F_{1,0,3}\,T_3,
\\
0
&
{\overset{\green{\bf 0,1,2}}{\,\,\,=\,\,\,}}
F_{1,1,2}\,T_1
+
\big(2\,F_{0,2,2}-2\big)\,T_2
+
3\,F_{1,0,3}\,T_3,
\\
0
&
{\overset{\green{\bf 0,0,3}}{\,\,\,=\,\,\,}}
F_{1,0,3}\,T_1
+
F_{0,1,3}\,T_2
+
4\,F_{0,0,4}\,T_4,
\endaligned
\]
which forces:
\[
\aligned
F_{3,0,1}
&
\,:=\,
0,
&
\ \ \ \ \ \ \
F_{2,1,1}
&
\,:=\,
0,
&
\ \ \ \ \ \ \
F_{2,0,2}
&
\,:=\,
0,
&
\ \ \ \ \ \ \
F_{1,1,2}
&
\,:=\,
0,
\\
F_{1,0,3}
&
\,:=\,
0,
&
\ \ \ \ \ \ \
F_{0,2,2}
&
\,:=\,
1,
&
\ \ \ \ \ \ \
F_{0,1,3}
&
\,:=\,
0
&
\ \ \ \ \ \ \
F_{0,0,4}
&
\,:=\,
0.
\endaligned
\]

Once these constraints are also set for the corresponding
$G$-coefficients in the target space:
\[
\aligned
G_{3,0,1}
&
\,:=\,
0,
&
\ \ \ \ \ \ \
G_{2,1,1}
&
\,:=\,
0,
&
\ \ \ \ \ \ \
G_{2,0,2}
&
\,:=\,
0,
&
\ \ \ \ \ \ \
G_{1,1,2}
&
\,:=\,
0,
\\
G_{1,0,3}
&
\,:=\,
0,
&
\ \ \ \ \ \ \
G_{0,2,2}
&
\,:=\,
1,
&
\ \ \ \ \ \ \
G_{0,1,3}
&
\,:=\,
0
&
\ \ \ \ \ \ \
G_{0,0,4}
&
\,:=\,
0,
\endaligned
\]
at order 4, one of the equations $\eqFG_{i,j,k}$ is:
\[
0
{\overset{\green{\bf 2,2,0}}{\,\,\,=\,\,\,}}
a_{1,1}^2\,
\Big(
-\,F_{2,2,0}
+
a_{1,1}^2\,G_{2,2,0}
+
b_3
\Big),
\]
so using $b_3$, we normalize:
\[
G_{2,2,0}
\,:=\,
0
\,=:\,
F_{2,2,0},
\]
and then to stabilize this normalization, it is necessary to set:
\[
b_3
\,:=\,
0.
\]

The affine group matrix is reduced to:
\[
\left[\,
\begin{array}{cccc}
a_{1,1} & \red{\bf 0} & \red{\bf 0} & \red{\bf 0}
\\
\red{\bf 0} & \red{a_{1,1}} & \red{\bf 0} & \red{\bf 0}
\\
\red{\bf 0} & \red{\bf 0} & \red{\bf 1} & \red{\bf 0}
\\
\red{\bf 0} & \red{\bf 0} & \red{\bf 0} & \red{a_{1,1}^2}
\end{array}
\!\right],
\]
and the 6 remaining order 4 coefficients in: 
\[
\aligned
u
&
\,=\,
x^2+y^2
+
y^2z
\\
&
\ \ \ \ \
+
F_{4,0,0}\,x^4
+
F_{3,1,0}\,x^3y
+
F_{1,3,0}\,xy^3
+
F_{1,2,1}\,xy^2z
+
F_{0,4,0}\,y^4
+
F_{0,3,1}\,y^3z
+
y^2z^2
+
{\rm O}_{x,y,z}(5),
\\
v
&
\,=\,
p^2+q^2
+
q^2r
\\
&
\ \ \ \ \
+
G_{4,0,0}\,p^4
+
G_{3,1,0}\,p^3q
+
G_{1,3,0}\,pq^3
+
G_{1,2,1}\,pq^2r
+
G_{0,4,0}\,q^4
+
G_{0,3,1}\,q^3r
+
q^2r^2
+
{\rm O}_{p,q,r}(5),
\endaligned
\]
are all relative invariants, 
organized as the linear representation
$6 \times 6$ matrix:
\[
\left(\!
\begin{array}{c}
G_{4,0,0}
\\
G_{3,1,0}
\\
G_{1,3,0}
\\
G_{1,2,1}
\\
G_{0,4,0}
\\
G_{0,3,1}
\end{array}
\!\right)
\,=\,
\left(\!
\def\arraystretch{1.25}
\begin{array}{cccccc}
\tfrac{1}{a_{1,1}^2} & 0 & 0 & 0 & 0 & 0
\\
0 & \tfrac{1}{a_{1,1}^2} & 0 & 0 & 0 & 0
\\
0 & 0 & \tfrac{1}{a_{1,1}^2} & 0 & 0 & 0
\\
0 & 0 & 0 & \tfrac{1}{a_{1,1}} & 0 & 0 
\\
0 & 0 & 0 & 0 & \tfrac{1}{a_{1,1}^2} & 0
\\
0 & 0 & 0 & 0 & 0 & \tfrac{1}{a_{1,1}} 
\end{array}
\!\right)\,
\left(\!
\begin{array}{c}
F_{4,0,0}
\\
F_{3,1,0}
\\
F_{1,3,0}
\\
F_{1,2,1}
\\
F_{0,4,0}
\\
F_{0,3,1}
\end{array}
\!\right).
\]

At order 4, the process of creation of non-intersecting branches
leads us to the 11 branches:
\[
\def\arraystretch{1.25}
\begin{array}{rcccccc}
\green{\bf 2b3e}\,\,\,\,
\green{\downarrow}\,\,
& F_{4,0,0} & F_{3,1,0} & F_{1,3,0} & F_{1,2,1} & F_{0,4,0} & F_{0,3,1}
\\
\green{\bf 4a} & 
0 & 0 & 0 & 0 & 0 & 0 
\\
\green{\bf 4b} & 
0 & 0 & 0 & 0 & 0 & 1
\\
\green{\bf 4c} & 
0 & 0 & 0 & 0 & 1 & F_{0,3,1}
\\
\green{\bf 4d} & 
0 & 0 & 0 & 0 & -1 & F_{0,3,1}
\\
\green{\bf 4e} & 
0 & 0 & 0 & 1 & F_{0,4,0} & F_{0,3,1}
\\
\green{\bf 4f} & 
0 & 0 & 1 & F_{1,2,1} & F_{0,4,0} & F_{0,3,1}
\\
\green{\bf 4g} & 
0 & 0 & -1 & F_{1,2,1} & F_{0,4,0} & F_{0,3,1}
\\
\green{\bf 4h} & 
0 & 1 & F_{1,3,0} & F_{1,2,1} & F_{0,4,0} & F_{0,3,1}
\\
\green{\bf 4i} & 
0 & -1 & F_{1,3,0} & F_{1,2,1} & F_{0,4,0} & F_{0,3,1}
\\
\green{\bf 4j} & 
1 & F_{3,1,0} & F_{1,3,0} & F_{1,2,1} & F_{0,4,0} & F_{0,3,1}
\\
\green{\bf 4k} & 
-1 & F_{3,1,0} & F_{1,3,0} & F_{1,2,1} & F_{0,4,0} & F_{0,3,1}
\end{array}
\]

Before entering branches,
to the normalization 
$F_{2,2,0} := 0$ by means of $b_3$, 
then at order 4 in the equation 
${\overset{\green{\bf 2,2,0}}{\,\,\,=\,\,\,}}$ of $\eqLF$, 
there corresponds the infinitesimal assignment:
\[
\aligned
B_3
&
\,:=\,
\big(
-\,4\,F_{1,2,1}\,F_{4,0,0}
-
3\,F_{3,2,0}
\big)\,T_1
+
\big(
6\,F_{1,2,1}\,F_{1,3,0}
-
4\,F_{1,2,1}\,F_{3,1,0}
-
3\,F_{2,3,0}
\big)\,T_2
\\
&
\ \ \ \ \
+
\big(
F_{1,2,1}^2
-
F_{2,2,1}
\big)\,T_3.
\endaligned
\]

\SectionHead{Branch~\green{\bf 2b3e4a}}
{branch-2b-3e-4a}

Let us treat the branch~{\green{\bf 2b3e4a}}, 
for which all 6 relative invariants vanish:
\[
\aligned
F_{4,0,0}
&
\,:=\,
0,
&
\ \ \ \ \ \ \
F_{3,1,0}
&
\,:=\,
0,
&
\ \ \ \ \ \ \
F_{1,3,0}
&
\,:=\,
0,
&
\ \ \ \ \ \ \
F_{1,2,1}
&
\,:=\,
0,
\ \ \ \ \ \ \
F_{0,4,0}
&
\,:=\,
0,
&
\ \ \ \ \ \ \
F_{0,3,1}
&
\,:=\,
0,
\\
G_{4,0,0}
&
\,:=\,
0,
&
\ \ \ \ \ \ \
G_{3,1,0}
&
\,:=\,
0,
&
\ \ \ \ \ \ \
G_{1,3,0}
&
\,:=\,
0,
&
\ \ \ \ \ \ \
G_{1,2,1}
&
\,:=\,
0,
\ \ \ \ \ \ \
G_{0,4,0}
&
\,:=\,
0,
&
\ \ \ \ \ \ \
G_{0,3,1}
&
\,:=\,
0.
\endaligned
\]

Then the equations of $\eqLF$:
\[
\aligned
0
&
{\overset{\green{\bf 4,0,0}}{\,\,\,=\,\,\,}}
5\,F_{5,0,0}\,T_1+F_{4,1,0}\,T_2+F_{4,0,1}\,T_3,
\\
0
&
{\overset{\green{\bf 3,1,0}}{\,\,\,=\,\,\,}}
4\,F_{4,1,0}\,T_1+2\,F_{3,2,0}\,T_2+F_{3,1,1}\,T_3,
\\
0
&
{\overset{\green{\bf 1,3,0}}{\,\,\,=\,\,\,}}
2\,F_{2,3,0}\,T_1+4\,F_{1,4,0}\,T_2+F_{1,3,1}\,T_3,
\\
0
&
{\overset{\green{\bf 0,4,0}}{\,\,\,=\,\,\,}}
\big(F_{1,4,0}-3\,F_{3,2,0}\big)\,T_1
+
\big(5\,F_{0,5,0}-3\,F_{2,3,0}\big)\,T_2
+
\big(F_{0,4,1}-F_{2,2,1}\big)\,T_3,
\\
0
&
{\overset{\green{\bf 3,0,1}}{\,\,\,=\,\,\,}}
4\,F_{4,0,1}\,T_1+F_{3,1,1}\,T_2+2\,F_{3,0,2}\,T_3,
\\
0
&
{\overset{\green{\bf 2,1,1}}{\,\,\,=\,\,\,}}
3\,F_{3,1,1}\,T_1+2\,F_{2,2,1}\,T_2+2\,F_{2,1,2}\,T_3,
\\
0
&
{\overset{\green{\bf 1,2,1}}{\,\,\,=\,\,\,}}
2\,F_{2,2,1}\,T_1+3\,F_{1,3,1}\,T_2+2\,F_{1,2,2}\,T_3,
\\
0
&
{\overset{\green{\bf 0,3,1}}{\,\,\,=\,\,\,}}
F_{1,3,1}\,T_1+4\,F_{0,4,1}\,T_2+2\,F_{0,3,2}\,T_3,
\\
0
&
{\overset{\green{\bf 2,0,2}}{\,\,\,=\,\,\,}}
3\,F_{3,0,2}\,T_1+F_{2,1,2}\,T_2+3\,F_{2,0,3}\,T_3,
\\
0
&
{\overset{\green{\bf 1,1,2}}{\,\,\,=\,\,\,}}
2\,F_{2,1,2}\,T_1+2\,F_{1,2,2}\,T_2+3\,F_{1,1,3}\,T_3,
\\
0
&
{\overset{\green{\bf 0,2,2}}{\,\,\,=\,\,\,}}
F_{1,2,2}\,T_1+3\,F_{0,3,2}\,T_2+\big(3\,F_{0,2,3}-3\big)\,T_3
\\
0
&
{\overset{\green{\bf 1,0,3}}{\,\,\,=\,\,\,}}
2\,F_{2,0,3}\,T_1+F_{1,1,3}\,T_2+4\,F_{1,0,4}\,T_3,
\\
0
&
{\overset{\green{\bf 0,1,3}}{\,\,\,=\,\,\,}}
F_{1,1,3}\,T_1+\big(2\,F_{0,2,3}-2\big)\,T_2+4\,F_{0,1,4}\,T_3,
\\
0
&
{\overset{\green{\bf 0,0,4}}{\,\,\,=\,\,\,}}
F_{1,0,4}\,T_1+F_{0,1,4}\,T_2+5\,F_{0,0,5}\,T_3,
\endaligned
\]
force:
\[
\aligned
{}&
\aligned
F_{5,0,0}
&
\,:=\,
0,
&
\ \ \ \ \ \ \
F_{4,1,0}
&
\,:=\,
0,
&
\ \ \ \ \ \ \
F_{4,0,1}
&
\,:=\,
0,
&
\ \ \ \ \ \ \
F_{3,2,0}
&
\,:=\,
0,
\ \ \ \ \ \ \
F_{3,1,1}
&
\,:=\,
0,
\\
F_{3,0,2}
&
\,:=\,
0,
&
\ \ \ \ \ \ \
F_{2,3,0}
&
\,:=\,
0,
&
\ \ \ \ \ \ \
F_{2,2,1}
&
\,:=\,
0,
&
\ \ \ \ \ \ \
F_{2,1,2}
&
\,:=\,
0,
\ \ \ \ \ \ \
F_{2,0,3}
&
\,:=\,
0,
\\
F_{1,3,1}
&
\,:=\,
0,
&
\ \ \ \ \ \ \
F_{1,2,2}
&
\,:=\,
0,
&
\ \ \ \ \ \ \
F_{1,1,3}
&
\,:=\,
0,
&
\ \ \ \ \ \ \
F_{1,0,4}
&
\,:=\,
0,
\ \ \ \ \ \ \
F_{0,5,0}
&
\,:=\,
0,
\\
F_{0,4,1}
&
\,:=\,
0,
&
\ \ \ \ \ \ \
F_{0,3,2}
&
\,:=\,
0,
&
\ \ \ \ \ \ \
F_{0,2,3}
&
\,:=\,
0,
&
\ \ \ \ \ \ \
F_{0,1,4}
&
\,:=\,
0,
\ \ \ \ \ \ \
F_{0,0,5}
&
\,:=\,
0,
\endaligned
\\
{}&
F_{1,4,0}
\,:=\,
3\,F_{3,2,0}.
\endaligned
\]

Importing in $\eqFG$ 
these equations and the same for the $G$-coefficients,
it then happens that at order 5, no more branching is required.

Pursuing up to order 6, the process stabilizes, and we receive:
\[
u
\,=\,
F
\,=\,
x^2+y^2+y^2z
+
y^2z^2
+
y^3z^3
+
y^2z^4
+
y^2z^5
+
{\rm O}_{x,y,z}(8),
\]
with a 4D transitive
Lie algebra of infinitesimal affine vector fields:
\[
\aligned
e_1
&
\,:=\,
\partial_x
+
2\,x\,\partial_u,
\\
e_2
&
\,:=\,
\big(1-z\big)\,\partial_y
+
2\,y\,\partial_u,
\\
e_3
&
\,=\,
-\,\tfrac{1}{2}\,y\,\partial_y
+
\big(1-z\big)\,\partial_z,
\\
e_4
&
\,:=\,
x\,\partial_x
+
y\,\partial_y
+
2\,u\,\partial_u,
\endaligned
\]
having (solvable) Lie structure:
\[
\footnotesize
\def\arraystretch{1.25}
\begin{array}{c|cccc}
{} & e_1 & e_2 & e_3 & e_4
\\
\hline
e_1 & 
0 & 0 & 0 & e_1
\\
e_2 &
0 & 0 & \tfrac{1}{2}e_2 & e_2 
\\
e_3 &
0 & -\tfrac{1}{2}e_2 & 0 & 0
\\
e_4 &
-e_1 & -e_2 & 0 & 0
\end{array}
\]

\SectionHead{Shortcut for Subbranches of Branch~\green{\bf 2b3e}}
{shortcut-branch-2b-3e}

Instead of exploring all the remaining $10 = 11 - 1$ 
branches~\green{\bf 2b3e4b} up to~\green{\bf 2b3e4k}
and all their potential subbranches,
let us proceed as follows.

At order 4, the equations of $\eqLF$ which
involve {\em only} $T_1$, $T_2$, $T_3$ are:
\[
\aligned
0
&
{\overset{\green{\bf 3,0,1}}{\,\,\,=\,\,\,}}
4\,F_{4,0,1}\,T_1
+
\big(-\,F_{3,1,0}+F_{3,1,1}\big)\,T_2
+
2\,F_{3,0,2}\,T_3,
\\
0
&
{\overset{\green{\bf 2,1,1}}{\,\,\,=\,\,\,}}
\big(-\,3\,F_{3,1,0}+3\,F_{3,1,1}\big)\,T_1
+
\big(-\,2\,F_{1,2,1}^2+2\,F_{2,2,1}\big)\,T_2
+
2\,F_{2,1,2}\,T_3,
\\
0
&
{\overset{\green{\bf 2,0,2}}{\,\,\,=\,\,\,}}
3\,F_{3,0,2}\,T_1+F_{2,1,2}\,T_2+3\,F_{2,0,3}\,T_3,
\\
0
&
{\overset{\green{\bf 1,1,2}}{\,\,\,=\,\,\,}}
2\,F_{2,1,2}\,T_1
+
\big(-\,4\,F_{1,2,1}+2\,F_{1,2,2}\big)\,T_2
+
F_{1,1,3}\,T_3,
\\
0
&
{\overset{\green{\bf 0,2,2}}{\,\,\,=\,\,\,}}
\big(-\,2\,F_{1,2,1}+F_{1,2,2}\big)\,T_1
+
\big(-\,9\,F_{0,3,1}+3\,F_{0,3,2}\big)\,T_2
+
\big(3\,F_{0,2,3}-3\big)\,T_3,
\\
0
&
{\overset{\green{\bf 1,0,3}}{\,\,\,=\,\,\,}}
2\,F_{2,0,3}\,T_1+F_{1,1,3}\,T_2+4\,F_{1,0,4}\,T_3,
\\
0
&
{\overset{\green{\bf 0,1,3}}{\,\,\,=\,\,\,}}
F_{1,1,3}\,T_1
+
\big(2\,F_{0,2,3}-2\big)\,T_2
+
F_{1,0,4}\,T_3,
\\
0
&
{\overset{\green{\bf 0,0,4}}{\,\,\,=\,\,\,}}
F_{1,0,4}\,T_1+F_{0,1,4}\,T_2+5\,F_{0,0,5}\,T_3,
\endaligned
\]
while the remaining ones are:
\[
\aligned
0
&
{\overset{\green{\bf 4,0,0}}{\,\,\,=\,\,\,}}
2\,F_{4,0,0}\,A_{1,1}+\cdots,
\\
0
&
{\overset{\green{\bf 3,1,0}}{\,\,\,=\,\,\,}}
2\,F_{3,1,0}\,A_{1,1}+\cdots,
\\
0
&
{\overset{\green{\bf 1,3,0}}{\,\,\,=\,\,\,}}
2\,F_{1,3,0}\,A_{1,1}+\cdots,
\\
0
&
{\overset{\green{\bf 1,2,1}}{\,\,\,=\,\,\,}}
F_{1,2,1}\,A_{1,1}+\cdots,
\\
0
&
{\overset{\green{\bf 0,4,0}}{\,\,\,=\,\,\,}}
2\,F_{0,4,0}\,A_{1,1}+\cdots,
\\
0
&
{\overset{\green{\bf 0,3,1}}{\,\,\,=\,\,\,}}
F_{0,3,1}\,A_{1,1}+\cdots.
\endaligned
\]

Therefore, in all branches:
\[
\aligned
F_{4,0,1}
&
\,:=\,
0,
&
\ \ \ \ \ \ \
F_{3,1,1}
&
\,:=\,
F_{3,1,0},
&
\ \ \ \ \ \ \
F_{3,0,2}
&
\,:=\,
0,
&
\ \ \ \ \ \ \
F_{2,2,1}
&
\,:=\,
F_{1,2,1}^2,
&
\ \ \ \ \ \ \
F_{2,1,2}
&
\,:=\,
0,
\\
F_{2,0,3}
&
\,:=\,
0,
&
\ \ \ \ \ \ \
F_{1,2,2}
&
\,:=\,
2\,F_{1,2,1},
&
\ \ \ \ \ \ \
F_{1,1,3}
&
\,:=\,
0,
&
\ \ \ \ \ \ \
F_{1,0,4}
&
\,:=\,
0,
&
\ \ \ \ \ \ \
F_{0,3,2}
&
\,:=\,
3\,F_{0,3,1},
\\
F_{0,2,3}
&
\,:=\,
1,
&
\ \ \ \ \ \ \
F_{0,1,4}
&
\,:=\,
0,
&
\ \ \ \ \ \ \
F_{0,0,5}
&
\,:=\,
0.
&
\ \ \ \ \ \ \
&
&
&
\endaligned
\]
and it remains:

\[
\aligned
0
&
{\overset{\green{\bf 4,0,0}}{\,\,\,=\,\,\,}} 
2\,F_{4,0,0}\,A_{1,1}
+
5\,F_{5,0,0}\,T_1
+
\big(-\,F_{1,2,1}\,F_{3,1,0}+F_{4,1,0}\big)\,T_2,
\\
0
&
{\overset{\green{\bf 3,1,0}}{\,\,\,=\,\,\,}}
2\,F_{3,1,0}\,A_{1,1}
+
4\,F_{4,1,0}\,T_1
+
\big(4\,F_{1,2,1}\,F_{4,0,0}+2\,F_{3,2,0}\big)\,T_2
+
\tfrac{1}{2}\,F_{3,1,0}\,T_3,
\\
0
&
{\overset{\green{\bf 1,3,0}}{\,\,\,=\,\,\,}}
2\,F_{1,3,0}\,A_{1,1}
+
\big(4\,F_{0,3,1}\,F_{4,0,0}-F_{1,2,1}\,F_{1,3,0}
+3\,F_{1,2,1}\,F_{3,1,0}+2\,F_{2,3,0}\big)\,T_1
\\
&
\ \ \ \ \
+
\big(-\,3\,F_{0,3,1}\,F_{1,3,0}+F_{0,3,1}\,F_{3,1,0}
-8\,F_{0,4,0}\,F_{1,2,1}+4\,F_{1,4,0}\big)\,T_2
\\
&
\ \ \ \ \
+
\big(-\,\tfrac{3}{2}\,F_{1,3,0}-2\,F_{0,3,1}\,F_{1,2,1}
+F_{1,3,1}\big)\,T_3,
\\
0
&
{\overset{\green{\bf 1,2,1}}{\,\,\,=\,\,\,}}
F_{1,2,1}\,A_{1,1}
+
\big(F_{1,2,1}^2+4\,F_{4,0,0}\big)\,T_1
+
\big(-\,6\,F_{0,3,1}\,F_{1,2,1}-9\,F_{1,3,0}
+3\,F_{1,3,1}+F_{3,1,0}\big)\,T_2,
\\
0
&
{\overset{\green{\bf 0,4,0}}{\,\,\,=\,\,\,}}
2\,A_{1,1}\,F_{0,4,0}
+
\big(-\,F_{0,3,1}\,F_{1,3,0}+3\,F_{0,3,1}\,F_{3,1,0}
-4\,F_{1,2,1}\,F_{4,0,0}+F_{1,4,0}-3\,F_{3,2,0}\big)\,T_1
\\
&
\ \ \ \ \
+
\big(-\,4\,F_{0,3,1}\,F_{0,4,0}+7\,F_{1,2,1}\,F_{1,3,0}
-4\,F_{1,2,1}\,F_{3,1,0}+5\,F_{0,5,0}-3\,F_{2,3,0}\big)\,T_2
\\
&
\ \ \ \ \
+
\big(-\,F_{0,3,1}^2-2\,F_{0,4,0}+F_{0,4,1}\big)\,T_3,
\\
0
&
{\overset{\green{\bf 0,3,1}}{\,\,\,=\,\,\,}}
F_{0,3,1}\,A_{1,1}
+
\big(-\,F_{0,3,1}\,F_{1,2,1}-2\,F_{1,3,0}+F_{1,3,1}\big)\,T_1
\\
&
\ \ \ \ \
+
\big(-\,3\,F_{0,3,1}^2+F_{1,2,1}^2-12\,F_{0,4,0}
+4\,F_{0,4,1}\big)\,T_2
+
\tfrac{3}{2}\,F_{0,3,1}\,T_3.
\endaligned
\]

Since the current reduced subgroup matrix is 1-dimensional:
\[
\left[\,
\begin{array}{cccc}
a_{1,1} & \red{\bf 0} & \red{\bf 0} & \red{\bf 0}
\\
\red{\bf 0} & \red{a_{1,1}} & \red{\bf 0} & \red{\bf 0}
\\
\red{\bf 0} & \red{\bf 0} & \red{\bf 1} & \red{\bf 0}
\\
\red{\bf 0} & \red{\bf 0} & \red{\bf 0} & \red{a_{1,1}^2}
\end{array}
\!\right],
\]
parametrized by the single group parameter $a_{1,1}$,
then similarly, at the infinitesimal level, 
together with the 3 transitivity parameters $T_1$, $T_2$, $T_3$, 
only the single 
isotropy parameter
$A_{1,1}$ is still free in the infinitesimal
transformation $L$ of
Section~{\ref{infinitesimal-affine-transformations}},
and in fact, at this stage of the process:
\[
L
\,=\,
X\,\partial_x
+
Y\,\partial_y
+
Z\,\partial_z
+
U\,\partial_u,
\]
with:
\[
\aligned
X
&
\,=\,
x\,A_{1,1}
+
\big(1-2\,F_{4,0,0}\,u\big)\,T_1
+
\big(
F_{1,2,1}\,y
-
\tfrac{1}{2}\,F_{3,1,0}\,u
\big)\,
T_2,
\\
Y
&
\,=\,
y\,A_{1,1}
+
\big(-\,\tfrac{3}{2}\,F_{3,1,0}\,u\big)\,T_1
+
\big(
1-F_{1,2,1}\,x-z
\big)\,T_2
-
\tfrac{1}{2}\,y\,T_3,
\\
Z
&
\,=\,
\big(
4\,F_{4,0,0}\,x
-
F_{1,3,0}\,y
+
3\,F_{3,1,0}\,y
-
F_{1,2,1}\,z
-
3\,F_{3,2,0}\,u
-
4\,F_{1,2,1}\,F_{4,0,0}\,u
\big)\,T_1
\\
&
\ \ \ \ \
+
\big(
F_{3,1,0}\,x
-
\,3\,F_{1,3,0}\,x
-
4\,F_{0,4,0}\,y
-
3\,F_{0,3,1}\,z
-
3\,F_{2,3,0}\,u
-
4\,F_{1,2,1}\,F_{1,3,0}\,u
+
6\,F_{1,2,1}\,F_{1,3,0}\,u
\big)\,T_2
\\
&
\ \ \ \ \
+
\big(
1-F_{1,2,1}\,x-F_{0,3,1}\,y-z
\big)\,T_3,
\\
U
&
\,=\,
2\,u\,A_{1,1}
+
2\,x\,T_1
+
2\,y\,T_2.
\endaligned
\]

The main idea is to realize that, 
in all possible remaining 10 
branches~\green{\bf 2b3e4b} up to~\green{\bf 2b3e4k},
where $a_{1,1}$ would be normalized in some way,
always, at the infinitesimal level,
the reduction would be of the linear form:
\[
A_{1,1}
\,:=\,
\Lambda\,T_1
+
\Phi\,T_2
+
\Psi\,T_3,
\]
for some expressions $\Lambda$, $\Phi$, $\Psi$,
depending on the branches, leading to some
{\em simply transitive} Lie algebra of affine vector fields,
generated by:
\[
\aligned
e_1
&
\,:=\,
\big(
1
+
\Lambda\,x
-
2\,F_{4,0,0}\,u
\big)\,
\partial_x
+
\big(
y\,\Lambda
-
\tfrac{3}{2}\,
F_{3,1,0}\,
u
\big)\,\partial_y
\\
&
\ \ \ \ \
+
\big(
4\,F_{4,0,0}\,x
+
3\,F_{3,1,0}\,y
-
F_{1,3,0}\,y
-
F_{1,2,1}\,z
-
3\,F_{3,2,0}\,u
-
4\,F_{1,2,1}\,F_{4,0,0}\,u
\big)\,\partial_z
+
\big(
2\,x+2\,\Lambda\,u
\big)\,\partial_u,
\\
e_2
&
\,:=\,
\big(
\Phi\,x
+
F_{1,2,1}\,y
-
\tfrac{1}{2}\,F_{3,1,0}\,u
\big)\,
\partial_x
+
\big(
1
-
F_{1,2,1}\,x
+
\Phi\,y
-
z
\big)\,\partial_y
\\
&
\ \ \ \ \
+
\big(
F_{3,1,0}\,x
-
3\,F_{1,3,0}\,x
-
4\,F_{0,4,0}\,y
-
3\,F_{0,3,1}\,z
-
3\,F_{2,3,0}\,u
-
4\,F_{1,2,1}\,F_{3,1,0}\,u
+
6\,F_{1,2,1}\,F_{1,3,0}\,u
\big)\,\partial_z
\\
&
\ \ \ \ \
+
\big(
2\,y+2\,\Phi\,u
\big)\,\partial_u,
\\
e_3
&
\,:=\,
\big(
\Psi\,x
\big)\,
\partial_x
+
\big(
-\,\tfrac{1}{2}\,y
+
\Psi\,y
\big)\,\partial_y
+
\big(
1
-
F_{1,2,1}\,x
-
F_{0,3,1}\,y
-
z
\big)\,\partial_z
+
\big(
2\,\Psi\,u
\big)\,\partial_u.
\endaligned
\]

These 3 vector fields on $\R_{x,y,z,u}^4$ should constitute
a 3D Lie algebra. The structure constants are uniquely
determined by the Lie brackets:
\[
\big[e_1,e_2\big],
\ \ \ \ \ \ \ \ \ \ \ \ \ \ \ \ \ \ \ \
\big[e_1,e_3\big],
\ \ \ \ \ \ \ \ \ \ \ \ \ \ \ \ \ \ \ \
\big[e_2,e_3\big],
\]
taken at the origin $(x,y,z,u) = (0,0,0,0)$, namely:
\[
\big[
e_i,\,e_j
\big]
\,=\,
\sum_{1\leqslant r\leqslant 3}\,
C_{i,j,r}\,e_r
\ \ \ \ \ \ \ \ \ \ \ \ \ \ \ \ \ \ \ \
{\scriptstyle{(1\,\leqslant\,i<j\leqslant\,3)}},
\]
and we find:
\[
\aligned
C_{1,2,1}
&
\,:=\,
\Phi,
\\
C_{1,2,2}
&
\,:=\,
-\,
F_{1,2,1}
-
\Lambda,
\\
C_{1,2,3}
&
\,:=\,
-\,2\,F_{1,3,0}
-
2\,F_{3,1,0},
\\
C_{1,3,1}
&
\,:=\,
\Psi,
\\
C_{1,3,2}
&
\,:=\,
0,
\\
C_{1,3,3}
&
\,:=\,
0,
\\
C_{2,3,1}
&
\,:=\,
0,
\\
C_{2,3,2}
&
\,:=\,
\tfrac{1}{2}
+
\Psi,
\\
C_{2,3,3}
&
\,:=\,
2\,F_{0,3,1}.
\endaligned
\]

Setting, for $1 \leqslant i < j \leqslant 3$
and for $1 \leqslant r \leqslant 3$:
\[
C_{i,i,r}
\,:=\,
3,
\ \ \ \ \ \ \ \ \ \ \ \ \ \ \ \ \ \ \ \ \ \ \ \ \ \
C_{j,i,r}
\,:=\,
-\,C_{i,j,r},
\]
we find that
all Jacobi identities are satisfied if and only if:
\[
\aligned
0
&
\,=\,
\Jac_1
\,:=\,
\Psi\,
\big(
F_{1,2,1}
+
\Lambda
\big),
\\
0
&
\,=\,
\Jac_2
\,:=\,
-\,\tfrac{1}{2}\,\Phi
-
\Phi\,\Psi
-
2\,F_{0,3,1}\,\Psi,
\\
0
&
\,=\,
\Jac_3
\,:=\,
F_{3,1,0}
+
F_{1,3,0}
+
4\,F_{3,1,0}\,\Psi
+
4\,F_{1,3,0}\,\Psi
-
2\,F_{0,3,1}\,\Lambda
-
2\,F_{0,3,1}\,F_{1,2,1}.
\endaligned
\]
However, these (only) three equations are not informative enough.

Therefore, we explore the requirement that:
\[
\big[
e_i,\,
e_j
\big]
-
\sum_{1\leqslant r\leqslant 3}\,
C_{i,j,r}\,e_r
\,=\,
0
+
{\rm O}_{x,y,z}(5),
\]
which gives us a lot more equations:
\[
\footnotesize
\aligned
Z_1
&
\,:=\,
\Psi\,
\big(
F_{1,2,1}
+
\Lambda
\big),
\\
Z_2
&
\,:=\,
\Psi\,
\big(\Lambda+F_{1,2,1}\big),
\\
Z_3
&
\,:=\,
-\,\Psi\,F_{1,2,1}
\\
Z_4
&
\,:=\,
-\,2\,\Psi\,\Lambda,
\\
Z_5
&
\,:=\,
\,\Psi\,\Lambda,
\\
Z_6
&
\,:=\,
\tfrac{3}{2}\,F_{0,3,1}\,
\big(3+2\,\Psi\big),
\\
Z_7
&
\,:=\,
-\,\tfrac{1}{2}\,\Phi
-
\Phi\,\Psi
-
2\,F_{0,3,1}\,\Psi,
\\
Z_8
&
\,:=\,
\tfrac{9}{2}\,F_{1,3,0}
-
\tfrac{3}{2}\,F_{3,1,0}
-
F_{1,2,1}\,\Phi
+
6\,\Psi\,F_{1,3,0}
-
2\,\Psi\,F_{3,1,0},
\\
Z_9
&
\,:=\,
-\,F_{0,3,1}\,\Lambda
+
\tfrac{1}{2}\,F_{1,3,0}
-
\tfrac{3}{2}\,F_{3,1,0}
+
2\,\Psi\,F_{1,3,0}
-
6\,\Psi\,F_{3,1,0}
-
F_{0,3,1}\,F_{1,2,1},
\
\\
Z_{10}
&
\,:=\,
-\,4\,F_{0,3,1}\,\Psi
-
2\,\Phi\,\Psi
-
\Phi,
\\
Z_{11}
&
\,:=\,
-\,F_{1,2,1}^2
-
F_{1,2,1}\,\Lambda
-
4\,F_{4,0,0},
\\
Z_{12}
&
\,:=\,
F_{1,2,1}^2
+
F_{1,2,1}\,\Lambda
+
4\,F_{4,0,0},
\\
Z_{13}
&
\,:=\,
-\,F_{1,2,1}^2
-
F_{1,2,1}\,\Lambda
-
8\,F_{4,0,0}\,\Psi
-
4\,F_{4,0,0},
\\
Z_{14}
&
\,:=\,
F_{1,2,1}\,\Phi
+
2\,F_{1,3,0}\,\Psi
+
2\,F_{3,1,0}\,\Psi
-
F_{3,1,0},
\\
Z_{15}
&
\,:=\,
2\,F_{1,2,1}\,\Phi
+
4\,F_{1,3,0}\,\Psi
+
4\,F_{3,1,0}\,\Psi
-
2\,F_{3,1,0},
\\
Z_{16}
&
\,:=\,
-\,F_{0,3,1}^2
-
F_{0,3,1}\,\Phi
+
8\,F_{0,4,0}\,\Psi
-
F_{1,2,1}^2
+
4\,F_{0,4,0},
\\
Z_{17}
&
\,:=\,
-\,3\,F_{0,3,1}\,F_{1,2,1}
-
3\,F_{0,3,1}\,\Lambda
+
F_{1,2,1}\,\Phi
-
3\,F_{1,3,0}
+
F_{3,1,0},
\\
Z_{18}
&
\,:=\,
-\,2\,F_{0,3,1}\,F_{1,2,1}
-
2\,F_{0,3,1}\,\Lambda
+
4\,F_{1,3,0}\,\Psi
+
4\,F_{3,1,0}\,\Psi
+
F_{1,3,0}
+
F_{3,1,0},
\\
Z_{19}
&
\,:=\,
-\,12\,F_{0,3,1}\,F_{4,0,0}
+
3\,F_{1,2,1}\,F_{1,3,0}
-
5\,F_{1,2,1}\,F_{3,1,0}
-
6\,F_{1,3,0}\,\Lambda
+
2\,F_{3,1,0}\,\Lambda
-
8\,F_{4,0,0}\,\Phi
-
6\,F_{2,3,0},
\\
Z_{20}
&
\,:=\,
F_{0,3,1}\,F_{1,3,0}
-
11\,F_{0,3,1}\,F_{3,1,0}
-
8\,F_{0,4,0}\,F_{1,2,1}
-
8\,F_{0,4,0}\,\Lambda
+
4\,F_{1,2,1}\,F_{4,0,0}
+
2\,F_{1,3,0}\,\Phi
-
6\,F_{3,1,0}\,\Phi
+
6\,F_{3,2,0}.
\endaligned
\]

We then apply the "{\footnotesize\sf EliminationIdeal}" procedure to
determine, in the ideal $\big\langle Z_i \big\rangle_{1 \leqslant
i \leqslant 20}$, the sub-ideal generated by
equations which involved {\em only} the
appearing $F$-variables:
\[
\big\{
F_{0,3,1},\,
F_{0,4,0},\,
F_{1,2,1},\,
F_{1,3,0},\,
F_{2,3,0},\,
F_{3,1,0},\,
F_{3,2,0},\,
F_{4,0,0}
\big\},
\]
hence involve {\em none} of the (unknown) variables (quantities):
\[
\big\{
\Lambda,\,
\Phi,\,
\Psi
\big\}.
\]

We obtain that this elimination ideal is generated by the 11
polynomials:
\[
\aligned
ZF_1
&
\,=\,
F_{1,3,0},
\\
ZF_2
&
\,=\,
F_{2,3,0},
\\
ZF_3
&
\,=\,
F_{3,1,0},
\\
ZF_4
&
\,=\,
F_{0,3,1}\,F_{1,2,1},
\\
ZF_5
&
\,=\,
F_{0,3,1}\,F_{4,0,0},
\\
ZF_6
&
\,=\,
F_{3,2,0}\,F_{0,3,1},
\\
ZF_7
&
\,=\,
F_{0,3,1}\,
\big(
F_{0,3,1}^2
-
4\,F_{0,4,0}
\big),
\\
ZF_8
&
\,=\,
F_{1,2,1}\,
\big(
F_{1,2,1}^2
-
4\,F_{0,4,0}
\big),
\\
ZF_9
&
\,=\,
F_{3,2,0}\,
\big(
F_{1,2,1}^2
-
4\,F_{0,4,0}
\big),
\\
ZF_{10}
&
\,=\,
2\,F_{1,2,1}\,F_{4,0,0}
+
F_{3,2,0},
\\
ZF_{11}
&
\,=\,
8\,F_{0,4,0}\,F_{4,0,0}
+
F_{1,2,1}\,F_{3,2,0}.
\endaligned
\]

Luckily, 2 of the 6 relative invariants must necessarily be
zero:
\[
\aligned
F_{1,3,0}
&
\,:=\,
0,
\\
F_{3,1,0}
&
\,:=\,
0,
\endaligned
\]
and, moreover, the following product of relative invariants:
\[
F_{0,3,1}\,
F_{1,2,1}
\,=\,
0,
\]
must always be zero!

This enables us to create only $4 < 10 = 11 - 1$ more subbranches
at order 4:
\[
\def\arraystretch{1.25}
\begin{array}{rcccccc}
\green{\bf 2b3e}\,\,\,\,
\green{\downarrow}\,\,
& F_{4,0,0} & F_{3,1,0} & F_{1,3,0} & F_{1,2,1} & F_{0,4,0} & F_{0,3,1}
\\
\green{\bf 4a} & 
0 & 0 & 0 & 0 & 0 & 0 
\\
\green{\bf 4b} & 
0 & 0 & 0 & 0 & 1 & 0
\\
\green{\bf 4c} & 
0 & 0 & 0 & 0 & -1 & 0
\\
\green{\bf 4d} & 
F_{4,0,0} & 0 & 0 & 1 & \tfrac{1}{4} & 0
\\
\green{\bf 4e} & 
0 & 0 & 0 & 0 & \tfrac{1}{4} & 1
\end{array}
\]

Once this is done, we launch again the standard method
in each of these 4 subbranches, and we find 4 (families of)
homogeneous models,
{\em see}~Section~{\ref{2b-models}}.

\SectionHead{Branch~{\green{\bf 2b3f}}}
{branch-2b-3f}

Unexpectedly, this branch will appear to be redundant with
Branch~{\green{\bf 2b3e}}.
At order 4, the linear representation is:
\[
\left(\!
\begin{array}{c}
G_{4,0,0}
\\
G_{3,1,0}
\\
G_{3,0,1}
\\
G_{2,1,1}
\\
G_{1,3,0}
\\
G_{0,4,0}
\end{array}
\!\right)
\,=\,
\left(\!
\def\arraystretch{1.25}
\begin{array}{cccccc}
\tfrac{1}{a_{1,1}^2} & 0 & 0 & 0 & 0 & 0
\\
0 & \tfrac{1}{a_{1,1}^2} & 0 & 0 & 0 & 0
\\
0 & 0 & \tfrac{1}{a_{1,1}} & 0 & 0 & 0
\\
0 & 0 & 0 & \tfrac{1}{a_{1,1}} & 0 & 0 
\\
0 & 0 & 0 & 0 & \tfrac{1}{a_{1,1}^2} & 0
\\
0 & 0 & 0 & 0 & 0 & \tfrac{1}{a_{1,1}^2} 
\end{array}
\!\right)\,
\left(\!
\begin{array}{c}
F_{4,0,0}
\\
F_{3,1,0}
\\
F_{3,0,1}
\\
F_{2,1,1}
\\
F_{1,3,0}
\\
F_{0,4,0}
\end{array}
\!\right).
\]
In Branch~{\green{\bf 2b3f}}, we set:
\[
F_{2,0,1}
\,:=\,
1, 
\ \ \ \ \ \ \ \
F_{1,1,1}
\,:=\,
0,
\ \ \ \ \ \ \ \
F_{3,0,0}
\,:=\,
0,
\ \ \ \ \ \ \ \
F_{0,3,0}
\,:=\,
0,
\]
while $F_{0,2,1}$ is an absolute invariant, {\em with the
hypothesis} $F_{0,2,1} \neq 1$.

At order 3, the equations of $\eqLF$ are:
\[
\aligned
0
&
{\overset{\green{\bf 3,0,0}}{\,\,\,=\,\,\,}}
-\,F_{0,2,1}\,A_{3,1}
+
A_{3,1}
+
\big(
F_{4,0,0}
-
2\,F_{2,2,0}
\big)\,
T_1
+
\big(
F_{3,1,0}
-
3\,F_{1,3,0}
\big)\,T_2
+
\big(
F_{3,0,1}
-
F_{1,2,1}
\big)\,T_3,
\\
0
&
{\overset{\green{\bf 0,3,0}}{\,\,\,=\,\,\,}}
F_{0,2,1}\,A_{3,2}
-
A_{3,2}
+
\big(
F_{1,3,0}
-
3\,F_{3,1,0}
\big)\,T_1
+
\big(
4\,F_{0,4,0}
-
2\,F_{2,2,0}
\big)\,T_2
+
\big(
F_{0,3,1}
-
F_{2,1,1}
\big)\,T_3,
\\
0
&
{\overset{\green{\bf 2,0,1}}{\,\,\,=\,\,\,}}
A_{3,3}
+
3\,F_{3,0,1}\,T_1
+
F_{2,1,1}\,T_2
+
\big(
2\,F_{2,0,2}
-
1
\big)\,T_3,
\\
0
&
{\overset{\green{\bf 1,1,1}}{\,\,\,=\,\,\,}}
2\,F_{0,2,1}\,A_{2,1}
-
2\,A_{2,1}
+
2\,F_{2,1,1}\,T_1
+
2\,F_{1,2,1}\,T_2
+
2\,F_{1,1,2}\,T_3,
\\
0
&
{\overset{\green{\bf 0,2,1}}{\,\,\,=\,\,\,}}
F_{0,2,1}\,A_{3,3}
+
F_{1,2,1}\,T_1
+
3\,F_{0,3,1}\,T_2
+
\big(
2\,F_{0,2,2}
-
F_{0,2,1}^2
\big)\,T_3,
\\
0
&
{\overset{\green{\bf 1,0,2}}{\,\,\,=\,\,\,}}
\big(
2\,F_{2,0,2}
-
2
\big)\,T_1
+
F_{1,1,2}\,T_2
+
3\,F_{1,0,3}\,T_3,
\\
0
&
{\overset{\green{\bf 0,1,2}}{\,\,\,=\,\,\,}}
F_{1,1,2}\,T_1
+
\big(
2\,F_{0,2,2}
-
2\,F_{0,2,1}^2
\big)\,T_2
+
3\,F_{0,1,3}\,T_3,
\\
0
&
{\overset{\green{\bf 0,0,3}}{\,\,\,=\,\,\,}}
F_{1,0,3}\,T_1
+
F_{0,1,3}\,T_2
+
4\,F_{0,0,4}\,T_3.
\endaligned
\]

We may solve 4 infinitesimal group parameters:
\[
\aligned
A_{3,3}
&
\,:=\,
-\,3\,F_{3,0,1}\,T_1
-
F_{2,1,1}\,T_2
+
\big(
1
-
2\,F_{2,0,2}
\big)\,T_3,
\\
A_{3,2}
&
\,:=\,
\frac{1}{F_{0,2,1}-1}\,
\Big\{
\big(
3\,F_{3,1,0}
-
F_{1,3,0}
\big)\,T_1
+
\big(
2\,F_{2,2,0}
-
4\,F_{0,4,0}
\big)\,T_2
+
\big(
F_{2,1,1}
-
F_{0,3,1}
\big)\,T_3
\Big\},
\\
A_{3,1}
&
\,:=\,
\frac{1}{F_{0,2,1}-1}\,
\Big\{
\big(
4\,F_{4,0,0}
-
2\,F_{2,2,0}
\big)\,T_1
+
\big(
F_{3,1,0}
-
3\,F_{1,3,0}
\big)\,T_2
+
\big(
F_{3,0,1}
-
F_{1,2,1}
\big)\,T_3
\Big\},
\\
A_{2,1}
&
\,:=\,
\frac{1}{F_{0,2,1}-1}\,
\Big\{
-\,F_{2,1,1}\,T_1
-
F_{1,2,1}\,T_2
-
F_{1,1,2}\,T_3
\Big\},
\endaligned
\]
and it remains the 4 equations:
\[
\aligned
0
&
{\overset{\green{\bf 0,2,1}}{\,\,\,=\,\,\,}}
\big(
F_{1,2,1}
-
3\,F_{0,2,1}\,F_{3,0,1}
\big)\,
T_1
+
\big(
3\,F_{0,3,1}
-
F_{0,2,1}\,F_{2,1,1}
\big)\,
T_2
\\
&
\ \ \ \ \
+
\big(
F_{0,2,1}
+
2\,F_{0,2,2}
-
F_{0,2,1}^2
-
2\,F_{0,2,1}\,F_{2,0,2}
\big)\,
T_3,
\\
0
&
{\overset{\green{\bf 1,0,2}}{\,\,\,=\,\,\,}}
\big(
2\,F_{2,0,2}
-
2
\big)\,T_1
+
F_{1,1,2}\,T_2
+
3\,F_{1,0,3}\,T_3,
\\
0
&
{\overset{\green{\bf 0,1,2}}{\,\,\,=\,\,\,}}
F_{1,1,2}\,T_1
+
\big(
2\,F_{0,2,2}
-
2\,F_{0,2,1}^2
\big)\,T_2
+
3\,F_{0,1,3}\,T_3,
\\
0
&
{\overset{\green{\bf 0,0,3}}{\,\,\,=\,\,\,}}
F_{1,0,3}\,T_1
+
F_{0,1,3}\,T_2
+
4\,F_{0,0,4}\,T_3,
\endaligned
\]
which imply, necessarily:
\[
\aligned
F_{1,0,3}
&
\,:=\,
0,
&
\ \ \ \ \ \ \ 
F_{1,1,2}
&
\,:=\,
0,
&
\ \ \ \ \ \ \ 
F_{1,2,1}
&
\,:=\,
3\,F_{0,2,1}\,F_{3,0,1},
&
\ \ \ \ \ \ \ 
F_{2,0,2}
&
\,:=\,
1,
\\
F_{0,1,3}
&
\,:=\,
0,
&
\ \ \ \ \ \ \ 
F_{0,2,2}
&
\,:=\,
F_{0,2,1}^2,
&
\ \ \ \ \ \ \
F_{0,3,1}
&
\,:=\,
\tfrac{1}{3}\,
F_{0,2,1}\,F_{2,1,1},
&
\ \ \ \ \ \ \
F_{0,0,4}
&
\,:=\,
0,
\endaligned
\]
and once these assignments are made, it remains:
\[
0
{\overset{\green{\bf 0,2,1}}{\,\,\,=\,\,\,}}
F_{0,2,1}\,
\big(
F_{0,2,1}
-
1
\big)\,
T_3,
\]
whence:
\[
F_{0,2,1}
\,:=\,
0.
\]

Consquently, Branch~{\green{\bf 2b3f}} reduces to:
\[
\def\arraystretch{1.25}
\begin{array}{rccccc}
\green{\bf 2b}\,\,\,\,
\green{\downarrow}\,\,
& 
F_{2,0,1} & F_{1,1,1} & F_{0,2,1} & F_{3,0,0} & F_{0,3,0}
\\
\green{\bf 3f} & 
1 & 0 & 0 & 0 & 0
\end{array}
\]
and since Branch~{\green{\bf 2b3e}} was:
\[
\def\arraystretch{1.25}
\begin{array}{rccccc}
\green{\bf 2b}\,\,\,\,
\green{\downarrow}\,\,
& 
F_{2,0,1} & F_{1,1,1} & F_{0,2,1} & F_{3,0,0} & F_{0,3,0}
\\
\green{\bf 3e} & 
0 & 0 & 1 & 0 & 0
\end{array}
\]
it is visible that, through the discrete symmetry:
\[
x
\,\,\,\longleftrightarrow\,\,\,
y,
\]
these two branches are the same. 
So Branch~{\green{\bf 2b3f}} disappears!

\SectionHead{Linear Representations Branch by Branch}
{linear-representations-branches}

Here, we present the occurring different branches 
and linear representations.

\Subsection{Branch~\green{\bf 2a}}
We recall that in the Branch \green{\bf 2a}, 
we found the linear representation:
\[
\left(\!
\begin{array}{c}
G_{2,0,1}
\\
G_{1,1,1}
\\
G_{0,2,1}
\\
G_{3,0,0}
\\
G_{0,3,0}
\end{array}
\!\right)
\,=\,
\left(\!
\def\arraystretch{1.25}
\begin{array}{ccccc}
\frac{a_{2,2}}{a_{1,1}a_{3,3}} & 0 & 0 & 0 & 0
\\
0 & \frac{1}{a_{3,3}} & 0 & 0 & 0
\\
0 & 0 & \frac{a_{1,1}}{a_{2,2}a_{3,3}} & 0 & 0
\\
-\frac{a_{2,2}a_{3,1}}{a_{1,1}^2a_{3,3}} & 0 & 0 & 
\frac{a_{2,2}}{a_{1,1}^2} & 0
\\
0 & 0 & -\frac{a_{1,1}a_{3,2}}{a_{2,2}^2a_{3,3}} & 0 &
\frac{a_{1,1}}{a_{2,2}^2}
\end{array}
\!\right)\,
\left(\!
\begin{array}{c}
F_{2,0,1}
\\
F_{1,1,1}
\\
F_{0,2,1}
\\
F_{3,0,0}
\\
F_{0,3,0}
\end{array}
\!\right),
\]
leading to the creation of $19$ branches: 
\[
\def\arraystretch{1.25}
\begin{array}{rccccc}
\green{\bf 2a}\,\,\,\,
\green{\downarrow}\,\,
& 
F_{2,0,1} & F_{1,1,1} & F_{0,2,1} & F_{3,0,0} & F_{0,3,0}
\\
\green{\bf 3a} & 
0 & 0 & 0 & 0 & 0
\\
\green{\bf 3b} & 
0 & 0 & 0 & 0 & 1
\\
\green{\bf 3c} & 
0 & 0 & 0 & 1 & 0
\\
\green{\bf 3d} &
0 & 0 & 0 & 1 & 1
\\
\green{\bf 3e} & 
0 & 0 & 1 & 0 & 0
\\
\green{\bf 3f} & 
0 & 0 & 1 & 1 & 0
\\
\green{\bf 3g} & 
0 & 1 & 0 & 0 & 0
\\
\green{\bf 3h} & 
0 & 1 & 0 & 0 & 1
\\
\green{\bf 3i} & 
0 & 1 & 0 & 1 & 0
\\
\green{\bf 3j} & 
0 & 1 & 0 & 1 & 1
\\
\green{\bf 3k} & 
0 & 1 & 1 & 0 & 0
\\
\green{\bf 3l} & 
0 & 1 & 1 & 1 & 0
\\
\green{\bf 3m} & 
1 & 0 & 0 & 0 & 0 
\\
\green{\bf 3n} &
1 & 0 & 0 & 0 & 1
\\
\green{\bf 3o} & 
1 & 0 & 1 & 0 & 0
\\
\green{\bf 3p} & 
1 & 0 & -1 & 0 & 0
\\
\green{\bf 3q} & 
1 & 1 & F_{0,2,1}^{\neq0} & 0 & 0
\\
\green{\bf 3r} & 
1 & 1 & 0 & 0 & 0
\\
\green{\bf 3s} & 
1 & 1 & 0 & 0 & 1
\end{array}
\]

Thanks to discrete symmetries like {\em e.g.}
the one $x \longleftrightarrow y$ used
to realize that Branch~{\green{\bf 2b3f}} is (unexpectedly) redundant
with Branch~{\green{\bf 2b3e}},
some of these 19 branches can be erased.
In this way, we can reduce the number of branches to $13$:
\[
\def\arraystretch{1.25}
\begin{array}{rccccc}
\green{\bf 2a}\,\,\,\,
\green{\downarrow}\,\,
& 
F_{2,0,1} & F_{1,1,1} & F_{0,2,1} & F_{3,0,0} & F_{0,3,0}
\\
\green{\bf 3a} & 
0 & 0 & 0 & 0 & 0
\\
\green{\bf 3b} & 
0 & 0 & 0 & 0 & 1
\\
\green{\bf 3c} & 
0 & 0 & 0 & 1 & 1
\\
\green{\bf 3d} &
0 & 0 & 1 & 0 & 0
\\
\green{\bf 3e} & 
0 & 0 & 1 & 1 & 0
\\
\green{\bf 3f} & 
0 & 1 & 0 & 0 & 0
\\
\green{\bf 3g} & 
0 & 1 & 0 & 0 & 1
\\
\green{\bf 3h} & 
0 & 1 & 0 & 1 & 1
\\
\green{\bf 3i} & 
0 & 1 & 1 & 0 & 0
\\
\green{\bf 3j} & 
0 & 1 & 1 & 1 & 0
\\
\green{\bf 3k} & 
1 & 0 & 1 & 0 & 0
\\
\green{\bf 3l} & 
1 & 0 & -1 & 0 & 0
\\
\green{\bf 3m} & 
1 & 1 & F_{0,2,1}^{\neq 0} & 0 & 0 
\end{array}
\]

Moreover, in the paragraphs below,
some branches like {\em e.g.} 
Branch~{\green{\bf 2a3c}} are {\em unmentioned},
and this means that no further sub-branches
exist at order $\geqslant 4$.

In Branch \green{\bf 2a3a}, 
at order 4, a single equation
remains, giving the linear representation:
\[
G_{2,2,0} := \frac{1}{a_{1,1}a_{2,2}}F_{2,2,0}.
\]
This leads us to the creation of $2$ branches, 
mutually inequivalent and having empty intersection:
\[
\def\arraystretch{1.25}
\begin{array}{rccc}
\green{\bf 2a3a}\,\,\,\,
\green{\downarrow}\,\,
& F_{2,2,0}
\\
\green{\bf 4a} & 
1
\\
\green{\bf 4b} & 
0
\end{array}
\]

In Branch \green{\bf 2a3b}, we find:
\[
\left(\!
\begin{array}{c}
G_{0,4,0}
\\
G_{1,3,0}
\\
G_{2,2,0}
\end{array}
\!\right)
\,=\,
\left(\!
\def\arraystretch{1.25}
\begin{array}{ccc}
\tfrac{1}{a_{2,2}} & 0 & 0
\\
0 & \tfrac{1}{a_{2,2}^2} & 0
\\
0 & 0 & \tfrac{1}{a_{2,2}^3}
\end{array}
\!\right)\,
\left(\!
\begin{array}{c}
F_{0,4,0}
\\
F_{1,3,0}
\\
F_{2,2,0}
\end{array}
\!\right),
\]
leading to:
\[
\def\arraystretch{1.25}
\begin{array}{rccc}
\green{\bf 2a3b}\,\,\,\,
\green{\downarrow}\,\,
& F_{0,4,0} & F_{1,3,0} & F_{2,2,0}
\\
\green{\bf 4a} & 
1 & F_{1,3,0} & F_{2,2,0}
\\
\green{\bf 4b} & 
0 & 1 & F_{2,2,0}
\\
\green{\bf 4c} & 
0 & -1 & F_{2,2,0}
\\
\green{\bf 4d} & 
0 & 0 & 1
\\
\green{\bf 4e} & 
0 & 0 & 0
\end{array}
\]

In Branch \green{\bf 2a3d}, we find:
\[
\left(\!
\begin{array}{c}
G_{2,2,0}
\\
G_{0,4,0}
\end{array}
\!\right)
\,=\,
\left(\!
\def\arraystretch{1.25}
\begin{array}{cc}
\tfrac{1}{a_{2,2}^2a_{3,3}} & 0
\\
0 & \tfrac{a_{3,3}}{a_{2,2}^2}
\end{array}
\!\right)\,
\left(\!
\begin{array}{c}
F_{2,2,0}
\\
F_{0,4,0}
\end{array}
\!\right),
\]
leading to:
\[
\def\arraystretch{1.25}
\begin{array}{rccc}
\green{\bf 2a3d}\,\,\,\,
\green{\downarrow}\,\,
& F_{2,2,0} & F_{0,4,0}
\\
\green{\bf 4a} & 
1 & 1
\\
\green{\bf 4b} & 
1 & -1
\\
\green{\bf 4c} & 
1 & 0
\\
\green{\bf 4d} & 
0 & 1
\\
\green{\bf 4e} & 
0 & 0
\end{array}
\]

\medskip

In Branch \green{\bf 2a3d4c}, we find:
\[
G_{2,3,0} := \frac{1}{a_{2,2}}F_{2,3,0},
\]
leading to:
\[
\def\arraystretch{1.25}
\begin{array}{rccc}
\green{\bf 2a3d4c}\,\,\,\,
\green{\downarrow}\,\,
& F_{2,3,0}
\\
\green{\bf 5a} & 
1
\\
\green{\bf 5b} & 
0
\end{array}
\]

In Branch \green{\bf 2a3d4d}, we find:
\[
\left(\!
\begin{array}{c}
G_{0,5,0}
\\
G_{0,4,1}
\\
G_{1,4,0}
\end{array}
\!\right)
\,=\,
\left(\!
\def\arraystretch{1.25}
\begin{array}{ccc}
\tfrac{1}{a_{2,2}} & 0 & 0
\\
0 & \tfrac{1}{a_{2,2}^2} & 0
\\
0 & 0 & \tfrac{1}{a_{2,2}^3}
\end{array}
\!\right)\,
\left(\!
\begin{array}{c}
F_{0,5,0}
\\
F_{0,4,1}
\\
F_{1,4,0}
\end{array}
\!\right),
\]
leading to:
\[
\def\arraystretch{1.25}
\begin{array}{rccc}
\green{\bf 2a3d4d}\,\,\,\,
\green{\downarrow}\,\,
& F_{0,5,0} & F_{0,4,1} & F_{1,4,0}
\\
\green{\bf 5a} & 
1 & F_{0,4,1} & F_{1,4,0}
\\
\green{\bf 5b} & 
0 & 1 & F_{1,4,0}
\\
\green{\bf 5c} & 
0 & -1 & F_{1,4,0}
\\
\green{\bf 5d} & 
0 & 0 & 1
\\
\green{\bf 5e} & 
0 & 0 & 0
\end{array}
\]

In Branch \green{\bf 2a3g}, we find:
\[
\left(\!
\begin{array}{c}
G_{1,4,0}
\\
G_{0,5,0}
\\
G_{2,3,0}
\end{array}
\!\right)
\,=\,
\left(\!
\def\arraystretch{1.25}
\begin{array}{ccc}
\tfrac{1}{a_{2,2}^3} & 0 & 0
\\
0 & \tfrac{1}{a_{2,2}^2} & 0
\\
0 & 0 & \tfrac{1}{a_{2,2}^4}
\end{array}
\!\right)\,
\left(\!
\begin{array}{c}
F_{1,4,0}
\\
F_{0,5,0}
\\
F_{2,3,0}
\end{array}
\!\right),
\]
leading to:
\[
\def\arraystretch{1.25}
\begin{array}{rccc}
\green{\bf 2a3g}\,\,\,\,
\green{\downarrow}\,\,
& F_{1,4,0} & F_{0,5,0} & F_{2,3,0}
\\
\green{\bf 5a} & 
1 & F_{0,5,0} & F_{2,3,0}
\\
\green{\bf 5b} & 
0 & 1 & F_{2,3,0}
\\
\green{\bf 5c} & 
0 & -1 & F_{2,3,0}
\\
\green{\bf 5d} & 
0 & 0 & 1
\\
\green{\bf 5e} & 
0 & 0 & 0
\end{array}
\]

In Branch \green{\bf 2a3i}, we find:
\[
\left(\!
\begin{array}{c}
G_{1,2,1}
\\
G_{0,4,0}
\\
G_{2,2,0}
\end{array}
\!\right)
\,=\,
\left(\!
\def\arraystretch{1.25}
\begin{array}{ccc}
\tfrac{1}{a_{1,1}} & 0 & 0
\\
0 & \tfrac{1}{a_{1,1}^2} & 0
\\
0 & 0 & \tfrac{1}{a_{1,1}^2}
\end{array}
\!\right)\,
\left(\!
\begin{array}{c}
F_{1,2,1}
\\
F_{0,4,0}
\\
F_{2,2,0}
\end{array}
\!\right),
\]
leading to:
\[
\def\arraystretch{1.25}
\begin{array}{rccc}
\green{\bf 2a3i}\,\,\,\,
\green{\downarrow}\,\,
& F_{1,2,1} & F_{0,4,0} & F_{2,2,0}
\\
\green{\bf 4a} & 
1 & F_{0,4,0} & F_{2,2,0}
\\
\green{\bf 4b} & 
0 & 1 & F_{2,2,0}
\\
\green{\bf 4c} & 
0 & -1 & F_{2,2,0}
\\
\green{\bf 4d} & 
0 & 0 & 1
\\
\green{\bf 4e} & 
0 & 0 & -1
\\
\green{\bf 4f} & 
0 & 0 & 0
\end{array}
\]

In Branch \green{\bf 2a3m}, we find:
\[
\left(\!
\begin{array}{c}
G_{3,0,1}
\\
G_{2,1,1}
\\
G_{0,4,0}
\\
G_{1,3,0}
\\
G_{2,2,0}
\\
G_{4,0,0}
\end{array}
\!\right)
\,=\,
\left(\!
\def\arraystretch{1.25}
\begin{array}{cccccc}
\tfrac{1}{a_{1,1}} & 0 & 0 & 0 & 0 & 0
\\
0 & \tfrac{1}{a_{1,1}} & 0 & 0 & 0 & 0
\\
0 & 0 & \tfrac{1}{a_{1,1}^2} & 0 & 0 & 0
\\
0 & 0 & 0 & \tfrac{1}{a_{1,1}^2} & 0 & 0
\\
0 & 0 & 0 & 0 & \tfrac{1}{a_{1,1}^2} & 0
\\
0 & 0 & 0 & 0 & 0 & \tfrac{1}{a_{1,1}^2}
\end{array}
\!\right)\,
\left(\!
\begin{array}{c}
F_{3,0,1}
\\
F_{2,1,1}
\\
F_{0,4,0}
\\
F_{1,3,0}
\\
F_{2,2,0}
\\
F_{4,0,0}
\end{array}
\!\right),
\]
leading to:
\[
\def\arraystretch{1.25}
\begin{array}{rcccccc}
\green{\bf 2a3m}\,\,\,\,
\green{\downarrow}\,\,
& F_{3,0,1} & F_{2,1,1} & F_{0,4,0} & F_{1,3,0} & F_{2,2,0} & F_{4,0,0}
\\
\green{\bf 4a} & 
1 & F_{2,1,1} & F_{0,4,0} & F_{1,3,0} & F_{2,2,0} & F_{4,0,0}
\\
\green{\bf 4b} & 
0 & 1 & F_{0,4,0} & F_{1,3,0} & F_{2,2,0} & F_{4,0,0}
\\
\green{\bf 4c} & 
0 & 0 & 1 & F_{1,3,0} & F_{2,2,0} & F_{4,0,0}
\\
\green{\bf 4d} & 
0 & 0 & -1 & F_{1,3,0} & F_{2,2,0} & F_{4,0,0}
\\
\green{\bf 4e} & 
0 & 0 & 0 & 1 & F_{2,2,0} & F_{4,0,0}
\\
\green{\bf 4f} & 
0 & 0 & 0 & -1 & F_{2,2,0} & F_{4,0,0}
\\
\green{\bf 4g} & 
0 & 0 & 0 & 0 & 1 & F_{4,0,0}
\\
\green{\bf 4h} & 
0 & 0 & 0 & 0 & -1 & F_{4,0,0}
\\
\green{\bf 4i} & 
0 & 0 & 0 & 0 & 0 & 1
\\
\green{\bf 4j} & 
0 & 0 & 0 & 0 & 0 & -1
\\
\green{\bf 4k} & 
0 & 0 & 0 & 0 & 0 & 0
\end{array}
\]

\Subsection{Branch~\green{\bf 2b}}
We recall from Section~{\ref{order-3-linear-representations}}
that there are 6 branches:
\[
\def\arraystretch{1.25}
\begin{array}{rccccc}
\green{\bf 2b}\,\,\,\,
\green{\downarrow}\,\,
& 
F_{2,0,1} & F_{1,1,1} & F_{0,2,1} & F_{3,0,0} & F_{0,3,0}
\\
\green{\bf 3a} & 
0 & 0 & 0 & 0 & 0
\\
\green{\bf 3b} & 
0 & 0 & 0 & 1 & 0
\\
\green{\bf 3c} & 
1 & 0 & 1 & 0 & 0
\\
\green{\bf 3d} & 
1 & 0 & 1 & 1 & 0
\\
\green{\bf 3e} & 
0 & 0 & 1 & 0 & 0
\\
\green{\bf 3f} & 
1 & 0 & F_{0,2,1}^{\neq0} & 0 & 0
\end{array}
\]

In Branch~\green{\bf 2b3a}, we find:
\[
G_{2,2,0} :=\frac{1}{a_{1,1}^2+a_{2,1}^2}F_{2,2,0},
\]
leading to:
\[
\def\arraystretch{1.25}
\begin{array}{rccc}
\green{\bf 2b3a}\,\,\,\,
\green{\downarrow}\,\,
& F_{2,2,0}
\\
\green{\bf 4a} & 
0
\\
\green{\bf 4b} & 
1
\\
\green{\bf 4c} & 
-1
\end{array}
\]

Branches~\green{\bf 2b3b}, \green{\bf 2b3c}, 
\green{\bf 2b3d} do not lead to further linear
representations at orders $\geqslant 4$.

In Section~{\ref{branch-2b-3f}}, we have already explained
how Branch~\green{\bf 2b3f} disappears.

Finally, in Branch~\green{\bf 2b3e}, 
we already explained that the linear representation:
\[
\left(\!
\begin{array}{c}
G_{4,0,0}
\\
G_{3,1,0}
\\
G_{1,3,0}
\\
G_{1,2,1}
\\
G_{0,4,0}
\\
G_{0,3,1}
\end{array}
\!\right)
\,=\,
\left(\!
\def\arraystretch{1.25}
\begin{array}{cccccc}
\tfrac{1}{a_{1,1}^2} & 0 & 0 & 0 & 0 & 0
\\
0 & \tfrac{1}{a_{1,1}^2} & 0 & 0 & 0 & 0
\\
0 & 0 & \tfrac{1}{a_{1,1}^2} & 0 & 0 & 0
\\
0 & 0 & 0 & \tfrac{1}{a_{1,1}} & 0 & 0 
\\
0 & 0 & 0 & 0 & \tfrac{1}{a_{1,1}^2} & 0
\\
0 & 0 & 0 & 0 & 0 & \tfrac{1}{a_{1,1}} 
\end{array}
\!\right)\,
\left(\!
\begin{array}{c}
F_{4,0,0}
\\
F_{3,1,0}
\\
F_{1,3,0}
\\
F_{1,2,1}
\\
F_{0,4,0}
\\
F_{0,3,1}
\end{array}
\!\right),
\]
joint with an "{\footnotesize\sf EliminationIdeal}" trick, 
leads to:
\[
\def\arraystretch{1.25}
\begin{array}{rcccccc}
\green{\bf 2b3e}\,\,\,\,
\green{\downarrow}\,\,
& F_{4,0,0} & F_{3,1,0} & F_{1,3,0} & F_{1,2,1} & F_{0,4,0} & F_{0,3,1}
\\
\green{\bf 4a} & 
0 & 0 & 0 & 0 & 0 & 0 
\\
\green{\bf 4b} & 
0 & 0 & 0 & 0 & 1 & 0
\\
\green{\bf 4c} & 
0 & 0 & 0 & 0 & -1 & 0
\\
\green{\bf 4d} & 
F_{4,0,0} & 0 & 0 & 1 & \tfrac{1}{4} & 0
\\
\green{\bf 4e} & 
0 & 0 & 0 & 0 & \tfrac{1}{4} & 1
\end{array}
\]

\SectionHead{\bf 2a Models}
{2a-models}

Here, we present the Branch~{\green{\bf 2a}} models.

\[
\text{\bf Model 2b3a4a}\ \ \ \ \
\aligned
u = 2x^3y^3+x^2y^2+xy+\cdots.
\endaligned
\]
\[
\def\arraystretch{1.25}
\begin{array}{llll}
e_1 := -(2u-1)\partial_x+y\partial_u, &
e_2 :=-(2u-1)\partial_y+x\partial_u, &
e_3 := \partial_z, &
e_4 := x\partial_x-y\partial_y,
\\
e_5 := x\partial_z, &
e_6 := y\partial_z, &
e_7 := z\partial_z, &
e_8 := u\partial_z,
\end{array}
\]
\[
\footnotesize
\def\arraystretch{1.25}
\begin{array}{c|cccccccc}
{} & e_1 & e_2 & e_3 & e_4 & e_5 & e_6 & e_7 & e_8 
\\
\hline
e_1 & 
0 & 2e_4 & 0 & e_1 & e_3-2e_8 & 0 & 0 & e_6
\\
e_2 &
-2e_4 & 0 & 0 & -e_2 & 0 & e_3-2e_8 & 0 & e_5
\\
e_3 &
0 & 0 & 0 & 0 & 0 & 0 & 0e_3 & 0
\\
e_4 &
-e_1 & e_2 & 0 & 0 & e_5 & -e_6 & 0 & 0
\\
e_5 &
-e_3+2e_8 & 0 & 0 & -e_5 & 0 & 0 & e_5 & 0
\\
e_6 &
0 & -e_3+2e_8 & 0 & e_6 & 0 & 0 & e_6 & 0
\\
e_7 &
0 & 0 & -e_3 & 0 & -e_5 & -e_6 & 0 & -e_8
\\
e_8 &
-e_6 & -e_5 & 0 & 0 & 0 & 0 & e_8 & 0
\end{array}
\]



\[
\text{\bf Model 2a3a4b}
\ \ \ \ \
\aligned
u
=
xy.
\endaligned
\]
\[
\def\arraystretch{1.25}
\begin{array}{llll}
e_1
\,:=\,
\partial_x+y\partial_u
, &
e_2
\,:=\,
\partial_y+x\partial_u
, &
& 
\\
e_3
\,:=\,
\partial_z
, &
e_4
\,:=\,
x\partial_x+u\partial_u
, &
e_5
\,:=\,
y\partial_y+u\partial_u
, &
e_6
\,:=\,
x\partial_z,
\\
e_7
\,:=\,
y\partial_z
, &
e_8
\,:=\,
z\partial_z
, &
e_9
\,:=\,
u\partial_z.
\end{array}
\]

\[
\footnotesize
\def\arraystretch{1.25}
\begin{array}{c|ccccccccc}
{} & e_1 & e_2 & e_3 & e_4 & e_5 & e_6 & e_7 & e_8 & e_9 
\\
\hline
e_1 & 
0 & 0 & 0 & e_1 & 0 & e_3 & 0 & 0 & e_7
\\
e_2 &
0 & 0 & 0 & 0 & e_2 & 0 & e_3 & 0 & e_6
\\
e_3 &
0 & 0 & 0 & 0 & 0 & 0 & 0 & e_3 & 0
\\
e_4 &
-e_1 & 0 & 0 & 0 & 0 & e_6 & 0 & 0 & e_9
\\
e_5 &
0 & -e_2 & 0 & 0 & 0 & 0 & e_7 & 0 & e_9
\\
e_6 &
-e_3 & 0 & 0 & -e_6 & 0 & 0 & 0 & e_6 & 0
\\
e_7 &
0 & -e_3 & 0 & 0 & -e_7 & 0 & 0 & e_7 & 0
\\
e_8 &
0 & 0 & -e_3 & 0 & 0 & -e_6 & -e_7 & 0 & -e_9
\\
e_9 &
-e_7 & -e_6 & 0 & -e_9 & -e_9 & 0 & 0 & e_9 & 0
\end{array}
\]


\[
\text{\bf Model 2a3b4a}
\ \ \ \ \
\left\{
\aligned
u
&
=
xy+y^3+F_{1,3,0}xy^3+y^4
+
\big(-F_{1,3,0}^2+F_{1,3,0}(\tfrac{5}{2}F_{0,5,0}F_{1,3,0}-3F_{1,3,0}^2
+
\\
&
\ \ \ \ \
-2F_{1,3,0})\big)x^2y^3
+
(\tfrac{5}{2}F_{0,5,0}F_{1,3,0}-3F_{1,3,0}^2-2F_{1,3,0})xy^4
+
F_{0,5,0}y^5
+
\\
&
\ \ \ \ \
+	(13F_{1,3,0}^3+25F_{1,3,0}^4+12F_{1,3,0}^5-\tfrac{125}{6}F_{0,5,0}F_{1,3,0}^3+\tfrac{25}{3}F_{0,5,0}^2F_{1,3,0}^3
+
\\
&
\ \ \ \ \
-20F_{0,5,0}F_{1,3,0}^4)x^3y^3
+
(9F_{1,3,0}^2+\tfrac{45}{2}F_{1,3,0}^3+\tfrac{27}{2}F_{1,3,0}^4-\tfrac{75}{4}F_{0,5,0}F_{1,3,0}^2
+
\\
&
\ \ \ \ \
-\tfrac{45}{2}F_{0,5,0}F_{1,3,0}^3+\tfrac{75}{8}F_{0,5,0}^2F_{1,3,0}^2)x^2y^4
+
(5F_{0,5,0}^2F_{1,3,0}-6F_{0,5,0}F_{1,3,0}^2
+
\\
&
\ \ \ \ \
-5F_{0,5,0}F_{1,3,0})xy^5
+
(-\tfrac{1}{2}F_{1,3,0}-\tfrac{3}{4}F_{0,5,0}F_{1,3,0}-\tfrac{3}{2}F_{1,3,0}^2+\tfrac{5}{3}F_{0,5,0}^2
+
\\
&
\ \ \ \ \
-\tfrac{2}{3}F_{0,5,0})y^6
+
\cdots.
\endaligned\right.
\]
\[
\def\arraystretch{1.25}
\begin{array}{ll}
e_1
&
\,:=\,
-(5F_{0,5,0}F_{1,3,0}x-6F_{1,3,0}^2x-7F_{1,3,0}x-1)\partial_x-(-3yF_{1,3,0}+\tfrac{5}{2}yF_{0,5,0}F_{1,3,0}-3yF_{1,3,0}^2)\partial_y
+
\\
&
\ \ \ \ \
-(-y-10uF_{1,3,0}+\tfrac{15}{2}uF_{0,5,0}F_{1,3,0}-9uF_{1,3,0}^2)\partial_u,
\\
e_2
&
\,:=\,
-(10F_{0,5,0}x+3F_{1,3,0}u-12F_{1,3,0}x-12x+3y)\partial_x-(5F_{0,5,0}y-6F_{1,3,0}y-4y-1)\partial_y
+
\\
&
\ \ \ \ \
-(15F_{0,5,0}u-18F_{1,3,0}u-16u-x)\partial_u,
\\
e_3
&
\,:=\,
\partial_z
,
\,\,\,\,\,\,\,\,\,\,\,\,\,\,\,\,\,\,\,\,\,\,\,\,
e_5
\,:=\,
y\partial_z
,
\,\,\,\,\,\,\,\,\,\,\,\,\,\,\,\,\,\,\,\,\,\,\,\,
e_7
\,:=\,
u\partial_z
,
\\
e_4
&
\,:=\,
x\partial_z
,
\,\,\,\,\,\,\,\,\,\,\,\,\,\,\,\,\,\,\,\,\,
e_6
\,:=\,
z\partial_z.
\end{array}
\]
Gröbner basis generators of 
moduli space core algebraic variety in 
$\R^2 \ni (F_{0,5,0},F_{1,3,0})$:
\[
\aligned
B_1 & := F_{1,3,0}^2(5F_{0,5,0}-6F_{1,3,0}-6),
\\
B_2 & := F_{1,3,0}(5F_{0,5,0}-8+6F_{1,3,0})(5F_{0,5,0}-6F_{1,3,0}-6).
\endaligned
\]
\[
\footnotesize
\def\arraystretch{1.25}
\begin{array}{c|ccccccc}
{} & e_1 & e_2 & e_3 & e_4 & e_5 & e_6 & e_7
\\
\hline
e_1 & 
0 &
\rotatebox[origin=c]{90}{
\begin{tabular}{p{4cm}}
$(-10F_{0,5,0}+12F_{1,3,0}+12)e_1+(-3F_{1,3,0}+\tfrac{5}{2}F_{0,5,0}F_{1,3,0}-3F_{1,3,0}^2)e_2$
\end{tabular}} &
0 &
\rotatebox[origin=c]{90}{
\begin{tabular}{p{4cm}}
$e_3+(-5F_{0,5,0}F_{1,3,0}+6F_{1,3,0}^2+7F_{1,3,0})e_4$
\end{tabular}} &
\rotatebox[origin=c]{90}{
\begin{tabular}{p{4cm}}
$(3F_{1,3,0}-(5/2)F_{0,5,0}F_{1,3,0}+3F_{1,3,0}^2)e_5$
\end{tabular}} &
0 &
\rotatebox[origin=c]{90}{
\begin{tabular}{p{4cm}}
$e_5+(10F_{1,3,0}-(15/2)F_{0,5,0}F_{1,3,0}+9F_{1,3,0}^2)e_7$
\end{tabular}}
\\
e_2 &
\rotatebox[origin=c]{90}{
\begin{tabular}{p{4cm}}
$-(-10F_{0,5,0}+12F_{1,3,0}+12)e_1-(-3F_{1,3,0}+\tfrac{5}{2}F_{0,5,0}F_{1,3,0}-3F_{1,3,0}^2)e_2$
\end{tabular}} &
0 & 0 &
\rotatebox[origin=c]{90}{
\begin{tabular}{p{4cm}}
$(-10F_{0,5,0}+12F_{1,3,0}+12)e_4-3e_5-3F_{1,3,0}e_7$
\end{tabular}} &
\rotatebox[origin=c]{90}{
\begin{tabular}{p{4cm}}
$e_3+(-5F_{0,5,0}+6F_{1,3,0}+4)e_5$
\end{tabular}} &
0 &
\rotatebox[origin=c]{90}{
\begin{tabular}{p{4cm}}
$e_4+(-15F_{0,5,0}+18F_{1,3,0}+16)e_7$
\end{tabular}}
\\
e_3 &
0 & 0 & 0 & 0 & 0 & e_3 & 0
\\
e_4 &
\rotatebox[origin=c]{90}{
\begin{tabular}{p{4cm}}
$-e_3-(-5F_{0,5,0}F_{1,3,0}+6F_{1,3,0}^2+7F_{1,3,0})e_4$
\end{tabular}} &
\rotatebox[origin=c]{90}{
\begin{tabular}{p{4cm}}
$-(-10F_{0,5,0}+12F_{1,3,0}+12)e_4+3e_5+3F_{1,3,0}e_7$
\end{tabular}}
 &
0 & 0 & 0 &
e_4 & 0
\\
e_5 &
\rotatebox[origin=c]{90}{
\begin{tabular}{p{4cm}}
$-(3F_{1,3,0}-\tfrac{5}{2}F_{0,5,0}F_{1,3,0}+3F_{1,3,0}^2)e_5$
\end{tabular}} &
\rotatebox[origin=c]{90}{
\begin{tabular}{p{4cm}}
$-e_3-(-5F_{0,5,0}+6F_{1,3,0}+4)e_5$
\end{tabular}} &
0 & 0 & 0 & e_5 & 0
\\
e_6 &
0 & 0 & -e_3 & -e_4 & -e_5 & 0 & -e_7
\\
e_7 &
\rotatebox[origin=c]{90}{
\begin{tabular}{p{4cm}}
$-e_5-(10F_{1,3,0}-\tfrac{15}{2}F_{0,5,0}F_{1,3,0}+9F_{1,3,0}^2)e_7$
\end{tabular}} &
\rotatebox[origin=c]{90}{
\begin{tabular}{p{4cm}}
$-e_4-(-15F_{0,5,0}+18F_{1,3,0}+16)e_7$
\end{tabular}} &
0 & 0 & 0 & e_7 & 0
\end{array}
\]


\[
\text{\bf Model 2a3b4b}
\ \ \ \ \
\left\{
\aligned
u
=
xy+y^3+xy^3+\tfrac{6}{5}y^5+\tfrac{6}{5}xy^5+\tfrac{51}{35}y^7
+
\cdots
.
\endaligned\right.
\]
\[
\def\arraystretch{1.25}
\begin{array}{llll}
e_1
\,:=\,
(x+1)\partial_x+(u+y)\partial_u
, &
e_2
\,:=\,
-(3u+3y)\partial_x+\partial_y+x\partial_u
, &
& 
\\
e_3
\,:=\,
\partial_z
, &
e_4
\,:=\,
x\partial_z
, &
e_5
\,:=\,
y\partial_z
,
\\
e_6
\,:=\,
z\partial_z,
&
e_7
\,:=\,
u\partial_z.
\end{array}
\]

\[
\footnotesize
\def\arraystretch{1.25}
\begin{array}{c|ccccccc}
{} & e_1 & e_2 & e_3 & e_4 & e_5 & e_6 & e_7
\\
\hline
e_1 & 
0 & 0 & 0 & e_3+e_4 & 0 & 0 & e_5+e_7
\\
e_2 &
0 & 0 & 0 & -3e_5-3e_7 & e_3 & 0 & e_4
\\
e_3 &
0 & 0 & 0 & 0 & 0 & e_3 & 0
\\
e_4 &
-e_3-e_4 & 3e_5+3e_7 & 0 & 0 & 0 & e_4 & 0
\\
e_5 &
0 & -e_3 & 0 & 0 & 0 & e_5 & 0
\\
e_6 &
0 & 0 & -e_3 & -e_4 & -e_5 & 0 & -e_7
\\
e_7 &
-e_5-e_7 & -e_4 & 0 & 0 & 0 & e_7 & 0
\end{array}
\]


\[
\text{\bf Model 2a3b4c}
\ \ \ \ \
\left\{
\aligned
u
=
xy+y^3-xy^3-\tfrac{6}{5}y^5+\tfrac{6}{5}xy^5
+
\cdots
.
\endaligned\right.
\]
\[
\def\arraystretch{1.25}
\begin{array}{llll}
e_1
\,:=\,
-(x-1)\partial_x-(u-y)\partial_u, &
e_2
\,:=\,
(3u-3y)\partial_x+\partial_y+x\partial_u
, &
& 
\\
e_3
\,:=\,
\partial_z
, &
e_4
\,:=\,
x\partial_z
, &
e_5
\,:=\,
y\partial_z
,
\\
e_6
\,:=\,
z\partial_z
, &
e_7
\,:=\,
u\partial_z.
\end{array}
\]

\[
\footnotesize
\def\arraystretch{1.25}
\begin{array}{c|ccccccc}
{} & e_1 & e_2 & e_3 & e_4 & e_5 & e_6 & e_7
\\
\hline
e_1 & 
0 & 0 & 0 & e_3-e_4 & 0 & 0 & e_5-e_7
\\
e_2 &
0 & 0 & 0 & -3e_5+3e_7 & e_3 & 0 & e_4
\\
e_3 &
0 & 0 & 0 & 0 & 0 & e_3 & 0
\\
e_4 &
-e_3+e_4 & 3e_5-3e_7 & 0 & 0 & 0 & e_4 & 0
\\
e_5 &
0 & -e_3 & 0 & 0 & 0 & e_5 & 0
\\
e_6 &
0 & 0 & -e_3 & -e_4 & -e_5 & 0 & -e_7
\\
e_7 &
-e_5+e_7 & -e_4 & 0 & 0 & 0 & e_7 & 0
\end{array}
\]


\[
\text{\bf Model 2a3b4e}
\ \ \ \ \
\left\{
u
\,=\,
y^3+xy
.
\right.
\]
\[
\def\arraystretch{1.25}
\begin{array}{llll}
e_1
\,:=\,
\partial_x+y\partial_u
, &
e_2
\,:=\,
-3y\partial_x+\partial_y+x\partial_u
, &
e_3
\,:=\,
\partial_z
,
\\
e_4
\,:=\,
2x\partial_x+y\partial_y+3u\partial_u
, &
e_5
\,:=\,
x\partial_z
, &
e_6
\,:=\,
y\partial_z
,
\\
e_7
\,:=\,
z\partial_z
, &
e_8
\,:=\,
u\partial_z
.
\end{array}
\]

\[
\footnotesize
\def\arraystretch{1.25}
\begin{array}{c|cccccccc}
{} & e_1 & e_2 & e_3 & e_4 & e_5 & e_6 & e_7 & e_8 
\\
\hline
e_1 & 
0 & 0 & 0 & 2e_1 & e_3 & 0 & 0 & e_6
\\
e_2 &
0 & 0 & 0 & e_2 & -3e_6 & e_3 & 0 & e_5
\\
e_3 &
0 & 0 & 0 & 0 & 0 & 0 & e_3 & 0
\\
e_4 &
-2e_1 & -e_2 & 0 & 0 & 2e_5 & e_6 & 0 & 3e_8
\\
e_5 &
-e_3 & 3e_6 & 0 & -2e_5 & 0 & 0 & e_5 & 0
\\
e_6 &
0 & -e_3 & 0 & -e_6 & 0 & 0 & e_6 & 0
\\
e_7 &
0 & 0 & -e_3 & 0 & -e_5 & -e_6 & 0 & -e_8
\\
e_8 &
-e_6 & -e_5 & 0 & -3e_8 & 0 & 0 & e_8 & 0
\end{array}
\]


\[
\!\!\!\!\!\!\!
\text{\bf Model 2a3c}
\ \ \ \ \
\left\{
\aligned
u
&
=
xy+x^3+y^3+F_{4,0,0}x^4+F_{3,1,0}x^3y
+
(\tfrac{2}{9}F_{0,4,0}F_{1,3,0}-\tfrac{1}{9}F_{3,1,0}F_{1,3,0}
+
\\
&
\ \ \ \ \
+\tfrac{8}{9}F_{0,4,0}F_{4,0,0})x^2y^2
+
F_{1,3,0}xy^3+F_{0,4,0}y^4
+
(\tfrac{1}{15}F_{4,0,0}F_{1,3,0}+\tfrac{4}{3}F_{4,0,0}^2
+
\\
&
\ \ \ \ \
+\tfrac{6}{5}F_{3,1,0})x^5
+
(\tfrac{7}{9}F_{0,4,0}F_{1,3,0}-\tfrac{2}{9}F_{3,1,0}F_{1,3,0}+\tfrac{28}{9}F_{0,4,0}F_{4,0,0}-9
+
\\
&
\ \ \ \ \
+\tfrac{4}{3}F_{3,1,0}F_{4,0,0})x^4y
+
\big(\tfrac{1}{3}F_{1,3,0}(\tfrac{2}{9}F_{0,4,0}F_{1,3,0}-\tfrac{1}{9}F_{3,1,0}F_{1,3,0}+\tfrac{8}{9}F_{0,4,0}F_{4,0,0})
+
\\
&
\ \ \ \ \
+\tfrac{4}{3}(\tfrac{2}{9}F_{0,4,0}F_{1,3,0}-\tfrac{1}{9}F_{3,1,0}F_{1,3,0}+\tfrac{8}{9}F_{0,4,0}F_{4,0,0})F_{4,0,0}+3F_{1,3,0}\big)x^3y^2
+
\\
&
\ \ \ \ \
+
(\tfrac{2}{3}F_{1,3,0}^2+\tfrac{4}{3}F_{4,0,0}F_{1,3,0}+6F_{3,1,0}+6F_{0,4,0})x^2y^3
+
(\tfrac{19}{9}F_{0,4,0}F_{1,3,0}
+
\\
&
\ \ \ \ \
-\tfrac{2}{9}F_{3,1,0}F_{1,3,0}+\tfrac{28}{9}F_{0,4,0}F_{4,0,0}-9)xy^4
+
(\tfrac{6}{5}F_{1,3,0}+\tfrac{1}{15}F_{3,1,0}F_{0,4,0}+\tfrac{4}{3}F_{0,4,0}^2)y^5
+
\cdots.
\endaligned\right.
\]
\[
\def\arraystretch{1.25}
\begin{array}{llll}
e_1
&
:=
-(-1+\tfrac{1}{3}xF_{1,3,0}+\tfrac{8}{3}xF_{4,0,0}-9u+\tfrac{4}{9}uF_{0,4,0}F_{1,3,0}-\tfrac{2}{9}uF_{3,1,0}F_{1,3,0}+\tfrac{16}{9}uF_{0,4,0}F_{4,0,0})\partial_x
+
\\
&
\ \ \ \ \
-
(3x+\tfrac{4}{3}yF_{4,0,0}+\tfrac{2}{3}yF_{1,3,0}+3uF_{3,1,0})\partial_y
-
(F_{1,3,0}u+4F_{4,0,0}u-y)\partial_u,
,
\\
e_2
&
:= 
-(\tfrac{4}{3}xF_{0,4,0}+\tfrac{2}{3}xF_{3,1,0}+3y+3F_{1,3,0}u)\partial_x
-
(-1+\tfrac{1}{3}yF_{3,1,0}+\tfrac{8}{3}yF_{0,4,0}-9u+\tfrac{4}{9}uF_{0,4,0}F_{1,3,0}
+
\\
&
\ \ \ \ \
-\tfrac{2}{9}uF_{3,1,0}F_{1,3,0}+\tfrac{16}{9}uF_{0,4,0}F_{4,0,0})\partial_y
-(4F_{0,4,0}u+F_{3,1,0}u-x)\partial_u,
\\
e_3 
&
:= \partial_z,
\\
e_4 & := x\partial_z,
\\
e_5 & := y\partial_z,
\\
e_6 & := z\partial_z,
\\
e_7 & := u\partial_z.
\end{array}
\]
Gröbner basis generators of 
moduli space core algebraic variety in 
$\R^4 \ni (F_{0,4,0},F_{1,3,0},F_{3,1,0},F_{4,0,0})$:
\[
\aligned
B_1 
&
:= F_{0,4,0}F_{1,3,0}-F_{3,1,0}F_{4,0,0},
\\
B_2 & := -2F_{0,4,0}F_{3,1,0}^2+F_{1,3,0}^3+2F_{1,3,0}^2F_{4,0,0}-F_{3,1,0}^3,
\\
B_3 & := 2F_{0,4,0}^2F_{3,1,0}+F_{0,4,0}F_{3,1,0}^2-F_{1,3,0}^2F_{4,0,0}-2F_{1,3,0}F_{4,0,0}^2,
\\
B_4 
&
:= -32F_{0,4,0}F_{4,0,0}^2+F_{1,3,0}^2F_{3,1,0}+2F_{1,3,0}F_{3,1,0}F_{4,0,0}-16F_{3,1,0}F_{4,0,0}^2+36F_{0,4,0}F_{3,1,0}+18F_{3,1,0}^2
+
\\
&
\ \ \ \ \
+81F_{1,3,0}+162F_{4,0,0},
\\
B_5 &
:= -32F_{0,4,0}^2F_{4,0,0}-16F_{0,4,0}F_{3,1,0}F_{4,0,0}+F_{1,3,0}F_{3,1,0}^2+2F_{3,1,0}^2F_{4,0,0}+18F_{1,3,0}^2+36F_{1,3,0}F_{4,0,0}
+
\\
&
\ \ \ \ \
+162F_{0,4,0}+81F_{3,1,0},
\\
B_6 
&
:= 32F_{0,4,0}^3F_{4,0,0}-10F_{0,4,0}F_{3,1,0}^2F_{4,0,0}+8F_{1,3,0}^2F_{4,0,0}^2+16F_{1,3,0}F_{4,0,0}^3-F_{3,1,0}^3F_{4,0,0}
+
\\
&
\ \ \ \ \
-18F_{1,3,0}F_{3,1,0}F_{4,0,0}-36F_{3,1,0}F_{4,0,0}^2-162F_{0,4,0}^2-81F_{0,4,0}F_{3,1,0},
\\
B_7 &
:= 2F_{0,4,0}F_{3,1,0}^3-16F_{1,3,0}F_{3,1,0}F_{4,0,0}^2+F_{3,1,0}^4-32F_{3,1,0}F_{4,0,0}^3+576F_{0,4,0}^2F_{4,0,0}
+
\\
&
\ \ \ \ \
+288F_{0,4,0}F_{3,1,0}F_{4,0,0}-243F_{1,3,0}^2-486F_{1,3,0}F_{4,0,0}-2916F_{0,4,0}-1458F_{3,1,0},
\endaligned
\]

\[
\footnotesize
\def\arraystretch{1.25}
\begin{array}{c|ccccccc}
{} & e_1 & e_2 & e_3 & e_4 & e_5 & e_6 & e_7
\\
\hline
e_1 & 
0 &\rotatebox[origin=c]{90}{
	\begin{tabular}{p{4cm}}
$(-\tfrac{4}{3}F_{0,4,0}-\tfrac{2}{3}F_{3,1,0})e_1+(\tfrac{2}{3}F_{1,3,0}+\tfrac{4}{3}F_{4,0,0})e_2$\end{tabular}}
& 0 &\rotatebox[origin=c]{90}{
	\begin{tabular}{p{4cm}}
$e_3+(-\tfrac{1}{3}F_{1,3,0}-\tfrac{8}{3}F_{4,0,0})e_4+(9-\tfrac{4}{9}F_{0,4,0}F_{1,3,0}+\tfrac{2}{9}F_{3,1,0}F_{1,3,0}-\tfrac{16}{9}F_{0,4,0}F_{4,0,0})e_7$\end{tabular}}
&\rotatebox[origin=c]{90}{
	\begin{tabular}{p{4cm}}
$-3e_4+(-\tfrac{4}{3}F_{4,0,0}-\tfrac{2}{3}F_{1,3,0})e_5-3F_{3,1,0}e_7$\end{tabular}}
& 0 &\rotatebox[origin=c]{90}{
	\begin{tabular}{p{4cm}}
$e_5+(-F_{1,3,0}-4F_{4,0,0})e_7$\end{tabular}}
\\
e_2 &\rotatebox[origin=c]{90}{
	\begin{tabular}{p{4cm}}
$-(-\tfrac{4}{3}F_{0,4,0}-\tfrac{2}{3}F_{3,1,0})e_1-(\tfrac{2}{3}F_{1,3,0}+\tfrac{4}{3}F_{4,0,0})e_2$\end{tabular}}
& 0 & 0 &\rotatebox[origin=c]{90}{
	\begin{tabular}{p{4cm}}
$(-\tfrac{4}{3}F_{0,4,0}-\tfrac{2}{3}F_{3,1,0})e_4-3e_5-3F_{1,3,0}e_7$\end{tabular}}
&\rotatebox[origin=c]{90}{
	\begin{tabular}{p{4cm}}
$e_3+(-\tfrac{1}{3}F_{3,1,0}-\tfrac{8}{3}F_{0,4,0})e_5+(9-\tfrac{4}{9}F_{0,4,0}F_{1,3,0}+\tfrac{2}{9}F_{3,1,0}F_{1,3,0}-\tfrac{16}{9}F_{0,4,0}F_{4,0,0})e_7$\end{tabular}}
& 0 &
\rotatebox[origin=c]{90}{
	\begin{tabular}{p{4cm}}
	$e_4+(-4F_{0,4,0}-F_{3,1,0})e_7$\end{tabular}}
\\
e_3 &
0 & 0 & 0 & 0 & 0 & e_3 & 0
\\
e_4 &\rotatebox[origin=c]{90}{
	\begin{tabular}{p{4cm}}
$-e_3-(-\tfrac{1}{3}F_{1,3,0}-\tfrac{8}{3}F_{4,0,0})e_4-(9-\tfrac{4}{9}F_{0,4,0}F_{1,3,0}+\tfrac{2}{9}F_{3,1,0}F_{1,3,0}-\tfrac{16}{9}F_{0,4,0}F_{4,0,0})e_7$\end{tabular}}
&\rotatebox[origin=c]{90}{
	\begin{tabular}{p{4cm}}
$-(-\tfrac{4}{3}F_{0,4,0}-\tfrac{2}{3}F_{3,1,0})e_4+3e_5+3F_{1,3,0}e_7$\end{tabular}}
& 0 & 0 & 0 & e_4 & 0
\\
e_5 &\rotatebox[origin=c]{90}{
	\begin{tabular}{p{4cm}}
$3e_4-(-\tfrac{4}{3}F_{4,0,0}-\tfrac{2}{3}F_{1,3,0})e_5+3F_{3,1,0}e_7$\end{tabular}}
&\rotatebox[origin=c]{90}{
	\begin{tabular}{p{4cm}}
$-e_3-(-\tfrac{1}{3}F_{3,1,0}-\tfrac{8}{3}F_{0,4,0}e_5-(9-\tfrac{4}{9}F_{0,4,0}F_{1,3,0}+\tfrac{2}{9}F_{3,1,0}F_{1,3,0}-\tfrac{16}{9}F_{0,4,0}F_{4,0,0})e_7$\end{tabular}}
& 0 & 0 & 0 & e_5 & 0
\\
e_6 &
0 & 0 & -e_3 & -e_4 & -e_5 & 0 & -e_7
\\
e_7 &
\rotatebox[origin=c]{90}{
	\begin{tabular}{p{4cm}}
$-e_5-(-F_{1,3,0}-4F_{4,0,0})e_7$\end{tabular}}
&
\rotatebox[origin=c]{90}{
	\begin{tabular}{p{4cm}}
	$-e_4-(-4F_{0,4,0}-F_{3,1,0})e_7$\end{tabular}}
& 0 & 0 & 0 & e_7 & 0
\end{array}
\]


\[
\!\!\!\!\!\!\!\!\!
\text{\bf Model 2a3d4d5a}
\ \ \ \ \
\left\{
\aligned
u
&
\,=\,
xy+y^2z+y^4+F_{0,4,1}y^4z+F_{1,4,0}xy^4+y^5
+
\big(-\tfrac{3}{2}F_{1,4,0}^2+\tfrac{3}{2}F_{1,4,0}(-\tfrac{4}{5}F_{0,4,1}^2
+
\\
&
\ \ \ \ \
-\tfrac{24}{5}F_{0,4,1}F_{1,4,0}+\tfrac{18}{5}F_{0,6,0}F_{1,4,0}-3F_{1,4,0})\big)x^2y^4
+
(-\tfrac{4}{5}F_{0,4,1}^2-\tfrac{24}{5}F_{0,4,1}F_{1,4,0}
+
\\
&
\ \ \ \ \
+\tfrac{18}{5}F_{0,6,0}F_{1,4,0}-3F_{1,4,0})xy^5
+
F_{0,6,0}y^6
+
\big(-2F_{0,4,1}F_{1,4,0}+2F_{0,4,1}(-\tfrac{4}{5}F_{0,4,1}^2
+
\\
&
\ \ \ \ \
-\tfrac{24}{5}F_{0,4,1}F_{1,4,0}+\tfrac{18}{5}F_{0,6,0}F_{1,4,0}-3F_{1,4,0})\big)xy^4z
+
(-\tfrac{16}{5}F_{0,4,1}^2+\tfrac{12}{5}F_{0,4,1}F_{0,6,0}
+
\\
&
\ \ \ \ \
-2F_{0,4,1}+\tfrac{6}{5}F_{1,4,0})y^5z
+
\cdots.
\endaligned\right.
\]
\[
\!\!\!\!\!\!\!\!
\def\arraystretch{1.25}
\aligned
e_1 & := 
(1+13xF_{1,4,0}+\tfrac{12}{5}xF_{0,4,1}^2+\tfrac{72}{5}xF_{0,4,1}F_{1,4,0}-\tfrac{54}{5}xF_{0,6,0}F_{1,4,0})\partial_x
+
(4yF_{1,4,0}+\tfrac{4}{5}yF_{0,4,1}^2
+
\\
&
\ \ \ \ \
+\tfrac{24}{5}yF_{0,4,1}F_{1,4,0}-\tfrac{18}{5}yF_{0,6,0}F_{1,4,0})\partial_y
+
(9zF_{1,4,0}+\tfrac{8}{5}zF_{0,4,1}^2+\tfrac{48}{5}zF_{0,4,1}F_{1,4,0}-\tfrac{36}{5}zF_{0,6,0}F_{1,4,0})\partial_z
+
\\
&
\ \ \ \ \
+
(y+17uF_{1,4,0}+\tfrac{16}{5}uF_{0,4,1}^2+\tfrac{96}{5}uF_{0,4,1}F_{1,4,0}-\tfrac{72}{5}uF_{0,6,0}F_{1,4,0})\partial_u,
\\
e_2 & := -(4F_{0,4,1}u-24F_{0,4,1}x+18F_{0,6,0}x-20x+2z)\partial_x
+
(8F_{0,4,1}y-6F_{0,6,0}y+5y+1)\partial_y
+
\\
&
\ \ \ \ \
+
(4F_{0,4,1}x+16F_{0,4,1}z-12F_{0,6,0}z-4F_{1,4,0}u-4y+15z)\partial_z
+
(32F_{0,4,1}u-24F_{0,6,0}u+25u+x)\partial_u,
\\
e_3 & := (10F_{0,4,1}x+\tfrac{48}{5}xF_{0,4,1}^2-\tfrac{36}{5}xF_{0,4,1}F_{0,6,0}-\tfrac{3}{5}xF_{1,4,0}-y)\partial_x
+
(3yF_{0,4,1}+\tfrac{16}{5}yF_{0,4,1}^2
+
\\
&
\ \ \ \ \
-\tfrac{12}{5}yF_{0,4,1}F_{0,6,0}-\tfrac{1}{5}yF_{1,4,0})\partial_y
+
(1+7zF_{0,4,1}+\tfrac{32}{5}zF_{0,4,1}^2-\tfrac{24}{5}zF_{0,4,1}F_{0,6,0}-\tfrac{2}{5}zF_{1,4,0})\partial_z
+
\\
&
\ \ \ \ \
+
(13F_{0,4,1}u+\tfrac{64}{5}uF_{0,4,1}^2-\tfrac{48}{5}uF_{0,4,1}F_{0,6,0}-\tfrac{4}{5}uF_{1,4,0})\partial_u.
\endaligned
\]
Gröbner basis generators of 
moduli space core algebraic variety in 
$\R^3 \ni (F_{0,4,1}, F_{0,6,0}, F_{1,4,0})$:
\[
\aligned
B_1 & := 8F_{0,4,1}^3-5F_{0,4,1}F_{1,4,0}+3F_{1,4,0}^2,
\\
B_2 & := F_{1,4,0}(2F_{0,4,1}^2+12F_{0,4,1}F_{1,4,0}-9F_{0,6,0}F_{1,4,0}+10F_{1,4,0}),
\\
B_3 & := 4F_{0,4,1}^2F_{0,6,0}-5F_{0,4,1}^2-3F_{0,4,1}F_{1,4,0}+2F_{1,4,0}^2,
\\
B_4 & := F_{1,4,0}(4F_{0,4,1}F_{0,6,0}+32F_{0,4,1}F_{1,4,0}-24F_{0,6,0}F_{1,4,0}-5F_{0,4,1}+27F_{1,4,0}),
\\
B_5 & := F_{1,4,0}(384F_{0,4,1}F_{1,4,0}+36F_{0,6,0}^2-288F_{0,6,0}F_{1,4,0}-6F_{0,4,1}-85F_{0,6,0}+328F_{1,4,0}+50),
\\
B_6 & := 12F_{0,4,1}F_{0,6,0}^2+6F_{0,4,1}^2-35F_{0,4,1}F_{0,6,0}-16F_{0,4,1}F_{1,4,0}+4F_{0,6,0}F_{1,4,0}+8F_{1,4,0}^2+25F_{0,4,1}
+
\\
&
\ \ \ \ \-5F_{1,4,0}.
\endaligned
\]
\[
\footnotesize
\def\arraystretch{1.25}
\begin{array}{c|ccc}
{} & e_1 & e_2 & e_3
\\
\hline
e_1 & 
0 &\rotatebox[origin=c]{90}{
	\begin{tabular}{p{6cm}}
$(24F_{0,4,1}-18F_{0,6,0}+20)e_1+(-4F_{1,4,0}-\tfrac{4}{5}F_{0,4,1}^2-\tfrac{24}{5}F_{0,4,1}F_{1,4,0}+\tfrac{18}{5}F_{0,6,0}F_{1,4,0})e_2+4F_{0,4,1}e_3$
\end{tabular}}
&
\rotatebox[origin=c]{90}{
	\begin{tabular}{p{6cm}}
$(10F_{0,4,1}+\tfrac{48}{5}F_{0,4,1}^2-\tfrac{36}{5}F_{0,4,1}F_{0,6,0}-\tfrac{3}{5}F_{1,4,0})e_1+(-\tfrac{48}{5}F_{0,4,1}F_{1,4,0}+\tfrac{36}{5}F_{0,6,0}F_{1,4,0}-\tfrac{8}{5}F_{0,4,1}^2-9F_{1,4,0})e_3$
\end{tabular}}
\\
e_2 & \rotatebox[origin=c]{90}{
	\begin{tabular}{p{6cm}}
$-(24F_{0,4,1}-18F_{0,6,0}+20)e_1-(-4F_{1,4,0}-\tfrac{4}{5}F_{0,4,1}^2-\tfrac{24}{5}F_{0,4,1}F_{1,4,0}+\tfrac{18}{5}F_{0,6,0}F_{1,4,0})e_2-4F_{0,4,1}e_3$
\end{tabular}}
&
0
&
\rotatebox[origin=c]{90}{
	\begin{tabular}{p{6cm}}
$e_1+(3F_{0,4,1}+\tfrac{16}{5}F_{0,4,1}^2-\tfrac{12}{5}F_{0,4,1}F_{0,6,0}-\tfrac{1}{5}F_{1,4,0})e_2+(-15-16F_{0,4,1}+12F_{0,6,0})e_3$
\end{tabular}}
\\
e_3 &
\rotatebox[origin=c]{90}{
	\begin{tabular}{p{6cm}}
$-(10F_{0,4,1}+\tfrac{48}{5}F_{0,4,1}^2-\tfrac{36}{5}F_{0,4,1}F_{0,6,0}-\tfrac{3}{5}F_{1,4,0})e_1-(-\tfrac{48}{5}F_{0,4,1}F_{1,4,0}+\tfrac{36}{5}F_{0,6,0}F_{1,4,0}-\tfrac{8}{5}F_{0,4,1}^2-9F_{1,4,0})e_3$
\end{tabular}}
&
\rotatebox[origin=c]{90}{
	\begin{tabular}{p{6cm}}
$-e_1-(3F_{0,4,1}+\tfrac{16}{5}F_{0,4,1}^2-\tfrac{12}{5}F_{0,4,1}F_{0,6,0}-\tfrac{1}{5}F_{1,4,0})e_2-(-15-16F_{0,4,1}+12F_{0,6,0})e_3$
\end{tabular}}
&
0
\end{array}
\]


\[
\text{\bf Model 2a3d4d5e}
\ \ \ \ \
\left\{
\aligned
u
&
\,=\,
y^4+y^2z+xy.
\endaligned\right.
\]
\[
\def\arraystretch{1.25}
\begin{array}{llll}
e_1
\,:=\,
\partial_x+y\partial_u
, &
e_2
\,:=\,
x\partial_u-2z\partial_x-4y\partial_z+\partial_y
, &
& 
\\
e_3
\,:=\,
-y\partial_x+\partial_z
, &
e_4
\,:=\,
4u\partial_u+3x\partial_x+y\partial_y+2z\partial_z
.
\end{array}
\]

\[
\footnotesize
\def\arraystretch{1.25}
\begin{array}{c|cccc}
{} & e_1 & e_2 & e_3 & e_4
\\
\hline
e_1 & 
0 & 0 & 0 & 3e_1
\\
e_2 &
0 & 0 & e_1 & e_2
\\
e_3 &
0 & -e_1 & 0 & 2e_3
\\
e_4 &
-3e_1 & -e_2 & -2e_3 & 0
\end{array}
\]


\[
\text{\bf Model 2a3d4e}
\ \ \ \ \
\left\{
\aligned
u
&
\,=\,
y^2z+xy.
\endaligned\right.
\]
\[
\def\arraystretch{1.25}
\begin{array}{llll}
e_1
\,:=\,
\partial_x+y\partial_u
, &
e_2
\,:=\,
-2z\partial_x+\partial_y+x\partial_u
, &
& 
\\
e_3
\,:=\,
-y\partial_x+\partial_z
, &
e_4
\,:=\,
x\partial_x+y\partial_y+2u\partial_u
, &
e_5
\,:=\,
x\partial_x+z\partial_z+u\partial_u
.
\end{array}
\]

\rotatebox[origin=c]{90}{
	\begin{tabular}{p{3cm}}
\end{tabular}}

\[
\footnotesize
\def\arraystretch{1.25}
\begin{array}{c|ccccc}
{} & e_1 & e_2 & e_3 & e_4 & e_5 
\\
\hline
e_1 & 
0 & 0 & 0 & e_1 & e_1
\\
e_2 &
0 & 0 & e_1 & e_2 & 0
\\
e_3 &
0 & -e_1 & 0 & 0 & e_3
\\
e_4 &
-e_1 & -e_2 & 0 & 0 & 0
\\
e_5 &
-e_1 & 0 & -e_3 & 0 & 0
\end{array}
\]


\[
\text{\bf Model 2a3f}
\ \ \ \ \
\left\{
\aligned
u
=
xyz^5+xyz^4+xyz^3+xyz^2+xyz+xy
+
\cdots
.
\endaligned\right.
\]
\[
\def\arraystretch{1.25}
\begin{array}{llll}
e_1
\,:=\,
-(z-1)\partial_x+y\partial_u
, &
e_2
\,:=\,
-(z-1)\partial_y+x\partial_u
, &
& 
\\
e_3
\,:=\,
-(z-1)\partial_z+u\partial_u
, &
e_4
\,:=\,
x\partial_x+u\partial_u
, &
e_5
\,:=\,
y\partial_y+u\partial_u
, &
e_6
\,:=\,
-u\partial_y+x\partial_z
,
\\
e_7
\,:=\,
-u\partial_x+y\partial_z
.
\end{array}
\]

\rotatebox[origin=c]{90}{
	\begin{tabular}{p{3cm}}
\end{tabular}}

\[
\footnotesize
\def\arraystretch{1.25}
\begin{array}{c|ccccccc}
{} & e_1 & e_2 & e_3 & e_4 & e_5 & e_6 & e_7
\\
\hline
e_1 & 
0 & 0 & e_1 & e_1 & 0 & e_3+e_4-e_5 & 0
\\
e_2 &
0 & 0 & e_2 & 0 & e_2 & 0 & e_3-e_4+e_5
\\
e_3 &
-e_1 & -e_2 & 0 & 0 & 0 & e_6 & e_7
\\
e_4 &
-e_1 & 0 & 0 & 0 & 0 & e_6 & 0
\\
e_5 &
0 & -e_2 & 0 & 0 & 0 & 0 & e_7
\\
e_6 &
-e_3-e_4+e_5 & 0 & -e_6 & -e_6 & 0 & 0 & 0
\\
e_7 &
0 & -e_3+e_4-e_5 & -e_7 & 0 & -e_7 & 0 & 0
\end{array}
\]


\[
\text{\bf Model 2a3g5b}
\ \ \ \ \
\left\{
\aligned
u
&
\,=\,
xy+xyz+y^3+xyz^2+2y^3z+F_{2,3,0}x^2y^3+xyz^3+y^5+3y^3z^2
+
\\
&
\ \ \ \ \
+F_{0,6,0}y^6+4F_{2,3,0}x^2y^3z+3F_{2,3,0}xy^5+xyz^4+4y^5z+4y^3z^3+F_{2,3,0}x^2y^5
+
\\
&
\ \ \ \ \
+(\tfrac{9}{7}F_{0,6,0}^2+\tfrac{15}{7}F_{2,3,0}+\tfrac{15}{7})y^7+5F_{0,6,0}y^6z+15F_{2,3,0}xy^5z+10F_{2,3,0}x^2y^3z^2
+
\\
&
\ \ \ \ \
+10y^5z^2+5y^3z^4+xyz^5
+
\cdots
.
\endaligned\right.
\]
\[
\def\arraystretch{1.25}
\aligned
e_1 & := -(z-1)\partial_x-F_{2,3,0}u\partial_y+F_{2,3,0}x\partial_z+y\partial_u,
\\
e_2 & := -(6F_{0,6,0}x-5u+3y)\partial_x-(3F_{0,6,0}y+z-1)\partial_y-(3F_{2,3,0}u+5y)\partial_z-(9F_{0,6,0}u-x)\partial_u,
\\
e_3 & := -x\partial_x-y\partial_y-(z-1)\partial_z-u\partial_u.
\endaligned
\]
Gröbner basis generator of 
moduli space core algebraic variety in 
$\R^2 \ni (F_{2,3,0},F_{0,6,0})$:
\[
B_1
:= 
F_{2,3,0}F_{0,6,0}.
\]
\[
\footnotesize
\def\arraystretch{1.25}
\begin{array}{c|ccc}
{} & e_1 & e_2 & e_3
\\
\hline
e_1 & 
0 & -6F_{0,6,0}e_1 & 0
\\
e_2 &
6F_{0,6,0}e_1 & 0 & 0
\\
e_3 &
0 & 0 & 0
\end{array}
\]

\rotatebox[origin=c]{90}{
	\begin{tabular}{p{3cm}}
\end{tabular}}


\[
\text{\bf Model 2a3g5c}
\ \ \ \ \
\left\{
\aligned
u
&
\,=\,
xy+xyz+y^3+xyz^2+2y^3z+F_{2,3,0}x^2y^3+xyz^3+y^5+3y^3z^2
+
\\
&
\ \ \ \ \
+F_{0,6,0}y^6+4F_{2,3,0}x^2y^3z+3F_{2,3,0}xy^5+xyz^4-4y^5z+4y^3z^3-F_{2,3,0}x^2y^5
+
\\
&
\ \ \ \ \
+(-\tfrac{9}{7}F_{0,6,0}^2+\tfrac{15}{7}F_{2,3,0}+\tfrac{15}{7})y^7+5F_{0,6,0}y^6z+15F_{2,3,0}xy^5z+10F_{2,3,0}x^2y^3z^2
+
\\
&
\ \ \ \ \
-10y^5z^2+5y^3z^4+xyz^5
+
\cdots
.
\endaligned\right.
\]
\[
\def\arraystretch{1.25}
\aligned
e_1 & := -(z-1)\partial_x-F_{2,3,0}u\partial_y+F_{2,3,0}x\partial_z+y\partial_u,
\\
e_2 & := (6F_{0,6,0}x-5u-3y)\partial_x+(3F_{0,6,0}y-z+1)\partial_y-(3F_{2,3,0}u-5y)\partial_z+(9F_{0,6,0}ux)\partial_u,
\\
e_3 & := -x\partial_x-y\partial_y-(z-1)\partial_z-u\partial_u.
\endaligned
\]
Gröbner basis generator of 
moduli space core algebraic variety in 
$\R^2 \ni (F_{2,3,0},F_{0,6,0})$:
\[
B_1 
:= 
F_{2,3,0}F_{0,6,0}.
\]
\[
\footnotesize
\def\arraystretch{1.25}
\begin{array}{c|ccc}
{} & e_1 & e_2 & e_3
\\
\hline
e_1 & 
0 & -6F_{0,6,0}e_1 & 0
\\
e_2 &
6F_{0,6,0}e_1 & 0 & 0
\\
e_3 &
0 & 0 & 0
\end{array}
\]

\rotatebox[origin=c]{90}{
	\begin{tabular}{p{3cm}}
\end{tabular}}


\[
\text{\bf Model 2a3g5d}
\ \ \ \ \
\left\{
\aligned
u
&
\,=\,
xy+xyz+y^3+xyz^2+2y^3z+x^2y^3+xyz^3+3y^3z^2+4x^2y^3z+3xy^5
+
\\
&
\ \ \ \ \
+xyz^4+4y^3z^3+\tfrac{15}{7}y^7+15xy^5z+10x^2y^3z^2+5y^3z^4+xyz^5
+
\cdots.
\endaligned\right.
\]
\[
\def\arraystretch{1.25}
\aligned
e_1 & := -(z-1)\partial_x-u\partial_y+x\partial_z+y\partial_u,
\\
e_2 & := -3y\partial_x-(z-1)\partial_y-3u\partial_z+x\partial_u,
\\
e_3 & := -x\partial_x-y\partial_y-(z-1)\partial_z-u\partial_u.
\endaligned
\]

\[
\footnotesize
\def\arraystretch{1.25}
\begin{array}{c|ccc}
{} & e_1 & e_2 & e_3
\\
\hline
e_1 & 
0 & 0 & 0
\\
e_2 &
0 & 0 & 0
\\
e_3 &
0 & 0 & 0
\end{array}
\]


\[
\text{\bf Model 2a3g5e}
\ \ \ \ \
\left\{
\aligned
u
&
\,=\,
xy+xyz+y^3+xyz^2+2y^3z-x^2y^3+xyz^3+3y^3z^2-4x^2y^3z-3xy^5
+
\\
&
\ \ \ \ \
+xyz^4+4y^3z^3-\tfrac{15}{7}y^7-15xy^5z-10x^2y^3z^2+5y^3z^4+xyz^5
+
\cdots
.
\endaligned\right.
\]
\[
\def\arraystretch{1.25}
\aligned
e_1 & := -(z-1)\partial_x+u\partial_y-x\partial_z+y\partial_u,
\\
e_2 & := -3y\partial_x-(z-1)\partial_y+3u\partial_z+x\partial_u,
\\
e_3 & := -x\partial_x-y\partial_y-(z-1)\partial_z-u\partial_u.
\endaligned
\]

\[
\footnotesize
\def\arraystretch{1.25}
\begin{array}{c|ccc}
{} & e_1 & e_2 & e_3
\\
\hline
e_1 & 
0 & 0 & 0
\\
e_2 &
0 & 0 & 0
\\
e_3 &
0 & 0 & 0
\end{array}
\]

\[
\text{\bf Model 2a3g5f}
\ \ \ \ \
\left\{
\aligned
u
&
\,=\,
xyz^5+5y^3z^4+xyz^4+4y^3z^3+xyz^3+3y^3z^2+xyz^2+2y^3z+xyz
+
\\
&
\ \ \ \ \
+y^3+xy
+
\cdots
.
\endaligned\right.
\]
\[
\def\arraystretch{1.25}
\begin{array}{llll}
e_1
\,:=\,
-(z-1)\partial_x+y\partial_u
, &
e_2
\,:=\,
-3y\partial_x-(z-1)\partial_y+x\partial_u
, &
& 
\\
e_3
\,:=\,
x\partial_x-(z-1)\partial_z+2u\partial_u
, &
e_4
\,:=\,
2x\partial_x+y\partial_y+3u\partial_u
.
\end{array}
\]

\[
\footnotesize
\def\arraystretch{1.25}
\begin{array}{c|cccc}
{} & e_1 & e_2 & e_3 & e_4
\\
\hline
e_1 & 
0 & 0 & 2e_1 & 2e_1
\\
e_2 &
0 & 0 & e_2 & e_2
\\
e_3 &
-2e_1 & -e_2 & 0 & 0
\\
e_4 &
-2e_1 & -e_2 & 0 & 0
\end{array}
\]


\[
\text{\bf Model 2a3h}
\ \ \ \ \
\left\{
\aligned
u
&
\,=\,
xy+2x^3z+2y^3z+3x^3z^2+3y^3z^2+4x^3z^3+4y^3z^3+x^3+xyz
+
\\
&
\ \ \ \ \
+F_{1,3,0}xy^3+F_{0,4,0}y^4+xyz^2+F_{5,0,0}x^5
+
(\tfrac{2}{9}F_{0,4,0}F_{1,3,0}-9)x^4y
+
\\
&
\ \ \ \ \
+
F_{3,2,0}x^3y^2
+
(\tfrac{2}{3}F_{1,3,0}^2-5F_{5,0,0}+6F_{0,4,0})x^2y^3
+
(-9+\tfrac{11}{9}F_{0,4,0}F_{1,3,0})xy^4
+
\\
&
\ \ \ \ \
+
(\tfrac{12}{5}F_{1,3,0}-\tfrac{1}{5}F_{3,2,0}+\tfrac{4}{3}F_{0,4,0}^2)y^5
+
3F_{1,3,0}xy^3z
+
xyz^3
+
(\tfrac{1}{3}F_{0,4,0}F_{1,3,0}-9
+
\\
&
\ \ \ \ \
+\tfrac{1}{9}F_{5,0,0}F_{1,3,0})x^6
+
\big(\tfrac{1}{5}F_{1,3,0}(\tfrac{2}{9}F_{0,4,0}F_{1,3,0}-9)-\tfrac{27}{5}F_{1,3,0}+3F_{3,2,0}\big)x^5y
+
\\
&
\ \ \ \ \
+
(\tfrac{1}{3}F_{3,2,0}F_{1,3,0}-15F_{5,0,0}+\tfrac{3}{2}F_{1,3,0}^2+\tfrac{27}{2}F_{0,4,0})x^4y^2
+
(-45+\tfrac{26}{3}F_{0,4,0}F_{1,3,0}
+
\\
&
\ \ \ \ \
-\tfrac{40}{9}F_{5,0,0}F_{1,3,0}+\tfrac{10}{27}F_{1,3,0}^3)x^3y^3
+
\big(-\tfrac{81}{2}F_{1,3,0}+3F_{3,2,0}-3F_{1,3,0}(\tfrac{2}{9}F_{0,4,0}F_{1,3,0}
+
\\
&
\ \ \ \ \
-9)+\tfrac{5}{3}F_{0,4,0}F_{1,3,0}^2-5F_{0,4,0}F_{5,0,0}+10F_{0,4,0}^2\big)x^2y^4
+
\big(-15F_{5,0,0}-45F_{0,4,0}
+
\\
&
\ \ \ \ \
+\tfrac{28}{5}F_{1,3,0}^2-\tfrac{7}{15}F_{3,2,0}F_{1,3,0}+\tfrac{28}{9}F_{0,4,0}^2F_{1,3,0}-4F_{0,4,0}(\tfrac{2}{9}F_{0,4,0}F_{1,3,0}-9)\big)xy^5
+
\\
&
\ \ \ \ \
+
(-9+\tfrac{101}{15}F_{0,4,0}F_{1,3,0}-\tfrac{29}{45}F_{3,2,0}F_{0,4,0}+\tfrac{56}{27}F_{0,4,0}^3)y^6
+
(\tfrac{8}{9}F_{0,4,0}F_{1,3,0}
+
\\
&
\ \ \ \ \
-36)x^4yz
+
4F_{3,2,0}x^3y^2z
+
(\tfrac{8}{3}F_{1,3,0}^2-20F_{5,0,0}+24F_{0,4,0})x^2y^3z
+
(-36
+
\\
&
\ \ \ \ \
+\tfrac{44}{9}F_{0,4,0}F_{1,3,0})xy^4z
+
(\tfrac{48}{5}F_{1,3,0}-\tfrac{4}{5}F_{3,2,0}+\tfrac{16}{3}F_{0,4,0}^2)y^5z
+
6F_{1,3,0}xy^3z^2
+
\\
&
\ \ \ \ \
+xyz^4+3F_{0,4,0}y^4z+4F_{5,0,0}x^5z+6F_{0,4,0}y^4z^2+y^3
+
\cdots
.
\endaligned\right.
\]
\[
\def\arraystretch{1.25}
\aligned
e_1 & := 
-(1+\tfrac{1}{3}xF_{1,3,0}+z-9u+\tfrac{4}{9}uF_{0,4,0}F_{1,3,0})\partial_x
-
(3x+2yF_{1,3,0}\tfrac{1}{3}-5uF_{5,0,0})\partial_y
+
\\
&
\ \ \ \ \
+
(-5xF_{5,0,0}+\tfrac{4}{9}yF_{0,4,0}F_{1,3,0}+9F_{1,3,0}u-3F_{3,2,0}u)\partial_z
-
(F_{1,3,0}u-y)\partial_u,
\\
e_2 & := -(\tfrac{4}{3}xF_{0,4,0}+3y-3F_{1,3,0}u+F_{3,2,0}u)\partial_x
+
(1-\tfrac{8}{3}yF_{0,4,0}-z+9u+\tfrac{2}{9}uF_{0,4,0}F_{1,3,0})\partial_y
+
\\
&
\ \ \ \ \
-
(\tfrac{2}{9}xF_{1,3,0}F_{0,4,0}+6yF_{1,3,0}-yF_{3,2,0}+2F_{1,3,0}^2u+18F_{0,4,0}u-15uF_{5,0,0})\partial_z
-
(4F_{0,4,0}u-x)\partial_u,
\\
e_3 & := -x\partial_x-y\partial_y-(-1+z)\partial_z-u\partial_u.
\endaligned
\]
Gröbner basis generators of 
moduli space core algebraic variety in 
$\R^4 \ni (F_{0,4,0},F_{1,3,0},F_{3,2,0},F_{5,0,0})$:
\[
\aligned
B_1 & := 2F_{0,4,0}F_{1,3,0}^2-120F_{0,4,0}F_{5,0,0}+81F_{1,3,0},
\\
B_2 & := 8F_{0,4,0}^2F_{1,3,0}+27F_{1,3,0}^2-6F_{1,3,0}F_{3,2,0}-81F_{0,4,0},
\\
B_3 & := 5F_{1,3,0}^3+108F_{0,4,0}F_{1,3,0}-24F_{0,4,0}F_{3,2,0}-60F_{1,3,0}F_{5,0,0},
\\
B_4 & := -800F_{0,4,0}^2F_{5,0,0}+10F_{1,3,0}^2F_{3,2,0}+1647F_{0,4,0}F_{1,3,0}-216F_{0,4,0}F_{3,2,0}-540F_{1,3,0}F_{5,0,0},
\\
B_5 & := 16F_{0,4,0}^2F_{3,2,0}-160F_{0,4,0}F_{1,3,0}F_{5,0,0}+378F_{1,3,0}^2-54F_{1,3,0}F_{3,2,0}-729F_{0,4,0},
\\
B_6 & := 288F_{0,4,0}F_{1,3,0}F_{3,2,0}-32F_{0,4,0}F_{3,2,0}^2-1440F_{1,3,0}^2F_{5,0,0}+320F_{1,3,0}F_{3,2,0}F_{5,0,0}
+
\\
&
\ \ \ \ \
-60480F_{0,4,0}F_{5,0,0}+44469F_{1,3,0},
\\
B_7 & := -14400F_{0,4,0}^2F_{5,0,0}-16000F_{0,4,0}F_{5,0,0}^2+20F_{1,3,0}F_{3,2,0}^2+14823F_{0,4,0}F_{1,3,0}-594F_{0,4,0}F_{3,2,0}
+
\\
&
\ \ \ \ \
+8640F_{1,3,0}F_{5,0,0},
\\
B_8 & := 1280F_{0,4,0}^3F_{5,0,0}+4320F_{0,4,0}F_{1,3,0}F_{5,0,0}-960F_{0,4,0}F_{3,2,0}F_{5,0,0}+729F_{1,3,0}^2-162F_{1,3,0}F_{3,2,0}
+
\\
&
\ \ \ \ \
-10935F_{0,4,0},
\\
B_9 & := -345600F_{0,4,0}^2F_{5,0,0}^2+32F_{0,4,0}F_{3,2,0}^3-256000F_{0,4,0}F_{5,0,0}^3+738720F_{0,4,0}F_{1,3,0}F_{5,0,0}
+
\\
&
\ \ \ \ \
-73440F_{0,4,0}F_{3,2,0}F_{5,0,0}+216000F_{1,3,0}F_{5,0,0}^2-400221F_{1,3,0}^2+60507F_{1,3,0}F_{3,2,0}.
\endaligned
\]
\[
\footnotesize
\def\arraystretch{1.25}
\begin{array}{c|ccc}
{} & e_1 & e_2 & e_3
\\
\hline
e_1 & 
0 & \rotatebox[origin=c]{0}{
	\begin{tabular}{p{7cm}}
	$-\tfrac{4}{3}F_{0,4,0}e_1+\tfrac{2}{3}F_{1,3,0}e_2-\tfrac{2}{3}F_{0,4,0}F_{1,3,0}e_3$
	\end{tabular}}
	& 0
\\
e_2 &
\rotatebox[origin=c]{0}{
	\begin{tabular}{p{7cm}}
$\tfrac{4}{3}F_{0,4,0}e_1-\tfrac{2}{3}F_{1,3,0}e_2+\tfrac{2}{3}F_{0,4,0}F_{1,3,0}e_3$
\end{tabular}}
 & 0 & 0
\\
e_3 &
0 & 0 & 0
\end{array}
\]


\[
\text{\bf Model 2a3i4a}
\ \ \ \ \
\left\{
\aligned
u
&
\,=\,
10x^3y^3z+10x^2y^4z-20x^2y^3z^2-10x^2y^2z^3-15xy^4z^2-20xy^3z^3
+
\\
&
\ \ \ \ \
+10xy^2z^4+xyz^5-5y^4z^3+20y^3z^4+5y^2z^5+2x^3y^3-4x^2y^3z
+
\\
&
\ \ \ \ \
-6x^2y^2z^2-10xy^3z^2+6xy^2z^3+xyz^4-2y^4z^2+8y^3z^3+4y^2z^4
+
\\
&
\ \ \ \ \
-3x^2y^2z-3xy^3z+3xy^2z^2+xyz^3+2y^3z^2+3y^2z^3-x^2y^2+xy^2z
+
\\
&
\ \ \ \ \
+xyz^2+2y^2z^2+xyz+y^2z+xy
+
\cdots
.
\endaligned\right.
\]
\[
\def\arraystretch{1.25}
\begin{array}{ll}
e_1
&
\,:=\,
(3u-z+1)\partial_x+(u-y)\partial_y-x\partial_z-(u-y)\partial_u
,
\\
e_2
&
\,:=\,
-(x+2z)\partial_x+(2u+y-z+1)\partial_y-2z\partial_z+x\partial_u
, 
\\
e_3
&
\,:=\,
-(3u-x+y)\partial_x-2u\partial_y+(u+2x-z+1)\partial_z+2u\partial_u
.
\end{array}
\]

\[
\footnotesize
\def\arraystretch{1.25}
\begin{array}{c|ccc}
{} & e_1 & e_2 & e_3
\\
\hline
e_1 & 
0 & -e_1+e_2 & 2e_1+2e_3
\\
e_2 &
e_1-e_2 & 0 & e_1+e_2+2e_3
\\
e_3 &
-2e_1-2e_3 & -e_1-e_2-2e_3 & 0
\end{array}
\]


\[
\text{\bf Model 2a3i4b}
\ \ \ \ \
\left\{
\aligned
u
&
\,=\,
xy+xyz+y^2z+xyz^2+y^4+2y^2z^2+F_{0,5,0}y^5+xyz^3+4y^4z
+
\\
&
\ \ \ \ \
+3y^2z^3+(2+\tfrac{5}{4}F_{0,5,0}^2)y^6+5F_{0,5,0}y^5z+10y^4z^2+xyz^4+4y^2z^4
+
\\
&
\ \ \ \ \
+(\tfrac{36}{7}F_{0,5,0}+\tfrac{25}{14}F_{0,5,0}^3)y^7+(12+\tfrac{15}{2}F_{0,5,0}^2)y^6z+15F_{0,5,0}y^5z^2+20y^4z^3
+
\\
&
\ \ \ \ \
+xyz^5+5y^2z^5
+
\cdots.
,
\endaligned\right.
\]
for any value of $F_{0,5,0}$.
\[
\def\arraystretch{1.25}
\begin{array}{ll}
e_1
&
\,:=\,
-(-1+z)\partial_x+y\partial_u
,
\\
e_2
&
\,:=\,
-(\tfrac{5}{2}F_{0,5,0}x+2z-4u)\partial_x-(-1+\tfrac{5}{2}yF_{0,5,0}+z)\partial_y-4y\partial_z-(5F_{0,5,0}u-x)\partial_u
, 
\\
e_3
&
\,:=\,
-y\partial_x-y\partial_y-(-1+z)\partial_z
.
\end{array}
\]

\[
\footnotesize
\def\arraystretch{1.25}
\begin{array}{c|ccc}
{} & e_1 & e_2 & e_3
\\
\hline
e_1 & 
0 & -\tfrac{5}{2}F_{0,5,0}e_1 & e_1
\\
e_2 &
\tfrac{5}{2}F_{0,5,0}e_1 & 0 & e_1
\\
e_3 &
 -e_1 & -e_1 & 0
\end{array}
\]


\[
\text{\bf Model 2a3i4c}
\ \ \ \ \
\left\{
\aligned
u
&
\,=\,
xy+xyz+y^2z+xyz^2-y^4+2y^2z^2+F_{0,5,0}y^5+xyz^3-4y^4z
+
\\
&
\ \ \ \ \
+3y^2z^3+(2-\tfrac{5}{4}F_{0,5,0}^2)y^6+5F_{0,5,0}y^5z-10y^4z^2+xyz^4+4y^2z^4
+
\\
&
\ \ \ \ \
+(-\tfrac{36}{7}F_{0,5,0}+\tfrac{25}{14}F_{0,5,0}^3)y^7+(12-\tfrac{15}{2}F_{0,5,0}^2)y^6z+15F_{0,5,0}y^5z^2-20y^4z^3
+
\\
&
\ \ \ \ \
+xyz^5+5y^2z^5
+
\cdots
,
\endaligned\right.
\]
for any value of $F_{0,5,0}$.
\[
\def\arraystretch{1.25}
\begin{array}{ll}
e_1
&
\,:=\,
-(-1+z)\partial_x+y\partial_u
,
\\
e_2
&
\,:=\,
(\tfrac{5}{2}F_{0,5,0}x-2z-4u)\partial_x+(1+\tfrac{5}{2}yF_{0,5,0}-z)\partial_y+4y\partial_z+(5F_{0,5,0}u+x)\partial_u
, 
\\
e_3
&
\,:=\,
-y\partial_x-y\partial_y-(-1+z)\partial_z
.
\end{array}
\]

\[
\footnotesize
\def\arraystretch{1.25}
\begin{array}{c|ccc}
{} & e_1 & e_2 & e_3
\\
\hline
e_1 & 
0 & \tfrac{5}{2}F_{0,5,0}e_1 & e_1
\\
e_2 &
-\tfrac{5}{2}F_{0,5,0}e_1 & 0 & e_1
\\
e_3 &
-e_1 & -e_1 & 0
\end{array}
\]


\[
\text{\bf Model 2a3i4f}
\ \ \ \ \
\left\{
\aligned
u
&
\,=\,
xyz^4+4y^2z^4+xyz^3+3y^2z^3+xyz^2+2y^2z^2+xyz+y^2z+xy
+
\cdots.
\endaligned\right.
\]
\[
\def\arraystretch{1.25}
\begin{array}{llll}
e_1
\,:=\,
-(z-1)\partial_x+y\partial_u
, &
e_2
\,:=\,
-2z\partial_x-(z-1)\partial_y+x\partial_u
,
\\
e_3
\,:=\,
-y\partial_x-y\partial_y-(z-1)\partial_z
, &
e_4
\,:=\,
x\partial_x+y\partial_y+2u\partial_u
.
\end{array}
\]

\[
\footnotesize
\def\arraystretch{1.25}
\begin{array}{c|cccc}
{} & e_1 & e_2 & e_3 & e_4
\\
\hline
e_1 & 
0 & 0 & e_1 & e_1
\\
e_2 &
0 & 0 & e_1 & e_2
\\
e_3 &
-e_1 & -e_1 & 0 & 0
\\
e_4 &
-e_1 & -e_2 & 0 & 0 
\end{array}
\]


\[
\text{\bf Model 2a3j}
\ \ \ \ \
\left\{
\aligned
u
&
\,=\,
xy+2x^3y-\tfrac{3}{2}xy^3+2x^3z+\tfrac{13}{2}y^3z+2y^2z^2+\tfrac{15}{4}x^4y+\tfrac{7}{4}x^3y^2+\tfrac{45}{4}x^2y^3
+
\\
&
\ \ \ \ \
+\tfrac{3}{2}xy^4+\tfrac{15}{2}x^4z+\tfrac{211}{4}y^4z+3x^3z^2+\tfrac{69}{2}y^3z^2+3y^2z^3+\tfrac{243}{8}x^5y+\tfrac{147}{2}x^4y^2
+
\\
&
\ \ \ \ \
+\tfrac{261}{4}x^3y^3+\tfrac{1521}{8}x^2y^4+84xy^5+\tfrac{165}{4}x^5z+\tfrac{4099}{8}y^5z+30x^4z^2+\tfrac{1903}{4}y^4z^2
+
\\
&
\ \ \ \ \
+4x^3z^3+107y^3z^3+4y^2z^4+x^3+\tfrac{3}{2}y^4+3x^5+\tfrac{63}{4}y^5+\tfrac{51}{8}x^6+\tfrac{1181}{8}y^6
+
\\
&
\ \ \ \ \
+xyz+y^2z+xyz^2+xyz^3+xyz^4+6x^2yz+\tfrac{15}{2}xy^2z+36x^3yz+\tfrac{141}{2}x^2y^2z
+
\\
&
\ \ \ \ \
+\tfrac{237}{4}xy^3z+18x^2yz^2+\tfrac{69}{2}xy^2z^2+\tfrac{885}{4}x^4yz+\tfrac{989}{2}x^3y^2z+\tfrac{1557}{2}x^2y^3z
+
\\
&
\ \ \ \ \
+\tfrac{4845}{8}xy^4z+192x^3yz^2+507x^2y^2z^2+\tfrac{1233}{2}xy^3z^2+36x^2yz^3+93xy^2z^3
+
\cdots
.
\endaligned\right.
\]
\[
\def\arraystretch{1.25}
\begin{array}{ll}
e_1
&
\,:=\,
(-\tfrac{3}{2}x+1-z+\tfrac{3}{2}u)\partial_x
-
(3x+3y+3u)\partial_y
-
(3x-\tfrac{3}{2}y+6z+\tfrac{21}{4}u)\partial_z
-
(-y+\tfrac{9}{2}u)\partial_u
,
\\
e_2
&
\,:=\,
(-\tfrac{13}{2}x-2z+\tfrac{57}{4}u)\partial_x
+
(1-11y-z+15u\tfrac{1}{4})\partial_y
+
(-\tfrac{15}{4}x-6y-15z+\tfrac{243}{8}u)\partial_z
+
\\
&
\ \ \ \ \
-
(-x+\tfrac{35}{2}u)\partial_u
, 
\\
e_3
&
\,:=\,
(-2x-y+\tfrac{7}{2}u)\partial_x
+
(-3y+\tfrac{3}{2}u)\partial_y
+
(1-\tfrac{9}{2}x-\tfrac{13}{2}y-z+\tfrac{33}{4}u)\partial_z
-
4u\partial_u
.
\end{array}
\]

\[
\footnotesize
\def\arraystretch{1.25}
\begin{array}{c|ccc}
{} & e_1 & e_2 & e_3
\\
\hline
e_1 & 
0 & -\tfrac{13}{2}e_1+3e_2-\tfrac{21}{4}e_3 & 
-e_1+\tfrac{3}{2}e_3
\\
e_2 &
\tfrac{13}{2}e_1-3e_2+\tfrac{21}{4}e_3
& 0 & e_1-2e_2+\tfrac{17}{2}e_3
\\
e_3 &
e_1-\tfrac{3}{2}e_3 & -e_1+2e_2-\tfrac{17}{2}e_3 & 0
\end{array}
\]


\[
\text{\bf Model 2a3m4k}
\ \ \ \ \
\left\{
\aligned
u
&
\,=\,
xy+x^2z+xyz+\tfrac{1}{4}y^2z+2x^2z^2+2xyz^2+\tfrac{1}{2}y^2z^2+4x^2z^3+4xyz^3
+
\\
&
\ \ \ \ \
+y^2z^3+8x^2z^4+8xyz^4+2y^2z^4+16x^2z^5+16xyz^5+4y^2z^5
+
\cdots
.
\endaligned\right.
\]
\[
\def\arraystretch{1.25}
\begin{array}{llll}
e_1
\,:=\,
-(-1+z)\partial_x-2z\partial_y+y\partial_u
, &
e_2
\,:=\,
-\tfrac{1}{2}z\partial_x-(-1+z)\partial_y+x\partial_u
,
\\
e_3
\,:=\,
-\tfrac{1}{4}y\partial_x-x\partial_y
-
(2z-1)\partial_z+u\partial_u, &
e_4
\,:=\,
x\partial_x+y\partial_y+2u\partial_u
.
\end{array}
\]

\[
\footnotesize
\def\arraystretch{1.25}
\begin{array}{c|cccc}
{} & e_1 & e_2 & e_3 & e_4
\\
\hline
e_1 & 
0 & 0 & e_1+e_2 & e_1
\\
e_2 &
0 & 0 & \tfrac{1}{4}e_1+e_2 & e_2
\\
e_3 &
-e_1-e_2 & -\tfrac{1}{4}e_1-e_2 & 0 & 0
\\
e_4 &
-e_1 & -e_2 & 0 & 0
\end{array}
\]

\SectionHead{\bf 2b Models}
{2b-models}

Here, we present the Branch~{\green{\bf 2b}} models.

\[
\text{\bf Model 2b3a4a}
\ \ \ \ \
\left\{
\aligned
u
&
\,=\,
x^2+y^2
.
\endaligned\right.
\]
\[
\def\arraystretch{1.25}
\begin{array}{llll}
e_1
\,:=\,
2x\partial_u+\partial_x
, &
e_2
\,:=\,
2y\partial_u+\partial_y
, &
& 
\\
e_3
\,:=\,
\partial_z
, &
e_4
\,:=\,
2u\partial_u+x\partial_x+y\partial_y
, &
e_5
\,:=\,
y\partial_x-x\partial_y
, &
e_6
\,:=\,
x\partial_z
,
\\
e_7
\,:=\,
y\partial_z
, &
e_8
\,:=\,
z\partial_z
, &
e_9
\,:=\,
u\partial_z
.
\end{array}
\]

\rotatebox[origin=c]{90}{
	\begin{tabular}{p{3cm}}
\end{tabular}}

\[
\footnotesize
\def\arraystretch{1.25}
\begin{array}{c|ccccccccc}
{} & e_1 & e_2 & e_3 & e_4 & e_5 & e_6 & e_7 & e_8 & e_9 
\\
\hline
e_1 & 
0 & 0 & 0 & e_1 & -e_2 & e_3 & 0 & 0 & 2e_6
\\
e_2 &
0 & 0 & 0 & e_2 & e_1 & 0 & e_3 & 0 & 2e_7
\\
e_3 &
0 & 0 & 0 & 0 & 0 & 0 & 0 & e_3 & 0
\\
e_4 &
-e_1 & -e_2 & 0 & 0 & 0 & e_6 & e_7 & 0 & 2e_9
\\
e_5 &
e_2 & -e_1 & 0 & 0 & 0 & e_7 & -e_6 & 0 & 0
\\
e_6 &
-e_3 & 0 & 0 & -e_6 & -e_7 & 0 & 0 & e_6 & 0
\\
e_7 &
0 & -e_3 & 0 & -e_7 & e_6 & 0 & 0 & e_7 & 0
\\
e_8 &
0 & 0 & -e_3 & 0 & 0 & -e_6 & -e_7 & 0 & -e_9
\\
e_9 &
-2e_6 & -2e_7 & 0 & -2e_9 & 0 & 0 & 0 & e_9 & 0
\end{array}
\]


\[
\text{\bf Model 2b3a4b}
\ \ \ \ \
\left\{
\aligned
u
&
=
x^2+y^2+\tfrac{1}{2}x^4+x^2y^2+\tfrac{1}{2}y^4+\tfrac{1}{2}x^6+\tfrac{3}{2}x^4y^2+\tfrac{3}{2}x^2y^4+\tfrac{1}{2}y^6
+
\cdots
.
\endaligned\right.
\]
\[
\def\arraystretch{1.25}
\begin{array}{llll}
e_1
\,:=\,
-(u-1)\partial_x+2x\partial_u
, &
e_2
\,:=\,
-(u-1)\partial_y+2y\partial_u
, &
& 
\\
e_3
\,:=\,
\partial_z
, &
e_4
\,:=\,
y\partial_x-x\partial_y
, &
e_5
\,:=\,
x\partial_z
, &
e_6
\,:=\,
y\partial_z
,
\\
e_7
\,:=\,
z\partial_z
, &
e_8
\,:=\,
u\partial_z
.
\end{array}
\]

\[
\footnotesize
\def\arraystretch{1.25}
\begin{array}{c|cccccccc}
{} & e_1 & e_2 & e_3 & e_4 & e_5 & e_6 & e_7 & e_8 
\\
\hline
e_1 & 
0 & 2e_4 & 0 & -e_2 & e_3-e_8 & 0 & 0 & 2e_5
\\
e_2 &
-2e_4 & 0 & 0 & e_1 & 0 & e_3-e_8 & 0 & 2e_6
\\
e_3 &
0 & 0 & 0 & 0 & 0 & 0 & e_3 & 0
\\
e_4 &
e_2 & -e_1 & 0 & 0 & e_6 & -e_5 & 0 & 0
\\
e_5 &
-e_3+e_8 & 0 & 0 & -e_6 & 0 & 0 & e_5 & 0
\\
e_6 &
0 & -e_3+e_8 & 0 & e_5 & 0 & 0 & e_6 & 0
\\
e_7 &
0 & 0 & -e_3 & 0 & -e_5 & -e_6 & 0 & -e_8
\\
e_8 &
-2e_5 & -2e_6 & 0 & 0 & 0 & 0 & e_8 & 0
\end{array}
\]


\[
\text{\bf Model 2b3a4c}
\ \ \ \ \
\left\{
\aligned
u
&
=
x^2+y^2-\tfrac{1}{2}x^4-x^2y^2-\tfrac{1}{2}y^4+\tfrac{1}{2}x^6+\tfrac{3}{2}x^4y^2+\tfrac{3}{2}x^2y^4+\tfrac{1}{2}y^6
+
\cdots
.
\endaligned\right.
\]
\[
\def\arraystretch{1.25}
\begin{array}{llll}
e_1
\,:=\,
(u+1)\partial_x+2x\partial_u
, &
e_2
\,:=\,
(u+1)\partial_y+2y\partial_u
, &
& 
\\
e_3
\,:=\,
\partial_z
, &
e_4
\,:=\,
y\partial_x-x\partial_y
, &
e_5
\,:=\,
x\partial_z
, &
e_6
\,:=\,
y\partial_z
,
\\
e_7
\,:=\,
z\partial_z
, &
e_8
\,:=\,
u\partial_z
.
\end{array}
\]

\[
\footnotesize
\def\arraystretch{1.25}
\begin{array}{c|cccccccc}
{} & e_1 & e_2 & e_3 & e_4 & e_5 & e_6 & e_7 & e_8 
\\
\hline
e_1 & 
0 & -2e_4 & 0 & -e_2 & e_3+e_8 & 0 & 0 & 2e_5
\\
e_2 &
2e_4 & 0 & 0 & e_1 & 0 & e_3+e_8 & 0 & 2e_6
\\
e_3 &
0 & 0 & 0 & 0 & 0 & 0 & e_3 & 0
\\
e_4 &
e_2 & -e_1 & 0 & 0 & e_6 & -e_5 & 0 & 0
\\
e_5 &
-e_3-e_8 & 0 & 0 & -e_6 & 0 & 0 & e_5 & 0
\\
e_6 &
0 & -e_3-e_8 & 0 & e_5 & 0 & 0 & e_6 & 0
\\
e_7 &
0 & 0 & -e_3 & 0 & -e_5 & -e_6 & 0 & -e_8
\\
e_8 &
-2e_5 & -2e_6 & 0 & 0 & 0 & 0 & e_8 & 0
\end{array}
\]

\[
\text{\bf Model 2b3b}
\ \ \ \ \
\left\{
\aligned
u
&
=
x^2+y^2+x^3+F_{0,4,0}y^4+F_{1,3,0}xy^3+F_{2,2,0}x^2y^2+F_{3,1,0}x^3y+F_{4,0,0}x^4
+
\\
&
\ \ \ \ \
+(F_{2,2,0}-\tfrac{4}{5}F_{2,2,0}F_{4,0,0}+\tfrac{8}{5}F_{4,0,0}^2-\tfrac{3}{5}F_{4,0,0}+\tfrac{1}{15}F_{3,1,0}F_{1,3,0}-\tfrac{1}{5}F_{3,1,0}^2)x^5
+
\\
&
\ \ \ \ \
+
(-\tfrac{9}{8}F_{3,1,0}+\tfrac{1}{4}F_{1,3,0}-\tfrac{1}{3}F_{4,0,0}F_{1,3,0}+3F_{4,0,0}F_{3,1,0}+\tfrac{1}{6}F_{2,2,0}F_{1,3,0}
+
\\
&
\ \ \ \ \
-\tfrac{3}{2}F_{2,2,0}F_{3,1,0})x^4y
+
(\tfrac{1}{3}F_{1,3,0}^2-\tfrac{4}{3}F_{3,1,0}F_{1,3,0}-\tfrac{4}{3}F_{2,2,0}^2+\tfrac{8}{3}F_{2,2,0}F_{4,0,0}
+
\\
&
\ \ \ \ \
+F_{3,1,0}^2-F_{2,2,0})x^3y^2
+
(-\tfrac{15}{4}F_{1,3,0}+4F_{4,0,0}F_{1,3,0}-\tfrac{7}{3}F_{2,2,0}F_{1,3,0}
+
\\
&
\ \ \ \ \
+F_{2,2,0}F_{3,1,0}+\tfrac{2}{3}F_{0,4,0}F_{1,3,0}-2F_{0,4,0}F_{3,1,0})x^2y^3
+
(-\tfrac{1}{3}F_{1,3,0}^2+F_{3,1,0}F_{1,3,0}
+
\\
&
\ \ \ \ \
-9F_{0,4,0}-4F_{0,4,0}F_{2,2,0}+8F_{0,4,0}F_{4,0,0})xy^4
+
(-\tfrac{22}{15}F_{0,4,0}F_{1,3,0}
+
\\
&
\ \ \ \ \
+\tfrac{2}{15}F_{2,2,0}F_{1,3,0}+\tfrac{2}{5}F_{0,4,0}F_{3,1,0})y^5
+
\cdots
.
\endaligned\right.
\]
\[
\def\arraystretch{1.25}
\begin{array}{ll}
e_1
&
\,:=\,
(2F_{2,2,0}x-4F_{4,0,0}x+3x+1+\tfrac{1}{3}yF_{1,3,0}-yF_{3,1,0}-uF_{2,2,0})\partial_x
-
(\tfrac{1}{3}xF_{1,3,0}-xF_{3,1,0}-\tfrac{9}{2}y
+
\\
&
\ \ \ \ \
-2yF_{2,2,0}+4yF_{4,0,0}+\tfrac{1}{2}uF_{1,3,0})\partial_y
+
(4F_{2,2,0}u-8F_{4,0,0}u+9u+2x)\partial_u
,
\\
e_2
&
\,:=\,
(3xF_{1,3,0}-xF_{3,1,0}+\tfrac{4}{3}yF_{0,4,0}-\tfrac{2}{3}yF_{2,2,0}-\tfrac{3}{2}uF_{1,3,0})\partial_x
-
(\tfrac{4}{3}xF_{0,4,0}-1-\tfrac{2}{3}F_{2,2,0}x
+
\\
&
\ \ \ \ \
-3yF_{1,3,0}+yF_{3,1,0}+2uF_{0,4,0})\partial_y
+
(6F_{1,3,0}u-2F_{3,1,0}u+2y)\partial_u
,
\\
e_3
&
\,:=\,
\partial_z
,
\,\,\,\,\,\,\,\,\,\,\,\,\,\,\,\,
e_4
\,:=\,
x\partial_z
,
\,\,\,\,\,\,\,\,\,\,\,\,\,\,\,\,
e_5
\,:=\,
y\partial_z
,
\,\,\,\,\,\,\,\,\,\,\,\,\,\,\,\,
e_6
\,:=\,
z\partial_z
,
\,\,\,\,\,\,\,\,\,\,\,\,\,\,\,\,
e_7
\,:=\,
u\partial_z.
\end{array}
\]

Gröbner basis generators of 
moduli space core algebraic variety in 
$\R^5 \ni (F_{0,4,0},F_{1,3,0},F_{2,2,0},F_{3,1,0},F_{4,0,0})$:
\[
\aligned
B_1 & := -32F_{0,4,0}F_{3,1,0}+4F_{1,3,0}F_{2,2,0}+136F_{1,3,0}F_{4,0,0}-20F_{2,2,0}F_{3,1,0}+24F_{3,1,0}F_{4,0,0}-174F_{1,3,0}
+
\\
&
\ \ \ \ \
-27F_{3,1,0},
\\
B_2 & := 8F_{0,4,0}F_{1,3,0}-24F_{0,4,0}F_{3,1,0}+88F_{1,3,0}F_{4,0,0}-16F_{2,2,0}F_{3,1,0}+24F_{3,1,0}F_{4,0,0}-111F_{1,3,0}
+
\\
&
\ \ \ \ \
-27F_{3,1,0},
\\
B_3 & := -4F_{0,4,0}F_{2,2,0}+8F_{0,4,0}F_{4,0,0}+F_{1,3,0}^2+5F_{1,3,0}F_{3,1,0}-2F_{2,2,0}^2+4F_{2,2,0}F_{4,0,0}-6F_{0,4,0}-3F_{2,2,0},
\\
B_4 & := 8F_{0,4,0}^2-12F_{0,4,0}F_{2,2,0}+8F_{0,4,0}F_{4,0,0}+32F_{1,3,0}F_{3,1,0}-12F_{2,2,0}^2+28F_{2,2,0}F_{4,0,0}+12F_{0,4,0}
+
\\
&
\ \ \ \ \
-21F_{2,2,0}.
\endaligned
\]
\[
\footnotesize
\def\arraystretch{1.25}
\begin{array}{c|ccccccc}
{} & e_1 & e_2 & e_3 & e_4 & e_5 & e_6 & e_7
\\
\hline
e_1 & 
0
&
\rotatebox[origin=c]{90}{
	\begin{tabular}{p{4cm}}$
\tfrac{8}{3}F_{1,3,0}e_1+(\tfrac{-9}{2}-\tfrac{4}{3}F_{2,2,0}+4F_{4,0,0}-\tfrac{4}{3}F_{0,4,0})e_2$\end{tabular}}
& 0 
&\rotatebox[origin=c]{90}{
	\begin{tabular}{p{4cm}}$
e_3+(2F_{2,2,0}-4F_{4,0,0}+3)e_4+(\tfrac{1}{3}F_{1,3,0}-F_{3,1,0})e_5-F_{2,2,0}e_7$\end{tabular}}
&\rotatebox[origin=c]{90}{
	\begin{tabular}{p{4cm}}$
(-\tfrac{1}{3}F_{1,3,0}+F_{3,1,0})e_4+(\tfrac{9}{2}+2F_{2,2,0}-4F_{4,0,0})e_5-\tfrac{1}{2}F_{1,3,0}e_7$\end{tabular}}
& 0 &\rotatebox[origin=c]{90}{
	\begin{tabular}{p{4cm}}$
2e_4+(4F_{2,2,0}-8F_{4,0,0}+9)e_7$\end{tabular}}
\\
e_2 &\rotatebox[origin=c]{90}{
	\begin{tabular}{p{4cm}}$
-\tfrac{8}{3}F_{1,3,0}e_1-(\tfrac{-9}{2}-\tfrac{4}{3}F_{2,2,0}+4F_{4,0,0}-\tfrac{4}{3}F_{0,4,0})e_2$\end{tabular}}
& 0 & 0 
&\rotatebox[origin=c]{90}{
	\begin{tabular}{p{3cm}}$
(3F_{1,3,0}-F_{3,1,0})e_4+(\tfrac{4}{3}F_{0,4,0}-\tfrac{2}{3}F_{2,2,0})e_5-\tfrac{3}{2}F_{1,3,0}e_7 $\end{tabular}}
&\rotatebox[origin=c]{90}{
	\begin{tabular}{p{4cm}}$
e_3+(-\tfrac{4}{3}F_{0,4,0}+\tfrac{2}{3}F_{2,2,0})e_4+(3F_{1,3,0}-F_{3,1,0})e_5-2F_{0,4,0}e_7$\end{tabular}}
& 0 &\rotatebox[origin=c]{90}{
	\begin{tabular}{p{4cm}}$
2e_5+(6F_{1,3,0}-2F_{3,1,0})e_7$\end{tabular}}
\\
e_3 &
0 & 0 & 0 & 0 & 0 & e_3 & 0
\\
e_4 &\rotatebox[origin=c]{90}{
	\begin{tabular}{p{4cm}}$
-e_3-(F_{2,2,0}-4F_{4,0,0}+3)e_4-(\tfrac{1}{3}F_{1,3,0}-F_{3,1,0})e_5+F_{2,2,0}e_7$\end{tabular}}
&\rotatebox[origin=c]{90}{
	\begin{tabular}{p{4cm}}$
-(3F_{1,3,0}-F_{3,1,0})e_4-(\tfrac{4}{3}F_{0,4,0}-\tfrac{2}{3}F_{2,2,0})e_5+\tfrac{3}{2}F_{1,3,0}e_7$\end{tabular}}
& 0 & 0 & 0 & e_4 & 0
\\
e_5 &\rotatebox[origin=c]{90}{
	\begin{tabular}{p{4cm}}$
-(-\tfrac{1}{3}F_{1,3,0}+F_{3,1,0})e_4-(\tfrac{9}{2}+2F_{2,2,0}-4F_{4,0,0})e_5+\tfrac{1}{2}F_{1,3,0}e_7$\end{tabular}}
&\rotatebox[origin=c]{90}{
	\begin{tabular}{p{4cm}}$
-e_3-(-\tfrac{4}{3}F_{0,4,0}+\tfrac{2}{3}F_{2,2,0})e_4-(3F_{1,3,0}-F_{3,1,0})e_5+2F_{0,4,0}e_7$\end{tabular}}
& 0 & 0 & 0 & e_5 & 0
\\
e_6 &
0 & 0 & -e_3 & -e_4 & -e_5 & 0 & -e_7
\\
e_7 &\rotatebox[origin=c]{90}{
	\begin{tabular}{p{4cm}}$
-2e_4-(4F_{2,2,0}-8F_{4,0,0}+9)e_7$\end{tabular}}
&\rotatebox[origin=c]{90}{
	\begin{tabular}{p{4cm}}$ 
-2e_5-(6F_{1,3,0}-2F_{3,1,0})e_7$\end{tabular}} & 0 & 0 & 0 & e_7 & 0
\end{array}
\]


\[
\text{\bf Model 2b3c}
\ \ \ \ \
\left\{
\aligned
u
&
=
x^2z^5+y^2z^5+x^2z^4+y^2z^4+x^2z^3+y^2z^3+x^2z^2+y^2z^2+x^2z+y^2z
+
\\
&
\ \ \ \ \
+x^2+y^2
+
\cdots
.
\endaligned\right.
\]
\[
\def\arraystretch{1.25}
\begin{array}{llll}
e_1
\,:=\,
-(-1+z)\partial_x+2x\partial_u
, &
e_2
\,:=\,
-(-1+z)\partial_y+2y\partial_u
, &
& 
\\
e_3
\,:=\,
-(-1+z)\partial_z+u\partial_u
, &
e_4
\,:=\,
2u\partial_u+x\partial_x+y\partial_y
, &
e_5
\,:=\,
-y\partial_x+x\partial_y
, &
e_6
\,:=\,
-\tfrac{1}{2}u\partial_x+x\partial_z
,
\\
e_7
\,:=\,
-\tfrac{1}{2}u\partial_y+y\partial_z
.
\end{array}
\]

\[
\footnotesize
\def\arraystretch{1.25}
\begin{array}{c|ccccccc}
{} & e_1 & e_2 & e_3 & e_4 & e_5 & e_6 & e_7
\\
\hline
e_1 & 
0 & 0 & e_1 & e_1 & e_2 & e_3 & -e_5
\\
e_2 &
0 & 0 & e_2 & e_2 & -e_1 & e_5 & e_3
\\
e_3 &
-e_1 & -e_2 & 0 & 0 & 0 & e_6 & e_7
\\
e_4 &
-e_1 & -e_2 & 0 & 0 & 0 & e_6 & e_7
\\
e_5 &
-e_2 & e_1 & 0 & 0 & 0 & -e_7 & e_6
\\
e_6 &
-e_3 & -e_5 & -e_6 & -e_6 & e_7 & 0 & 0
\\
e_7 &
e_5 & -e_3 & -e_7 & -e_7 & -e_6 & 0 & 0
\end{array}
\]

\[
\text{\bf Model 2b3d}
\ \ \ \ \
\left\{
\aligned
u
&
\,=\,
-8F_{1,3,0}F_{4,0,0}zyx^4+16F_{1,3,0}F_{4,0,0}zy^3x^2+zx^2+zy^2+F_{4,0,0}x^4
+
\\
&
\ \ \ \ \
+z^2x^2+z^2y^2+F_{0,5,0}y^5+z^3x^2+z^3y^2+2zx^3-\tfrac{13}{90}F_{1,3,0}^2x^5+\tfrac{1}{5}F_{3,2,0}x^5
+
\\
&
\ \ \ \ \
+\tfrac{8}{5}F_{4,0,0}^2x^5-\tfrac{3}{5}F_{4,0,0}x^5+3z^2x^3-\tfrac{11}{120}F_{1,3,0}^2x^6+\tfrac{3}{10}F_{3,2,0}x^6+\tfrac{3}{5}F_{4,0,0}x^6
+
\\
&
\ \ \ \ \
-\tfrac{14}{5}F_{4,0,0}^2x^6+\tfrac{16}{5}F_{4,0,0}^3x^6-\tfrac{5}{8}F_{1,3,0}^2y^6+4z^3x^3+z^4x^2+z^4y^2+F_{1,3,0}F_{3,2,0}y^5x
+
\\
&
\ \ \ \ \
-2F_{1,3,0}F_{4,0,0}yx^4+4F_{1,3,0}F_{4,0,0}y^3x^2+3F_{1,3,0}zy^3x+12F_{1,3,0}F_{4,0,0}yx^5
+
\\
&
\ \ \ \ \
+\tfrac{1}{15}F_{1,3,0}F_{3,2,0}yx^5+\tfrac{20}{3}F_{0,5,0}F_{4,0,0}yx^5-8F_{1,3,0}F_{4,0,0}^2yx^5-\tfrac{55}{18}F_{0,5,0}F_{1,3,0}y^2x^4
+
\\
&
\ \ \ \ \
+\tfrac{83}{18}F_{4,0,0}F_{1,3,0}^2y^2x^4+4F_{4,0,0}F_{3,2,0}y^2x^4-29F_{1,3,0}F_{4,0,0}y^3x^3-\tfrac{4}{3}F_{1,3,0}F_{3,2,0}y^3x^3
+
\\
&
\ \ \ \ \
+16F_{1,3,0}F_{4,0,0}^2y^3x^3-12F_{4,0,0}F_{1,3,0}^2y^4x^2+\tfrac{5}{2}F_{0,5,0}F_{1,3,0}y^4x^2+12F_{0,5,0}F_{4,0,0}y^5x
+
\\
&
\ \ \ \ \
+\tfrac{20}{3}F_{0,5,0}zyx^4+9F_{1,3,0}zyx^4+4F_{3,2,0}zy^2x^3-15F_{1,3,0}zy^3x^2-6F_{1,3,0}^2zy^4x
+
\\
&
\ \ \ \ \
+6F_{1,3,0}z^2y^3x+x^2+y^2+x^3-\tfrac{10}{9}F_{1,3,0}^3y^3x^3+3F_{4,0,0}zx^4+\tfrac{2}{3}F_{4,0,0}F_{1,3,0}^2y^6
+
\\
&
\ \ \ \ \
+4F_{0,5,0}zy^5+\tfrac{27}{2}F_{1,3,0}^2y^4x^2-\tfrac{3}{2}F_{1,3,0}^2y^4x+F_{1,3,0}y^3x-\tfrac{9}{4}F_{3,2,0}y^2x^4-\tfrac{26}{45}F_{1,3,0}^2zx^5
+
\\
&
\ \ \ \ \
-\tfrac{27}{8}F_{1,3,0}yx^5+\tfrac{4}{5}F_{3,2,0}zx^5-\tfrac{27}{2}F_{0,5,0}y^5x+\tfrac{105}{8}F_{1,3,0}y^3x^3+F_{3,2,0}y^2x^3
+
\\
&
\ \ \ \ \
-\tfrac{83}{16}F_{1,3,0}^2y^2x^4-\tfrac{37}{18}F_{0,5,0}F_{1,3,0}y^6+\tfrac{2}{5}F_{4,0,0}F_{3,2,0}x^6+\tfrac{13}{270}F_{1,3,0}^3yx^5+6F_{4,0,0}z^2x^4
+
\\
&
\ \ \ \ \
-\tfrac{15}{4}F_{1,3,0}y^3x^2-\tfrac{2}{5}F_{4,0,0}F_{1,3,0}^2x^6+\tfrac{9}{4}F_{1,3,0}yx^4+\tfrac{32}{5}F_{4,0,0}^2zx^5+\tfrac{5}{54}F_{0,5,0}F_{1,3,0}x^6
+
\\
&
\ \ \ \ \
+\tfrac{5}{3}F_{0,5,0}yx^4+\tfrac{17}{18}F_{1,3,0}^3y^5x-\tfrac{5}{2}F_{0,5,0}yx^5-\tfrac{12}{5}F_{4,0,0}zx^5
+
\cdots
.
\endaligned\right.
\]
\[
\def\arraystretch{1.25}
\aligned
e_1 & := 
-(4F_{4,0,0}x-3x-1-\tfrac{1}{3}F_{1,3,0}y+z+\tfrac{4}{9}uF_{1,3,0}^2+uF_{3,2,0})\partial_x+(-\tfrac{1}{3}F_{1,3,0}x+\tfrac{9}{2}y-4yF_{4,0,0}
+
\\
&
\ \ \ \ \
+\tfrac{7}{2}uF_{1,3,0}-\tfrac{10}{3}uF_{1,3,0}F_{4,0,0}+\tfrac{10}{3}uF_{0,5,0})\partial_y-(-\tfrac{8}{9}xF_{1,3,0}^2-2xF_{3,2,0}-\tfrac{20}{3}F_{1,3,0}yF_{4,0,0}
+
\\
&
\ \ \ \ \
+\tfrac{20}{3}yF_{0,5,0}+8F_{1,3,0}y-\tfrac{7}{6}uF_{1,3,0}^2)\partial_z-(8F_{4,0,0}u-9u-2x)\partial_u,
\\
e_2 & := (3F_{1,3,0}x+\tfrac{10}{3}uF_{0,5,0}-4uF_{1,3,0}F_{4,0,0}+\tfrac{15}{4}uF_{1,3,0})\partial_x+(3F_{1,3,0}y+F_{3,2,0}u-z+1)\partial_y
+
\\
&
\ \ \ \ \
-(\tfrac{20}{3}xF_{0,5,0}+\tfrac{21}{2}F_{1,3,0}x-8F_{1,3,0}xF_{4,0,0}+2yF_{3,2,0}+5uF_{0,5,0})\partial_z+(6F_{1,3,0}u+2y)\partial_u,
\\
e_3 & := -x\partial_x-y\partial_y-(-1+z)\partial_z-u\partial_u.
\endaligned
\]
Gröbner basis generators of 
moduli space core algebraic variety in 
$\R^4 \ni (F_{0,5,0},F_{1,3,0},F_{3,2,0},F_{4,0,0})$:
\[
\aligned
B_1 & := 60672F_{1,3,0}F_{4,0,0}^2-17920F_{0,5,0}F_{4,0,0}+12288F_{1,3,0}F_{3,2,0}-136512F_{1,3,0}F_{4,0,0}+20160F_{0,5,0}
+
\\
&
\ \ \ \ \
+74655F_{1,3,0},
\\
B_2 & := 632F_{1,3,0}^2F_{4,0,0}-640F_{0,5,0}F_{1,3,0}-711F_{1,3,0}^2-288F_{3,2,0}F_{4,0,0}+324F_{3,2,0},
\\
B_3 & := 43520F_{0,5,0}F_{1,3,0}F_{4,0,0}+12288F_{1,3,0}^2F_{3,2,0}+27648F_{3,2,0}F_{4,0,0}^2-48960F_{0,5,0}F_{1,3,0}-2133F_{1,3,0}^2
+
\\
&
\ \ \ \ \
-62208F_{3,2,0}F_{4,0,0}+34992F_{3,2,0},
\\
B_4 & := F_{0,5,0}(43520F_{0,5,0}F_{4,0,0}+16384F_{1,3,0}F_{3,2,0}-48960F_{0,5,0}-2133F_{1,3,0}),
\\
B_5 & := 158F_{1,3,0}^3-1080F_{0,5,0}F_{4,0,0}-72F_{1,3,0}F_{3,2,0}+1215F_{0,5,0}.
\endaligned
\]
\[
\footnotesize
\def\arraystretch{1.25}
\begin{array}{c|ccc}
{} & e_1 & e_2 & e_3
\\
\hline
e_1 & 
0 
&
\rotatebox[origin=c]{0}{
	\begin{tabular}{p{7cm}}
	$\tfrac{8}{3}F_{1,3,0}e_1+(\tfrac{-9}{2}+4F_{4,0,0})e_2+(\tfrac{4}{3}F_{1,3,0}F_{4,0,0}-\tfrac{5}{2}F_{1,3,0})e_3$
	\end{tabular}}
& 0
\\
e_2 &
\rotatebox[origin=c]{0}{
	\begin{tabular}{p{7cm}}
	$-\tfrac{8}{3}F_{1,3,0}e_1-(\tfrac{-9}{2}+4F_{4,0,0})e_2-(\tfrac{4}{3}F_{1,3,0}F_{4,0,0}-\tfrac{5}{2}F_{1,3,0})e_3$
	\end{tabular}}
 & 0 & 0
\\
e_3 &
0 & 0 & 0
\end{array}
\]


\[
\text{\bf Model 2b3e4a}
\ \ \ \ \
\left\{
\aligned
u
&
\,=\,
x^2+y^2
+
y^2z
+
y^2z^2
+
y^2z^3
+
y^2z^4
+
y^2z^5
+
\cdots
.
\endaligned\right.
\]
\[
\def\arraystretch{1.25}
\begin{array}{llll}
e_1
\,:=\,
2x\partial_u+\partial_x
, &
e_2
\,:=\,
-(-1+z)\partial_y+2y\partial_u
,
\\
e_3
\,:=\,
-\tfrac{y}{2}\partial_y-(-1+z)\partial_z
, &
e_4
\,:=\,
2u\partial_u+x\partial_x+y\partial_y
.
\end{array}
\]

\[
\footnotesize
\def\arraystretch{1.25}
\begin{array}{c|cccc}
{} & e_1 & e_2 & e_3 & e_4
\\
\hline
e_1 & 
0 & 0 & 0 & e_1
\\
e_2 &
0 & 0 & \tfrac{1}{2}e_2 & e_2 
\\
e_3 &
0 & -\tfrac{1}{2}e_2 & 0 & 0
\\
e_4 &
-e_1 & -e_2 & 0 & 0
\end{array}
\]


\[
\text{\bf Model 2b3e4b}
\ \ \ \ \
\left\{
\aligned
u
&
\,=\,
x^2+y^2+zy^2+y^4+z^2y^2+F_{0,5,0}y^5+3zy^4+z^3y^2+\tfrac{5}{4}F_{0,5,0}^2y^6
+
\\
&
\ \ \ \ \
+2y^6+4F_{0,5,0}zy^5+6z^2y^4+z^4y^2+\tfrac{25}{14}F_{0,5,0}^3y^7+\tfrac{36}{7}F_{0,5,0}y^7
+
\\
&
\ \ \ \ \
+\tfrac{25}{4}F_{0,5,0}^2zy^6+10zy^6+10F_{0,5,0}z^2y^5+10z^3y^4+z^5y^2
+
\cdots
,
\endaligned\right.
\]
for any value for $F_{0,5,0}$.
\[
\def\arraystretch{1.25}
\aligned
e_1 & := \partial_x+2x\partial_u,
\\
e_2 & := -(\tfrac{5}{2}xF_{0,5,0})\partial_x-(-1+\tfrac{5}{2}yF_{0,5,0}+z)\partial_y-4y\partial_z-(5F_{0,5,0}u-2y)\partial_u,
\\
e_3 & := -\tfrac{1}{2}x\partial_x-y\partial_y-(z-1)\partial_z-u\partial_u,
\endaligned
\]

\[
\footnotesize
\def\arraystretch{1.25}
\begin{array}{c|ccc}
{} & e_1 & e_2 & e_3
\\
\hline
e_1 & 
0 & -\tfrac{5}{2}F_{0,5,0}e_1 & -\tfrac{1}{2}e_1
\\
e_2 &
\tfrac{5}{2}F_{0,5,0}e_1 & 0 & 0
\\
e_3 &
\tfrac{1}{2}e_1 & 0 & 0
\end{array}
\]


\[
\text{\bf Model 2b3e4c}
\ \ \ \ \
\left\{
\aligned
u
&
\,=\,
x^2+y^2+zy^2-y^4+z^2y^2+F_{0,5,0}y^5-3zy^4+z^3y^2-\tfrac{5}{4}F_{0,5,0}^2y^6
+
\\
&
\ \ \ \ \
+2y^6+4F_{0,5,0}zy^5-6z^2y^4+z^4y^2+\tfrac{25}{14}F_{0,5,0}^3y^7-\tfrac{36}{7}F_{0,5,0}y^7
+
\\
&
\ \ \ \ \
-\tfrac{25}{4}F_{0,5,0}^2zy^6+10zy^6+10F_{0,5,0}z^2y^5-10z^3y^4+z^5y^2
+
\cdots
,
\endaligned\right.
\]
for any value for $F_{0,5,0}$.
\[
\def\arraystretch{1.25}
\aligned
e_1 & := \partial_x+2x\partial_u,
\\
e_2 & := (\tfrac{5}{2}xF_{0,5,0})\partial_x+(1+\tfrac{5}{2}yF_{0,5,0}-z)\partial_y+4y\partial_z+(5F_{0,5,0}u+2y)\partial_u,
\\
e_3 & := -\tfrac{x}{2}\partial_x-y\partial_y
-(z-1)\partial_z-u\partial_u,
\endaligned
\]

\[
\footnotesize
\def\arraystretch{1.25}
\begin{array}{c|ccc}
{} & e_1 & e_2 & e_3
\\
\hline
e_1 & 
0 & \tfrac{5}{2}F_{0,5,0}e_1 & -\tfrac{1}{2}e_1
\\
e_2 &
-\tfrac{5}{2}F_{0,5,0}e_1 & 0 & 0
\\
e_3 &
\tfrac{1}{2}e_1 & 0 & 0
\end{array}
\]


\[
\text{\bf Model 2b3e4d}
\ \ \ \ \
\left\{
\aligned
u
&
\,=\,
x^2+y^2+zy^2+F_{4,0,0}x^4+\tfrac{1}{4}y^4+zy^2x+z^2y^2+\tfrac{2}{5}F_{4,0,0}x^5+\tfrac{8}{5}F_{4,0,0}^2x^5
+
\\
&
\ \ \ \ \
-2F_{4,0,0}y^2x^3+\tfrac{1}{2}y^4x+zy^2x^2+\tfrac{1}{2}zy^4+2z^2y^2x+z^3y^2+\tfrac{18}{5}F_{4,0,0}^2x^6
+
\\
&
\ \ \ \ \
+\tfrac{16}{5}F_{4,0,0}^3x^6+\tfrac{1}{5}F_{4,0,0}x^6-4F_{4,0,0}^2y^2x^4-3F_{4,0,0}y^2x^4+\tfrac{3}{4}y^4x^2+\tfrac{3}{2}F_{4,0,0}y^4x^2
+
\\
&
\ \ \ \ \
+zy^2x^3-2F_{4,0,0}zy^2x^3+\tfrac{3}{2}zy^4x+3z^2y^2x^2+\tfrac{3}{4}z^2y^4+3z^3y^2x+z^4y^2
+
\\
&
\ \ \ \ \
+\tfrac{116}{35}F_{4,0,0}^2x^7+\tfrac{464}{35}F_{4,0,0}^3x^7+\tfrac{256}{35}F_{4,0,0}^4x^7+\tfrac{4}{35}F_{4,0,0}x^7-\tfrac{48}{5}F_{4,0,0}^3y^2x^5
+
\\
&
\ \ \ \ \
-\tfrac{74}{5}F_{4,0,0}^2y^2x^5-\tfrac{18}{5}F_{4,0,0}y^2x^5+y^4x^3+4F_{4,0,0}^2y^4x^3+4F_{4,0,0}y^4x^3
+
\\
&
\ \ \ \ \
-\tfrac{1}{2}F_{4,0,0}y^6x+zy^2x^4-4F_{4,0,0}^2zy^2x^4-5F_{4,0,0}zy^2x^4+3zy^4x^2
+
\\
&
\ \ \ \ \
+3F_{4,0,0}zy^4x^2+4z^2y^2x^3-2F_{4,0,0}z^2y^2x^3+3z^2y^4x+6z^3y^2x^2+z^3y^4
+
\\
&
\ \ \ \ \
+4z^4y^2x+z^5y^2
+
\cdots
,
\endaligned\right.
\]
for any value of $F_{4,0,0}$.
\[
\def\arraystretch{1.25}
\aligned
e_1 & := -(2F_{4,0,0}u+4F_{4,0,0}x+x-1)\partial_x-(4F_{4,0,0}y+y)\partial_y+(2F_{4,0,0}u+4F_{4,0,0}x-z)\partial_z
+
\\
&
\ \ \ \ \
-(8F_{4,0,0}u+2u-2x)\partial_u,
\\
e_2 & := y\partial_x-(x-1+z)\partial_y-y\partial_z+2y\partial_u,
\\
e_3 & := -\tfrac{y}{2}\partial_y-(x-1+z)\partial_z,
\endaligned
\]

\[
\footnotesize
\def\arraystretch{1.25}
\begin{array}{c|ccc}
{} & e_1 & e_2 & e_3
\\
\hline
e_1 &
0 & 4F_{4,0,0}e_2 & 0
\\
e_2 &
-4F_{4,0,0}e_2 & 0 & \tfrac{1}{2}e_2
\\
e_3 &
0 & -\tfrac{1}{2}e_2 & 0
\end{array}
\]


\[
\text{\bf Model 2b3e4e}
\ \ \ \ \
\left\{
\aligned
u
&
\,=\,
x^2+y^2+zy^2+\tfrac{1}{4}y^4+zy^3+z^2y^2+\tfrac{1}{2}y^5+\tfrac{9}{4}zy^4+3z^2y^3+z^3y^2+\tfrac{9}{8}y^6
+
\\
&
\ \ \ \ \
+\tfrac{23}{4}zy^5+\tfrac{39}{4}z^2y^4+6z^3y^3+z^4y^2+\tfrac{11}{4}y^7+\tfrac{127}{8}zy^6+\tfrac{131}{4}z^2y^5+\tfrac{115}{4}z^3y^4
+
\\
&
\ \ \ \ \
+10z^4y^3+z^5y^2
+
\cdots
.
\endaligned\right.
\]
\[
\def\arraystretch{1.25}
\begin{array}{ll}
e_1 & := \partial_x+2x\partial_u,
\\
e_2 & := -3x\partial_x-(3y+z-1)\partial_y-(y+3z)\partial_z-(6u-2y)\partial_u
\\
e_3 & := -(\tfrac{3}{2}x)\partial_x-2y\partial_y-(y+z-1)\partial_z-3u\partial_u.
\end{array}
\]

\[
\footnotesize
\def\arraystretch{1.25}
\begin{array}{c|ccc}
{} & e_1 & e_2 & e_3
\\
\hline
e_1 & 
0 & -3e_1 & -\tfrac{3}{2}e_1
\\
e_2 &
3e_1 & 0 & -e_2+2e_3
\\
e_3 &
\tfrac{3}{2}e_1 & e_2-2e_3 & 0
\end{array}
\]



\vfill\end{document}